\RequirePackage{rotating} 
\documentclass[opre,nonblindrev]{informs3} 

\OneAndAHalfSpacedXI

\usepackage{natbib}
 \bibpunct[, ]{(}{)}{,}{a}{}{,}%
 \def\bibfont{\footnotesize}%
 \def\bibsep{\smallskipamount}%
 \def\newblock{\ }%
\usepackage[utf8]{inputenc}
\usepackage[english]{babel}
\usepackage[margin=0.85in]{geometry}
\usepackage{bbding}
\usepackage[shortlabels]{enumitem}
\usepackage{{appendix}}
\usepackage[super]{nth}
\usepackage{bm, bbm}
\usepackage{amsmath, amssymb}
\usepackage{mathtools}
\usepackage{subcaption, placeins}
\usepackage{graphicx}
\usepackage[hidelinks]{hyperref}
\usepackage{pgfplots}

\usepackage{algorithm}
\usepackage{algorithmic}
\usepackage{subcaption}
\usepackage{booktabs}
\usepackage{multirow}
\usepackage{placeins}
\usepackage{booktabs}
\usepackage{accents}
\usepackage{graphics}

\usepackage{tikz}
\usepackage{placeins}
\usepackage{float}
\usepackage{graphicx}
\usepackage{accents}
\usepackage{rotating}
\usepackage{graphics}
\usepackage{todonotes}

\usepackage{color}
\definecolor{blue}{rgb}{0,0,1}

\newtheorem{appendixlemma}{Lemma}[section]
\setlist[itemize]{leftmargin=*, nosep}

\newcommand{\argminE}{\mathop{\mathrm{argmin}}}

\usepackage{etoolbox}
\newcommand{\zerodisplayskips}{%
  \setlength{\abovedisplayskip}{7pt}%
  \setlength{\belowdisplayskip}{5pt}%
  \setlength{\abovedisplayshortskip}{0pt}%
  \setlength{\belowdisplayshortskip}{0pt}}
\appto{\normalsize}{\zerodisplayskips}
\appto{\small}{\zerodisplayskips}
\appto{\footnotesize}{\zerodisplayskips}
\usepackage{xfrac}

\TheoremsNumberedThrough     
\ECRepeatTheorems

\EquationsNumberedThrough    

\begin{document}
\RUNAUTHOR{Bertsimas, Cory-Wright, and Pauphilet}
\RUNTITLE{Mixed-Projection Conic Optimization}

\TITLE{Mixed-Projection Conic Optimization: \\
A New Paradigm for Modeling Rank Constraints}

\ARTICLEAUTHORS{%
\AUTHOR{Dimitris Bertsimas}
\AFF{Sloan School of Management and Operations Research Center, Massachusetts Institute of Technology, Cambridge, MA, USA.\\ ORCID: \href{https://orcid.org/0000-0002-1985-1003}{$0000$-$0002$-$1985$-$1003$}\\ \EMAIL{dbertsim@mit.edu}} 
\AUTHOR{Ryan Cory-Wright}
\AFF{Operations Research Center, Massachusetts Institute of Technology, Cambridge, MA, USA, \\ORCID: \href{https://orcid.org/0000-0002-4485-0619}{$0000$-$0002$-$4485$-$0619$}\\\EMAIL{ryancw@mit.edu}}
\AUTHOR{Jean Pauphilet}
\AFF{London Business School, London, UK, \\ORCID: \href{https://orcid.org/0000-0001-6352-0984}{$0000$-$0001$-$6352$-$0984$}\\\EMAIL{jpauphilet@london.edu}}
}

\ABSTRACT{
We propose a framework for modeling and solving low-rank optimization problems to certifiable optimality. We introduce symmetric projection matrices that satisfy $\bm{Y}^2=\bm{Y}$, the matrix analog of binary variables that satisfy $z^2=z$, to model rank constraints. By leveraging regularization and strong duality, we prove that this modeling paradigm yields tractable convex optimization problems over the non-convex set of orthogonal projection matrices. Furthermore, we design outer-approximation algorithms to solve low-rank problems to certifiable optimality, compute lower bounds via their semidefinite relaxations, and provide near optimal solutions through rounding and local search techniques. We implement these numerical ingredients and, for the first time, solve low-rank optimization problems to certifiable optimality. Our algorithms also supply certifiably near-optimal solutions for larger problem sizes and outperform existing heuristics, {by deriving an alternative to the popular nuclear norm relaxation which generalizes the perspective relaxation from vectors to matrices.} {\color{black} Using currently available spatial branch-and-bound codes, not tailored to projection matrices, we can scale our exact (resp. near-exact) algorithms to matrices with up to $30$ (resp. $600$) rows/columns.} {All in all, our framework, which we name \emph{Mixed-Projection Conic Optimization}, solves low-rank problems to certifiable optimality in a tractable and unified fashion.}
}
\maketitle
\vspace{-10mm}
\SUBJECTCLASSname{programming: non-linear: quadratic, algorithms, applications}\\
\AREAOFREVIEWname{Optimization}
\section{Introduction}
Many central problems in optimization, machine learning, and control theory are equivalent to optimizing a low-rank matrix over a convex set. For instance, low-rank constraints successfully model notions of minimal complexity, low dimensionality, or orthogonality in a system. However, while rank constraints offer unparalleled modeling flexibility, no generic code currently solves these problems to certifiable optimality at even moderate sizes. This state of affairs has led influential works on low-rank optimization \citep[][]{candes2010matrix, recht2010guaranteed} to characterize low-rank optimization as intractable and advocate convex relaxations or heuristics which do not enjoy assumption-free optimality guarantees.

The manner in which rank constrained optimization is regarded today is reminiscent of how mixed-integer conic optimization (MICO), which can model NP-complete problems, was originally considered. After decades of research effort, however, algorithms and software for MICO are now widely available \citep[see, e.g.,][]{bonami2008algorithmic, coey2020outer} and solve large instances of disparate non-convex problems such as best subset selection \citep{bertsimas2017sparse} or {\color{black}sparse portfolio selection \citep{frangioni2007sdp,zheng2014improving, bertsimas2018scalable}} to certifiable optimality. Unfortunately, rank constraints cannot be represented using {\color{black}mixed-integer convex optimization} \citep[Lemma 4.1]{lubin2017regularity} and do not benefit from these advances.

In this work, we characterize the complexity of rank constrained optimization and propose a new, more general framework, which we term \textit{Mixed-Projection Conic Optimization} (MPCO). Our proposal generalizes MICO, by replacing binary variables $z$ which satisfy $z^2=z$ with symmetric orthogonal projection matrices $\bm{Y}$ which satisfy $\bm{Y}^2=\bm{Y}$, and offers the following advantages over existing state-of-the-art methods: First, it supplies certificates of (near) optimality for low-rank problems. Second, it demonstrates that some of the best ideas in MICO, such as decomposition methods, cutting-planes, relaxations, and random rounding schemes, admit straightforward extensions to MPCO. Finally, we implement a near-optimal rounding strategy and a globally optimal cutting-plane algorithm that improve upon the state-of-the-art for matrix completion and sensor location problems. We hope that MPCO gives rise to exciting new challenges for the optimization community to tackle.

\subsection{Scope of the Framework}\label{sec:examples}
Formally, we consider the problem:
\begin{equation}\label{rankminproblem}
\begin{aligned}
   \text{} \ \min_{\bm{X} \in \mathbb{R}^{n \times m}} \quad & \langle \bm{C}, \bm{X}\rangle+\lambda \cdot \mathrm{Rank}(
\bm{X}) \quad \text{s.t.} \quad \bm{A}\bm{X}=\bm{B},\  \mathrm{Rank}(\bm{X}) \leq k,\  \bm{X} \in \mathcal{K},
\end{aligned}
\end{equation}
where $\lambda$ (resp. $k$) prices (bounds) the rank of $\bm{X}$, $(\bm{A}, \bm{B}) \in \mathbb{R}^{\ell \times n} \times \mathbb{R}^{\ell \times m}$ defines an affine subspace, and $\mathcal{K}$ is a proper cone in the sense of \citet{boyd2004convex}, i.e., closed, convex, solid and pointed. Observe that Problem \eqref{rankminproblem} offers significant modeling flexibility, as it allows arbitrary conic constraints on $\bm{X}$. As a result, linear, convex quadratic, semidefinite, exponential, and power constraints and objectives can be captured by letting $\mathcal{K}$ be an appropriate product of the non-negative orthant and the second order, semidefinite, exponential, and power cones.

We now {\color{black}introduce our notation and }present some central problems from the optimization and machine learning literature which admit low-rank formulations and fall within our framework.

{\color{black}
\subsubsection*{Notation: }
We let nonbold face characters such as $b$ denote scalars, lowercase bold-faced characters such as $\bm{x}$ denote vectors, uppercase bold-faced characters such as $\bm{X}$ denote matrices, and calligraphic uppercase characters such as $\mathcal{Z}$ denote sets. We let $[n]$ denote the set of running indices $\{1, ..., n\}$. We let $\mathbf{e}$ denote a vector of all $1$'s, $\bm{0}$ denote a vector of all $0$'s, and $\mathbb{I}$ denote the identity matrix.

We also use an assortment of matrix operators. We let $\sigma_i(\bm{X})$ denote the $i$th largest singular value of a matrix $\bm{X}$, $\langle \cdot,\cdot \rangle$ denote the Euclidean inner product between two vectors or matrices of the same dimension, $\bm{X}^\dag$ denote the Moore-Penrose pseudoinverse of a matrix $\bm{X}$, $\Vert \cdot \Vert_F$ denote the Frobenius norm of a matrix, $\Vert \cdot \Vert_\sigma$ denote the spectral norm of a matrix, and $\Vert \cdot \Vert_*$ denote the nuclear norm of a matrix; see \citet{johnson1985matrix} for a general theory of matrix operators.

Finally, we use a wide variety of convex cones. We let $S^n$ denote the $n \times n$ cone of symmetric matrices, and $S^n_+$ denote the $n \times n$ positive semidefinite cone.}

\subsubsection{Low-Rank Matrix Completion} Given a sub-sample $(A_{i,j}: (i,j) \in \mathcal{I} \subseteq [n] \times {\color{black}[m]})$ of a matrix ${\bm{A} \in \mathbb{R}^{n \times {\color{black}m}}}$, the matrix completion problem is to recover the entire matrix, by assuming $\bm{A}$ is low rank and seeking a rank-$k$ matrix $\bm{X}$ which approximately fits the observed values. This problem arises in recommender system applications at Netflix and Amazon and admits the formulation:
\begin{align}\label{lowrankmatrixorig}
    \min_{\bm{X} \in \mathbb{R}^{n \times {\color{black}m}}} \quad & \frac{1}{2}\sum_{(i,j) \in \mathcal{I}} \left(X_{i,j}-A_{i,j} \right)^2 \quad \text{s.t.} \quad \mathrm{Rank}(\bm{X}) \leq k.
\end{align}

Since there are $(n+{\color{black}m})k$ degrees of freedom in a singular value decomposition of a rank-$k$ matrix $\bm{X} \in \mathbb{R}^{n \times {\color{black}m}}$, Problem \eqref{lowrankmatrixorig} is not well-defined unless $\vert \mathcal{I}\vert \geq (n+{\color{black}m})k$. 

\subsubsection{Minimum Dimension Euclidean Distance Embedding}
Given a set of pairwise distances $d_{i,j}$, the Euclidean Distance Embedding (EDM) problem is to determine the lowest dimensional space which the distances can be embedded in, such that the distances 
correspond to Euclidean distances. 
{\color{black}As discussed by} \citet{blekherman2012semidefinite} Theorem 2.49, a set of distances $d_{i,j}$ can be embedded in a Euclidean space of dimension $k$ if and only if there exists some Gram matrix $\bm{G} \succeq \bm{0}$ of rank $k$ such that $d_{i,j}^2 = G_{i,i} + G_{j,j} - 2 G_{i,j}$, on all pairs $(i,j)$ where $d_{i,j}$ is supplied.
Denoting $D_{i,j}=d_{i,j}^2$, we write these constraints in matrix form,
$\bm{D}=\mathrm{Diag}(\bm{G})\bm{e}^\top +\bm{e}\mathrm{Diag}(\bm{G})^\top -2\bm{G}$, where the equality is implicitly imposed only for pairs $(i,j)$ where $d_{i,j}$ is supplied. This is equivalent to:
\begin{equation}\label{prob:edmorig}
    \begin{aligned}
    \min_{\bm{G} \in S_+^n}  \ \mathrm{Rank}(\bm{G})\quad \text{s.t.} \quad \mathrm{Diag}(\bm{G})\bm{e}^\top +\bm{e}\mathrm{Diag}(\bm{G})^\top -2\bm{G}= \bm{D}.
\end{aligned}
\end{equation}

\subsubsection{Quadratically Constrained Quadratic Optimization}\label{sssec:qcqop}

A quadratically constrained quadratic {\color{black} optimization} problem ({\color{black}QCQO}) seeks an $\bm{x} \in \mathbb{R}^n$ which solves:
\begin{equation}\label{prob:qcqop}
\begin{aligned}
    \min_{\bm{x}\in \mathbb{R}^n} \ \bm{x}^\top \bm{Q}_0 \bm{x}+\bm{q}_0^\top \bm{x}\quad \text{s.t.} \quad \bm{x}^\top \bm{Q}_i \bm{x} + \bm{q}_i^\top \bm{x} \leq r_i \quad \forall i \in [m],
\end{aligned}
\end{equation}
where $\bm{Q}_0$, $\bm{Q}_i$ ,$\bm{q}_0$ $\bm{q}_i$, $r_i$ are given problem data. We assume that $\bm{Q}_0, \bm{Q}_i$ are symmetric matrices, but do not assume that they are positive semidefinite. Therefore, this problem is non-convex, and encompasses binary quadratic optimization \citep[]{goemans1995improved} and alternating current optimal power flow problems \citep{lavaei2011zero}. The fundamental difficulty in Problem \eqref{prob:qcqop} is the potential non-convexity of the outer product $\bm{x}\bm{x}^\top$. However, we can isolate this non-convexity by introducing a rank-one matrix $\bm{X}$ to model the outer product $\bm{x}\bm{x}^\top$. This leads to the following reformulation:
\begin{equation}\label{prob:qcqop_withrank}
\begin{aligned}
    \min_{\bm{x}\in \mathbb{R}^n, \bm{X} \in S^n} \ \langle \bm{Q}_0, \bm{X}\rangle+\langle\bm{q}_0, \bm{x}\rangle\quad \text{s.t.} \quad \langle \bm{Q}_i, \bm{X}\rangle+ \langle\bm{q}_i, \bm{x}\rangle \leq r_i \quad \forall i \in [m], \ \mathrm{Rank}\begin{pmatrix} 1 & \bm{x}^\top \\ \bm{x} & \bm{X}\end{pmatrix}=1.
\end{aligned}
\end{equation}

We have established that {\color{black}QCQO}s are rank constrained problems. Notably however, the converse is also true: rank constrained problems with linear, second-order cone, or semidefinite constraints are {\color{black}QCQO}s. Indeed, the constraint $\mathrm{Rank}(\bm{X}) \leq k$ is equivalent to requiring that $\bm{X}=\bm{U}\bm{V}^\top: \bm{U} \in \mathbb{R}^{n \times k}, \bm{V} \in \mathbb{R}^{m \times k}$, i.e., imposing $m \times n$ non-convex quadratic equalities. As modern solvers such as \verb|Gurobi| can now solve non-convex {\color{black}QCQO}s to global optimality, this {\color{black}QCQO} formulation can be used to solve low-rank problems, although it is not particularly scalable; we expand on this point in Section \ref{ssec:qcqopgurobi}.

\subsection{Background and Literature Review}

Our work arises at the intersection of three complementary areas of the low-rank optimization literature: (a) global optimization algorithms for non-convex quadratically constrained problems, (b) the interplay of convex relaxations and their dual side, randomized rounding methods, and (c) heuristics which provide high-quality solutions to non-convex problems in an efficient fashion.

\subsubsection{Global Optimization Techniques}
\paragraph{Branch-and-bound:} A broad class of global optimization algorithms have been proposed for {\color{black}QCQO}s, since \cite{mccormick1976computability} observed that convex envelopes of non-convex regions supply globally valid lower bounds. This gives rise to a numerical strategy where one recursively partitions the {\color{black}QCQO}'s feasible region into subregions, constructs convex envelopes for each subregion and uses these envelopes to construct iteratively improving lower bounds. This approach is known as spatial branch-and-bound; see {\color{black}\cite{lee2014optimal} for a scheme which decomposes a matrix into a sparse matrix plus a low-rank matrix,} \cite{kocuk2016strong} for a modern implementation in alternating current optimal power flow{\color{black}, and \cite{bertsimas2017certifiably} for an exact branch-and-bound approach to low-rank factor analysis.}

\paragraph{Branch-and-cut:} In a complementary direction, several branch-and-cut methods \citep{audet2000branch, linderoth2005simplicial} have been proposed for solving non-convex {\color{black}QCQO}s, by borrowing decomposition schemes from the mixed-integer nonlinear optimization literature \citep[][]{duran1986outer}. While often efficient in practice, a common theme in these methods is that the more efficient decomposition schemes used for MINLOs cannot be applied out-of-the-box, because they may fail to converge to a globally optimal solution \citep[see][for a counterexample]{grossmann2002review}. As a result, non-convex problems need to be preprocessed in an expensive fashion. This preprocessing step has inhibited the use of global optimization methods for low-rank problems; indeed, we are not aware of any works which apply branch-and-cut techniques to solve low-rank problems to certifiable optimality.

{\color{black}\paragraph{Complementarity:} In an opposite direction, several authors have proposed applying general nonlinear optimization techniques to address low-rank problems, since \citet{ding2014first} observed that a low-rank constraint is equivalent to a complementarity constraint over the positive semidefinite cone, and thus can be addressed by general techniques for mathematical programs with equilibrium constraints \citep[see][for a general theory]{luo1996mathematical}. Among others, \citet{bai2016conic} invoked the complementarity observation to design a completely positive reformulation of low-rank SDOs, and \citet{bi2020multi} developed a multi-stage convex relaxation of the complementarity constraint. We are not aware of any works which use these ideas to solve low-rank problems exactly where say $n \geq 10$.

}

\paragraph{Algebraic:} {By taking an algebraic view of rank constraints, several algebraic geometry techniques have been proposed for addressing low-rank SDOs. Among others, \cite{d2003semidefinite} proposed reformulating low-rank constraints as systems of polynomial equations which can be addressed via the sum-of-squares hierarchy \citep[][]{lasserre2001global}.
More recently, \citet{naldi2018solving} proposed a semi-algebraic reformulation of rank-constrained SDOs, which can be optimized over via Gr{\"o}bner basis computation \citep[][]{cox2013ideals}. Unfortunately, {\color{black}algebraic approaches do not scale well in practice.}
Indeed, as observed by \citet{recht2010guaranteed}, it seems unlikely that algebraic approaches can solve low-rank SDOs when $n>10$.

}

\subsubsection{Convex Relaxations and Random Rounding Methods for Low-Rank Problems}
\paragraph{Convex relaxations:} A number of authors have studied convex relaxations of low-rank problems, since \citet{fazel2002matrix} 
observed that the nuclear norm of a matrix is the convex envelope of a rank constraint on the set of matrices with spectral norm at most $M$, i.e.,
\begin{align}
    \mathrm{Conv}\left(\left\{\bm{X}
\in \mathbb{R}^{n \times m}:\ \Vert \bm{X}\Vert_\sigma \leq M, \mathrm{Rank}(\bm{X}) \leq k\right\}\right)=\left\{\bm{X}
\in \mathbb{R}^{n \times m}:\ \Vert \bm{X}\Vert_\sigma \leq M, \Vert \bm{X}\Vert_* \leq k M\right\}.
\end{align}
Because the epigraph of a nuclear norm is semidefinite representable, this gives rise to semidefinite relaxations of low-rank problems which can be computed in polynomial time. 

\paragraph{Rounding methods:} A complementary line of work aims to supply certifiably near-optimal solutions to low-rank problems, by rounding their semidefinite relaxations. Initiated by \citet{goemans1995improved} in the context of binary quadratic optimization, who established that randomly rounding an SDO relaxation supplies a $0.878$-approximation, it has evolved into a successful framework for solving rank-one optimization problems; see \citet[]{nemirovski1999maximization} for a unified approach in the rank-one case. However, this line of work has a key drawback. Namely, existing rounding methods do not address rank-$k$ problems such as matrix completion, due to the analytic difficulty of constructing a rounding mechanism which preserves both feasibility and near-optimality in the rank-$k$ case.


\subsubsection{Heuristic Methods} \label{ssec:heurmethods}
Due to the computational difficulty of solving Problem \eqref{rankminproblem} to certifiable optimality, and the analytic difficulty of deriving a high-quality randomized rounding procedure, a variety of heuristic methods have been proposed for solving Problem \eqref{rankminproblem}, originating with methods for solving low-rank linear matrix inequalities in the optimal control literature \citep{boyd1994linear}.

Although slow and somewhat ad-hoc in their original implementations, heuristic methods were moved front-and-center by the works of \citet{fazel2002matrix,burer2003nonlinear, burer2005local}. \cite{fazel2002matrix} observed that low-rank positive semidefinite matrices lie on the boundary of the PSD cone, and used this observation to justify a ``log-det'' heuristic, where a rank minimization objective is replaced with the function $\log\det(\bm{X}+\delta\mathbb{I})$. \cite{burer2003nonlinear,burer2005local} proposed implicitly modeling a rank constraint $\mathrm{Rank}(\bm{X}) \leq k$ by applying the non-linear reformulation $\bm{X}=\bm{U}\bm{V}^\top$, where $\bm{U}, \bm{V} \in \mathbb{R}^{n \times k}$ and eliminating $\bm{X}$, to obtain a problem which is non-convex in $(\bm{U}, \bm{V})$. Although originally solved using augmented Lagrangian techniques, modern implementations of the Burer-Monterio heuristic typically use alternating minimization \citep{jain2013low}, to enhance the methods scalability and ensure convergence towards a second-order critical point. 

The modern era of heuristics methods for low-rank matrix optimization was initiated by matrix completion and the Netflix competition \citep{bell2007lessons}. The magnitude of the data made available, challenged the aforementioned methods and led to new techniques such as \citet{wen2012solving}, who proposed a nonlinear successive over-relaxation approach that scales to $10,000 \times 10,000$ matrices,  \citet{recht2013parallel}, who proposed a stochastic gradient descent method which scales to $1,000,000 \times 1,000,000$ matrices; see \citet{udell2016generalized,nguyen2019low} for reviews of heuristic approaches.

\subsection{Contributions and Structure}
The key contributions of the paper are threefold: First, we propose using orthogonal projection matrices which satisfy $\bm{Y}^2=\bm{Y}$, the matrix analogue of binary variables which satisfy $z^2=z$, to model low-rank constraints via the non-linear equation $\bm{X}=\bm{Y}\bm{X}$. Under this lens, low-rank problems admit reformulations as optimization problems where some decision variables comprise a projection matrix. We term this family of problems \emph{Mixed-Projection Conic Optimization} (MPCO), in reference to mixed-integer optimization. To our knowledge, our approach is the first mathematical framework which solves low-rank optimization problems to certifiable optimality. Second, by leveraging regularization and strong duality we rewrite low-rank optimization problems as saddle-point problems over the space of orthogonal projection matrices, and propose an outer-approximation method to solve the saddle-point problem to certifiable optimality. Third, by analyzing the saddle-point problem, we derive new convex relaxations and rounding schemes which provide certifiably near optimal solutions in polynomial time in theory and rapidly in practice. {\color{black} Using a generic spatial branch-and-bound code, we are able to solve low-rank optimization problems exactly for matrices with $30$ rows and columns, and find near-exact solutions for matrices with up to $600$ rows and columns. We believe branch-and-bound schemes tailored to the set of projection matrices and dedicated semi-definite codes for the relaxations have to potential to increase the numerical scalability of MPCO even further and constitute an exciting future research direction.}

{\color{black}We note that the idea of using projection matrices to model low-rank constraints is not entirely new, as \citet{peng2005new} have proposed reformulating $k$-means clustering as a semidefinite optimization problem over the set of orthogonal projection matrices. However, our proposal has several key points of difference. Namely, (1) we consider optimizing over projection matrices directly, while \citet{peng2005new} use projection matrices as a vehicle to derive semidefinite relaxations, and (2) we use projection matrices to solve low-rank optimization problems that do not admit mixed-integer reformulations, while $k$-means clustering certainly admits a mixed-integer reformulation \citep{grotschel1989cutting}.}

The rest of the paper is laid out as follows:

In Section \ref{sec:unifyingperspective}, we show that projection matrices are a natural generalization of binary vectors to matrices. Inspired by a common tactic in cardinality constrained optimization, namely introducing binary variables to encode the support of the decision vector, we propose introducing a projection matrix to encode the image of the decision matrix and thereby model rank. We also investigate the complexity of low-rank optimization problems and show that rank minimization is as hard as the existential theory of the reals (i.e., deciding whether a semi-algebraic set is non-empty), and thus in PSPACE.

In Section \ref{sec:examples_plus_saddlepoint}, we derive the MPCO formulations of the aforementioned rank optimization problems.
By introducing a constraint on the spectral norm of $\bm{X}$ or a penalty on its Frobenius norm - the matrix analogs of big-$M$ constraints and perspective formulations \citep{gunluk2010perspective} respectively, we leverage strong duality, reformulate Problem \eqref{rankminproblem} as a saddle-point problem, and prove the resulting optimization problem admits a convex objective. 

We propose numerical algorithms to solve these MPCO problem to provable (near) optimality in Section \ref{sec:algorithmic}, by extending some of the most successful techniques from MICO. First, we propose an outer-approximation scheme for solving Problem \eqref{rankminproblem} exactly. 
Then, we obtain valid lower-bounds from solving its convex relaxations and propose an alternating minimization algorithm to do so. In addition, we prove that a singular value decomposition (SVD) followed by greedily rounding the eigenvalues provides certifiably near-optimal solutions in polynomial time. Finally, we propose a local-search strategy to improve the quality of the greedily rounded solution.

In Section \ref{sec:numeric}, we implement and numerically evaluate our proposed algorithms. On examples from matrix completion and sensor location, we demonstrate that methods proposed in this paper solve instances of Problem \eqref{rankminproblem} to certifiable optimality in minutes for $n$ in the tens. To our knowledge, our work is the first to demonstrate that moderately sized rank constrained problems can be solved to provable optimality in a tractable fashion. For $n$ in the hundreds, our proposal scales and provides in minutes solutions of higher quality than existing heuristics, such as nuclear norm minimization.

\section{From Cardinality to Rank: A Unifying Nonlinear Perspective}\label{sec:unifyingperspective}
Low rank constraints $\mathrm{Rank}(\bm{X})\leq k$ are a natural generalization of cardinality constraints $\Vert\bm{x}\Vert_0\leq k$ from vectors to matrices. Indeed, if $\bm{X}$ is a diagonal matrix then $\mathrm{Rank}(\bm{X})\leq k$ if and only if $\Vert\bm{X}\Vert_0\leq k$, and more generally $\mathrm{Rank}(\bm{X})\leq k$ if and only if $\Vert\sigma(\bm{X})\Vert_0\leq k$, where $\sigma(\bm{X})$ is the vector of singular values of $\bm{X}$. However, while cardinality and rank constraints are intimately linked, they are addressed using different algorithms. Namely, we can solve cardinality constrained problems with $100,000$s of variables to optimality \citep{bertsimas2017sparse}, while low-rank problems are dramatically harder and have not yet been solved to certifiable optimality for $n>10$ \citep{naldi2018solving}.

In our opinion, the difference between the community's understanding of cardinality and rank constraints has arisen because of two algorithmic barriers.
The first barrier is that rank constraints belong to a harder complexity class.
The second barrier arises because cardinality constraints can be represented using binary variables, while rank constraints cannot \citep[Corollary 4.1]{lubin2017regularity}. This presents a challenge for researchers, who have developed scalable methods for cardinality constraints by exploiting advances in mixed-integer {\color{black}conic} optimization {\color{black}(MICO)}, but cannot use these advances to address rank constraints. In this section, we question these barriers by characterizing the complexity of low-rank problems and proposing a new framework for modeling rank that generalizes {\color{black}(MICO)}.

\subsection{On the Complexity of Rank-Constrained Optimization}\label{ssec:complexity}
Existing studies of Problem \eqref{rankminproblem} typically claim that it is intractable, and support this claim by proving that it is NP-hard, by reduction from an NP-complete problem such as Boolean linear programming \citep[see, e.g.,][Section 7.3]{vandenberghe1996semidefinite}. In our opinion, this argument needs to be revisited, for two separate reasons. First, NP-hardness is a worst-case analysis statement. In practice, NP-hard optimization problems are often tractable. For instance, sparse regression can usually be solved to certifiable optimality with $100,000$s of features in minutes \citep{bertsimas2017sparse}. Second, there is no evidence that Problem \eqref{rankminproblem} is even in NP. Indeed, Problem \eqref{rankminproblem} cannot be represented using mixed-integer convex optimization \citep[Corollary 4.1]{lubin2017regularity}, while all $21$ of Karp's NP-complete problems can, and the best-known algorithms for Problem \eqref{rankminproblem} run in EXPTIME \citep{chistov1984complexity, naldi2018solving}.

In this section, we provide a more complete characterization of Problem \eqref{rankminproblem}'s complexity than is currently available in the literature. First, we demonstrate that it belongs to a different class than NP. In particular, it is  \textit{existential theory of the reals}-hard ($\exists\mathbb{R}$-hard; see \citet{renegar1992computational} for a general theory), i.e., as hard as any polynomial optimization problem, which implies that, if NP$ \subsetneqq \exists \mathbb{R}$, Problem \eqref{rankminproblem} is strictly harder than NP-complete problems. Second, we prove that Problem \eqref{rankminproblem} is actually in $\exists \mathbb{R}$.

We now demonstrate that Problem \eqref{rankminproblem} is existential theory of the reals complete (i.e., $\exists \mathbb{R}$-complete). We begin by reminding the reader of the definition of the $\exists\mathbb{R}$ complexity class \citep[c.f.][]{schaefer2013realizability}:
\begin{definition}
A decision problem belongs to the existential theory of the reals complexity class if it reduces to deciding whether a statement ``$
        \left(\exists x_1, ..., x_n \right)\phi(x_1, ..., x_n)
$''
is true or false, where $\phi(\cdot)$ is a quantifier-free Boolean formula involving polynomials equalities and inequalities, for instance, deciding the emptiness of a semi-algebraic set. We say a problem is $\exists \mathbb{R}$-hard if any problem in $\exists \mathbb{R}$ reduces to it.
\end{definition}

Note that 3-SAT $\in \exists \mathbb{R}$, so NP$\subseteq \exists \mathbb{R}$, and any statement in $\exists \mathbb{R}$ can be decided in PSPACE \citep{canny1988some}, so $\exists \mathbb{R} \subseteq$ PSPACE. 
To establish that Problem \eqref{rankminproblem} is $\exists\mathbb{R}$ hard, we require the following proposition, which is essentially a restatement of \citep[Theorem 3.1]{schaefer2013realizability} in the language of optimization.

\begin{proposition}\label{prop:graphreduction}
Let $G:=(V, E)$ be a graph, and $\ell(e)$ be the length of edge $e$. Then, deciding if $G$ can be embedded in $\mathbb{R}^2$ is $\exists\mathbb{R}$ complete, even when all edges have unit length.
\end{proposition}

By reducing Proposition \ref{prop:graphreduction}'s planar embedding problem to a Euclidean Distance Embedding problem, we obtain the following result (proof deferred to Appendix \ref{ssec:existentialtheoryofrealshard}):
\begin{theorem}\label{existentialtheoryofrealshard}
Problem \eqref{rankminproblem} is $\exists \mathbb{R}$-hard.
\end{theorem}

{Theorem \ref{existentialtheoryofrealshard} demonstrates that Problem \eqref{rankminproblem} is, from a traditional complexity theory perspective, at least as hard as any problem in $\exists \mathbb{R}$. However, its complexity status remains unresolved. Indeed, while \citet{candes2010matrix} have observed that Problem \eqref{rankminproblem} is in EXPTIME, it seems likely that $\exists \mathbb{R} \subset$ EXPTIME. We now address this matter, by proving that if $\mathcal{K}$ represents the semidefinite cone then Problem \eqref{rankminproblem} is in $\exists \mathbb{R}$, and hence $\exists \mathbb{R}$-complete; note that the examples listed in Section \ref{sec:examples} can all be rewritten as low-rank SDOs, so this result applies to all of the aforementioned examples (proof of theorem deferred to Appendix \ref{ssec:proofofthmexistentialcomp}).

\begin{theorem}\label{existentialtheorycompleteness}
Let $\mathcal{K}=S^n_+$ denote the $n \times n$ positive semidefinite cone. Then, Problem \eqref{rankminproblem} is in $\exists \mathbb{R}$, and hence $\exists \mathbb{R}$-complete.
\end{theorem}
\begin{remark}
Since $\exists \mathbb{R} \subseteq$ PSPACE $\subseteq$ EXPTIME, this upper bound improves upon the EXPTIME bound on Problem \eqref{rankminproblem}'s complexity stated by \citet{recht2010guaranteed, candes2010matrix} among others. Moreover, it seems unlikely to us that this bound can be further improved without settling fundamental questions in complexity theory (e.g., characterizing NP vs. $\exists \mathbb{R}$ vs. PSPACE vs. EXPTIME).
\end{remark}

\begin{remark}
Imposing an additional integrality constraint dramatically changes the complexity of Problem \eqref{rankminproblem}. Indeed, under the constraint $\bm{X} \in \mathbb{Z}^{n \times n}$, we can use a reduction from Hilbert's $\nth{10}$ problem to show that we \textit{cannot} decide in finite time whether Problem \eqref{rankminproblem}'s optimal objective is $0$, even if we know the objective is binary; see Appendix \ref{ssec:undecidableproof} for a proof.
\end{remark}

}

\subsection{Projection Matrices for Modeling Rank}\label{ssec:mixedprojconic}
As previously discussed, rank constraints can be seen as a generalization to the matrix case of cardinality constraints. For a vector $\bm{x} \in \mathbb{R}^n$, the cardinality constraint $\Vert \bm{x}\Vert_0 \leq k$ ensures that at most $k$ coordinates of $\bm{x}$ are non-zero, and can be modeled by introducing a vector of binary variables since
\begin{align}
    \Vert \bm{x}\Vert_0 \leq k \ \iff \exists \bm{z} \in \{0, 1\}^n:\  \bm{e}^\top \bm{z} \leq k,\ \bm{x}=\bm{z} \circ \bm{x},
\end{align}
where $\bm{z} \circ \bm{x}$ denotes the component-wise product of $\bm z$ and $\bm x$. Actually, non-linear constraints of the form ``$\bm{x}=\bm{z} \circ \bm{x}$'' where $\bm{z}$ is binary and $\bm{x}$ is continuous occur in a variety of mixed-integer optimization problems, far beyond cardinality constrained optimization. \citet{bertsimas2019unified} observed that such non-linear constraints ``$\bm{x}=\bm{z} \circ \bm{x}$'' actually lead to tractable optimization problems, provided that the overall objective is appropriately regularized. In particular, big-$M$ constraints \citep[][]{glover1975improved} and perspective reformulations \citep[][]{gunluk2010perspective} can be seen as appropriate regularizers. By building upon this observation and the work of several other authors \citep{fischetti2016redesigning, bertsimas2017sparse}, they successfully solve cardinality constrained problems at scale via a combination of branch-and-cut, randomized rounding and heuristic methods.

Unfortunately, rank constraints cannot be modeled using {\color{black}mixed-integer convex optimization} \citep[Corollary 4.1]{lubin2017regularity} and {\color{black}therefore MICO techniques cannot be applied ``out-of-the-box'' to address rank constraints.} Therefore, we now propose a new framework to model rank in optimization problems. Instead of a binary vector $\bm{z}$ to encode the support of $\bm{x}$, we introduce a projection matrix $\bm{Y}$ to capture the column space of $\bm{X}$ and obtain a similar non-linear reformulation.
{\color{black}
\begin{definition}
A matrix $\bm{Y} \in \mathbb{R}^{n \times n}$ is called a projection matrix if it satisfies the equality $\bm{Y}^2 = \bm{Y}$. In addition, if $\bm{Y}$ is symmetric, $\bm{Y}$ is called an orthogonal projection matrix.
\end{definition}
As symmetric matrices, orthogonal projection matrices are diagonalizable and their eigenvalues satisfy $\lambda_i^2 = \lambda_i$, i.e., are binary. As a result, the  Moore-Penrose pseudoinverse of an orthogonal projection $\bm{Y}$ is $\bm{Y}$ itself ($\bm{Y}=\bm{Y}^\dag$). In addition, since its eigenvalues are binary, the trace of $\bm{Y}$ equals the number of non-zero eigenvalues, i.e.,  $\mathrm{Rank}(\bm{Y})=\mathrm{tr}(\bm{Y})$. We are now in a position to link projection matrices and rank constraints.}
\begin{proposition}\label{prop:rankconstraintreformulation}
For any $\bm{X} \in \mathbb{R}^{n \times m}$,
$
 \mathrm{Rank}( \bm{X} ) \leq k \iff
  \exists \bm{Y} \in\mathcal{Y}_n \mbox{ : } \text{tr}(\bm{Y}) \leq k, \  \bm{X} = \bm{Y}\bm{X},
$
where $\mathcal{Y}_n:=\{\bm{P} \in S^n: \bm{P}^2 =\bm{P}\}$ is the set of $n\times n$ orthogonal projection matrices.
\end{proposition}
\proof{Proof of Proposition \ref{prop:rankconstraintreformulation}}
We prove the two implications successively.
\begin{itemize}
    \item Let $\bm{X}=\bm{U} \bm{\Sigma} \bm{V}^\top$, with ${\bm{U}}\in \mathbb{R}^{n\times k},\ {\bm{\Sigma}}\in \mathbb{R}^{k \times k},\ {\bm{V}}\in \mathbb{R}^{m \times k}$, be a singular value decomposition of $\bm{X}$ and
    define $\bm{Y}=\bm{U}\left(\bm{U}^\top \bm{U}\right)^{-1}\bm{U}^\top=\bm{U}\bm{U}^\top$. By construction, $\bm{X} = \bm{Y} \bm{X}$, since $\bm{U}^\top \bm{U}=\mathbb{I}$. Moreover, $\text{tr}(\bm{Y}) = \text{rank}(\bm{Y}) = \text{rank}(\bm{X}) \leq k$.
    \item Since $\bm{X} = \bm{Y} \bm{X}$,
 $\text{rank}(\bm{X}) \leq \text{rank}(\bm{Y}) =  \text{tr}(\bm{Y})  \leq k$.\hfill\Halmos
\end{itemize}
\endproof
\begin{remark}
In Proposition \ref{prop:rankconstraintreformulation}, the rank constraint is expressed via a trace constraint on $\bm{Y}$, the orthogonal projection onto the image or column space of $\bm{X}$. Alternatively, one could model the rank constraint via a matrix $\bm{Y}' \in \mathcal{Y}_m$ such that $\text{tr}(\bm{Y}') \leq k$ and $ \bm{X} = \bm{X}\bm{Y}'$. In this case, $\bm{Y}'$ encodes the projection onto the row space of $\bm{X}$. In practice, one could introduce both $\bm{Y}$ and $\bm{Y}'$ and obtain tighter formulations, at the price of introducing additional notation. We explore this idea in {\color{black}Appendix \ref{sec:A.linear.conjugate}}.
\end{remark}

Proposition \ref{prop:rankconstraintreformulation} suggests that projection matrices are to rank constraints what binary variables are to cardinality constraints. Indeed, similarities between the two are evident: binary variables $z$ are idempotent scalars which solve $z^2=z$, while projection matrices $\bm{Y}$ are idempotent matrices which solve $\bm{Y}^2=\bm{Y}$. Also, if $\bm{X}$ and $\bm{Y}$ are diagonal, Proposition \ref{prop:rankconstraintreformulation} recovers cardinality constrained optimization.

Over the past decades, extensive efforts have been devoted to improving the scalability of mixed-integer optimization. We believe that similar achievements can be obtained for rank constrained problems by adapting techniques from MICO to MPCO. In this direction, Table \ref{tab:analogybetweenbinaryandprojection} establishes a dictionary linking cardinality and rank constraints, and demonstrates that many of the techniques developed for binary convex optimization admit generalizations to MPCO, including the main results from our recent work \citep{bertsimas2019unified}. Note that we have not yet established most of the connections claimed in Table \ref{tab:analogybetweenbinaryandprojection}; this is the focus of the next two sections of the paper.

\begin{table}[h!]
\centering
\caption{Analogy between mixed-integer conic and mixed-projection conic optimization.}
\begin{tabular}{l l l}\toprule
 Framework & \citet{bertsimas2019unified} & This paper\\\midrule
 Parsimony concept & cardinality & rank\\
 Non-convex outer set & binaries & orthogonal projection matrices\\
 Strongly convex regularizer & $\ell_2^2$ penalty & Frobenius norm squared  \\
 Boundedness regularizer & $\ell_\infty$ norm & spectral norm\\
 Non-linear formulation & $\bm{x}=\bm{x} \circ \bm{z}$; $\bm{z} \in \{0, 1\}^n$  & $\bm{X}=\bm{Y}\bm{X}$, $\bm{Y} \in \mathcal{Y}_n$ \\\addlinespace[1mm]
Big-M formulation & $-M \bm{z} \leq \bm{x} \leq M \bm{z}$ & $\begin{pmatrix} M\bm{Y} & \bm{X}\\ \bm{X}^\top & M\mathbb{I} \end{pmatrix} \succeq \bm{0}$ \\\addlinespace[2mm]
Perspective formulation & $\begin{pmatrix} \theta_i & x_i\\ x_i & z_i\end{pmatrix} \succeq \bm{0}$ & $\begin{pmatrix} \bm{\theta} & \bm{X}\\ \bm{X}^\top & \bm{Y} \end{pmatrix} \succeq \bm{0}$\\\addlinespace[1mm]
Convex relaxation complexity & linear/second-order cone & semidefinite\\
Greedy rounding mechanism & coordinate-wise & singular value decomposition\\
 \bottomrule
\end{tabular}
\label{tab:analogybetweenbinaryandprojection}
\end{table}

\section{{\color{black} Regularization} and a Saddle-Point Reformulation}\label{sec:examples_plus_saddlepoint}
{\color{black} In this section, we prove that \eqref{rankminproblem_proj} can be reformulated as a saddle-point mixed-projection problem by leveraging regularization terms analogous to the big-$M$ and ridge regularization techniques from MICO, and derive their semidefinite relaxations, as summarized in Table \ref{tab:analogybetweenbinaryandprojection}. }

Throughout this paper, we let $\mathcal{Y}_n:=\{\bm{P} \in S^n: \bm{P}^2 =\bm{P}\}$ denote the set of $n \times n$ orthogonal projection matrices and $\mathcal{Y}_n^k:=\{\bm{P} \in S^n: \bm{P}^2 =\bm{P}, \mathrm{tr}(\bm{P}) \leq k\}$ denote projection matrices with rank at most $k$. {\color{black}
Although $\mathcal{Y}_n$ and $\mathcal{Y}^k_n$ do not commonly appear in the optimization literature, their convex hulls are well-studied, as we now remind the reader, by restating \citep[][Theorem 3]{overton1992sum}:
\begin{lemma}\label{lemma:convhull}
Let $\mathcal{Y}_n$ denote the $n \times n$ orthogonal projection matrices and $\mathcal{Y}_n^k$ denote the low-rank orthogonal projection matrices. Then, $\mathrm{Conv}(\mathcal{Y}_n)=\{\bm{P}: 0 \preceq \bm{P} \preceq \mathbb{I}\}$ and $\mathrm{Conv}(\mathcal{Y}_n^k)=\{\bm{P}: 0 \preceq \bm{P} \preceq \mathbb{I}, \mathrm{tr}(\bm{Y}) \leq k\}$. Moreover, the extreme points of $\mathrm{Conv}(\mathcal{Y}_n)$ are $\mathcal{Y}_n$, and the extreme points of $\mathrm{Conv}(\mathcal{Y}_n^k)$ are $\mathcal{Y}_n^k$.
\end{lemma}
}

\subsection{A Regularization Assumption}\label{ssec:regassumption}
By invoking Proposition \ref{prop:rankconstraintreformulation}, we rewrite Problem \eqref{rankminproblem} as the following mixed-projection conic problem:
\begin{equation}\label{rankminproblem_proj_noreg}
\begin{aligned}
    \min_{\bm{Y} \in \mathcal{Y}_n^k} \min_{\bm{X} \in \mathbb{R}^{n \times m}} \quad & \langle \bm{C}, \bm{X}\rangle+\lambda \cdot \mathrm{tr}(
\bm{Y}) \quad \text{s.t.} \quad \bm{A}\bm{X}=\bm{B},\ \bm{X}=\bm{Y}\bm{X},\ \bm{X} \in \mathcal{K}.
\end{aligned}
\end{equation}

Observe that Problem \eqref{rankminproblem_proj_noreg} has a two-stage structure which involves first selecting a low-rank projection matrix $\bm{Y}$ and second selecting a matrix $\bm{X}$ under the constraint $\bm{X}=\bm{Y}\bm{X}$. Moreover, selecting an optimal $\bm{X}$ given $\bm{Y}$ is \textit{easy}, because it involves solving a conic optimization problem under the linear constraint $\bm{X}=\bm{Y}\bm{X}$, while selecting an optimal $\bm{Y}$ is \textit{hard}, because $\mathcal{Y}_n^k$ is a non-convex set. Therefore, our modeling framework isolates the hardness of Problem \eqref{rankminproblem_proj_noreg} in $\mathcal{Y}_n^k$.

To cope with the non-linear constraints $\bm{X} = \bm{Y} \bm{X}$ in a tractable fashion, we augment the objective function in  \eqref{rankminproblem_proj_noreg} with a regularization term. Namely, we consider
\begin{equation}\label{rankminproblem_proj}
\begin{aligned}
    \min_{\bm{Y} \in \mathcal{Y}_n^k} \min_{\bm{X} \in \mathbb{R}^{n \times m}} \quad & \langle \bm{C}, \bm{X}\rangle+\Omega(\bm{X})+\lambda \cdot \mathrm{tr}(
\bm{Y}) \quad \text{s.t.} \quad \bm{A}\bm{X}=\bm{B},\ \bm{X}=\bm{Y}\bm{X},\ \bm{X} \in \mathcal{K}.
\end{aligned}
\end{equation}
where the regularization term $\Omega(\bm{X})$ satisfies the following assumption:
\begin{assumption}\label{assumption:regularizer}
In Problem \eqref{rankminproblem_proj}, the regularization term $\Omega(\bm{X})$ is one of:
\begin{itemize}
    \item A spectral norm penalty, $\Omega(\bm{X})=0$ if $\Vert \bm{X}\Vert_{\sigma} \leq M$ and $\Omega(\bm{X})=+\infty$ otherwise.
    \item A Frobenius norm penalty, $\Omega(\bm{X})=\frac{1}{2\gamma}\Vert \bm{X}\Vert_F^2$.
\end{itemize}
\end{assumption}
{\color{black}As we demonstrate in Section \ref{lowrankprojectionmatrix}, Assumption \ref{assumption:regularizer} is crucial for developing efficient low-rank algorithms, for the regularizer drives the convexity (see Theorem \ref{thm:innermax}) and smoothness (see Lemma \ref{lemma:lipschitz}) of the problem, and also make computationally cheap to evaluate subgradients readily accessible (Table \ref{tab:reg_conj}). The idea of leveraging regularization to optimize possibly non-smooth functions by obtaining computationally useful subgradients is classical \citep{nesterov2005smooth, nesterov2007smoothing} and therefore its effectiveness should not be surprising.}

The two regularizers are matrix analogues of the popular big-M constraints (constraints on the $\ell_\infty$ norm of the continuous variables) and ridge regularization (penalty on the $\ell_2^2$ norm) for vectors. In mixed-integer optimization, such regularization terms can efficiently cope with non-linear constraints between continuous and binary variables \citep{bertsimas2019unified} and motivate our current approach. Practically speaking, regularization can be a natural component of the original problem  \eqref{rankminproblem_proj_noreg}, otherwise we advocate for introducing it artificially, for it leads to tractable algorithms with moderate impact on the resulting solution. For instance, if $M$ is large enough so that the optimal solution to Problem \eqref{rankminproblem_proj_noreg}, $\bm{X}^\star$, satisfies $\Vert \bm{X}^\star\Vert_{\sigma} \leq M$, Problems \eqref{rankminproblem_proj} and \eqref{rankminproblem_proj_noreg} are equivalent.{\color{black}In Section \ref{ssec:computingM}, we develop a disciplined technique for computing such an $M$}. With the Frobenius norm penalty, the gap between  Problem \eqref{rankminproblem_proj}'s and \eqref{rankminproblem_proj_noreg}'s objective is at most $\tfrac{1}{2\gamma} \| \bm{X}^\star \|_F^2$, which can certainly be bounded whenever $\mathrm{tr}(\bm{X})$ is bounded, as often occurs in practice.

For ease of notation, we let
\begin{align*}
    g(\bm{X})=\langle \bm{C}, \bm{X}\rangle+\begin{cases} 0, \quad & \text{if} \ \bm{A}\bm{X}=\bm{B}, \ \bm{X} \in \mathcal{K},\\ +\infty, \quad & \text{otherwise}, \end{cases}
\end{align*}
denote the unregularized second-stage cost for a given $\bm{X}$. Therefore, Problem
\eqref{rankminproblem_proj} can be written as:
\begin{align}\label{outerproblem}
    \min_{\bm{Y} \in \mathcal{Y}_n^k} \ f(\bm{Y})+\lambda \cdot \mathrm{tr}(\bm{Y}),
\end{align}
\begin{align}\label{prob:definesubproblem}
    \text{where} \ f(\bm{Y}):=\min_{\bm{X} \in \mathbb{R}^{n \times m}} g(\bm{X})+\Omega(\bm{X}) \quad \text{s.t.}\quad \bm{X}=\bm{Y}\bm{X}
\end{align} yields a best choice of $\bm{X}$ given $\bm{Y}$. {\color{black}As we establish in this section, this turns out to be a computationally useful reformulation, for $f$ is convex in $\bm{Y}$ (see Theorem \ref{thm:innermax}) and Lipschitz continuous (see Lemma \ref{lemma:lipschitz}), and therefore the non-convexity in the problem has been isolated within the set $\mathcal{Y}^k_n$.}

Observe that both regularizers are coercive (i.e., ``blow up'' to $+\infty$ as $\Vert\bm{X}\Vert \rightarrow \infty$), and therefore render all unbounded solutions infeasible and ensure the compactness of the level sets of $\bm{X} \mapsto g(\bm{X})+\Omega(\bm{X})$. This alleviates two of the major issues with conic duality \citep[][Theorem 2.4.1]{ben2001lectures}. First, regularization ensures that optimal solutions to conic problems are attained \citep[see][Example 2.27, for a regularization-free counterexample]{blekherman2012semidefinite}.
Second, regularization ensures that infeasibility of a conic system is \textit{certifiable}\footnote{Unless the conic dual is also infeasible, this case is unimportant for our purposes, because it only arises when the original problem is itself infeasible for any $\bm{Y}$, which can be checked a priori.}, i.e., there is either a feasible solution or a certificate of infeasibility.
In general, such a procedure is not possible because a conic system could be infeasible but asymptotically feasible, i.e.,
    $$\nexists \bm{X}: \bm{A}\bm{X}=\bm{B}, \bm{X} \in \mathcal{K} \ \text{but} \ \exists \{ \bm{X}_t\}_{t=1}^\infty: \bm{X}_t \in \mathcal{K} \quad \forall t \ \text{with}\  \Vert \bm{A}\bm{X}_t -\bm{B}\Vert \rightarrow 0.$$
Here, the regularization term ensures that the set of feasible $\bm{X}$ (with objective at most $\theta_0 \in \mathbb{R}$) is a closed convex compact set. Therefore, $f(\bm{Y})$ cannot generate an asymptotically feasible problem.

{\color{black}
Finally, the two regularization functions in Assumption \ref{assumption:regularizer} satisfy a non-trivial property which turns out to be crucial in both proving that $f(\bm{Y})$ is convex and deriving our overall algorithmic strategy:
\begin{lemma}\label{lemma:conj.reg} Consider a regularization function $\Omega(\bm{X})$ satisfying Assumption \ref{assumption:regularizer}. There, there exists a Fenchel conjugate $\Omega^\star$ \citep[see, e.g.,][Chap. 3.3.1]{boyd2004convex} such that, for any projection matrix $\bm{Y} \in \mathcal{Y}_n$ and any matrix $\bm{\alpha}$, we have
\begin{align*}
    \min_{\bm{X}} \left\{ \Omega(\bm{Y} \bm{X}) + \langle \bm{\alpha},\bm{Y}\bm{X}  \rangle \right\} = \max_{\bm{V}_{11},\bm{V}_{22}} \: - \Omega^\star(\bm{\alpha}, \bm{Y}, \bm{V}_{11}, \bm{V}_{22}),
\end{align*}
and $\Omega^\star$ is linear in $\bm{Y}$ (see Table \ref{tab:reg_conj} for its explicit definition).
\end{lemma}
\proof{Proof of Lemma \ref{lemma:conj.reg}} We start with the Frobenius regularization case, $\Omega(\bm{X}) = \tfrac{1}{2\gamma} \| \bm{X} \|_F$ and
$\min_{\bm{X}} \left\{ \Omega(\bm{Y} \bm{X}) + \langle \bm{\alpha},\bm{Y}\bm{X}  \rangle \right\} =\tfrac{1}{2\gamma} \| \bm{Y} \bm{X} \|_F + \langle \bm{\alpha},\bm{Y}\bm{X}  \rangle$. Any solution to the minimization problem satisfies the first-order condition $\tfrac{1}{\gamma} \bm{Y} \bm{X} + \bm{Y} \bm{\alpha} = 0$. Hence, since $\bm{Y}^2 = Y$, $\bm{X}^\star = - \gamma \bm{Y} \bm{\alpha}$ satisfies the first-order condition and the optimal objective value is $- \Omega^\star(\alpha, \bm{Y}, \bm{V}_{11}, \bm{V}_{22}) =  -\tfrac{\gamma}{2} \langle \bm{\alpha}, \bm{Y} \bm{\alpha} \rangle$.

The spectral case is technically more challenging and detailed proofs are deferred to Appendix \ref{sec:A.linear.conjugate}. In the rectangular case, Lemma \ref{lemma:conj_spectral_rect} with $\bm{Y}' = \mathbb{I}_m$ yields
\begin{align*}
   \min_{\bm{X}} \left\{ \Omega(\bm{Y} \bm{X}) + \langle \bm{\alpha},\bm{Y}\bm{X}  \rangle \right\} = \max_{\bm{V}_{11},\bm{V}_{22}} \: -\frac{M}{2} \langle \bm{Y}, \bm{V}_{11} \rangle + \frac{M}{2} \langle \bm{I}_m, \bm{V}_{22} \rangle \mbox{ s.t. } \begin{pmatrix} \bm{V}_{11} & \bm{\alpha} \\ \bm{\alpha}^\top & \bm{V}_{22}\end{pmatrix} \succeq \bm{0}.
\end{align*}
In the symmetric case, Lemma \ref{lemma:conj_spectral_symmetric} states that
\begin{align*}
   \min_{\bm{X}} \left\{ \Omega(\bm{Y} \bm{X}) + \langle \bm{\alpha},\bm{Y}\bm{X}  \rangle \right\} = \max_{\bm{V}_{11},\bm{V}_{22} \succeq \bm{0}} \: -M \langle \bm{Y}, \bm{V}_{11} + \bm{V}_{22} \rangle \mbox{ s.t. } \bm{\alpha} = \bm{V}_{11} - \bm{V}_{22}. \quad \Halmos
\end{align*}
\endproof
}
\begin{table}
    \centering
    \caption{Regularization penalties and conjugates{\color{black}, as defined in Lemma \ref{lemma:conj.reg}.}}
    \begin{tabular}{c c c c c}\toprule
  Penalty   &  $\Omega(\bm{X})$ & $\Omega^\star(\bm{\alpha}, \bm{Y}, \bm{V}_{11}, \bm{V}_{22})$ & $\frac{\partial}{\partial Y_{i,j}} \Omega^\star(\bm{\alpha}, \bm{Y}, \bm{V}_{11}, \bm{V}_{22})$ \\\midrule
  Spectral norm ($\bm{X}$ rectangular) & $\begin{cases} 0, \ & \text{if} \ \Vert \bm{X}\Vert_\sigma \leq M,\\ +\infty, \ & \text{o.w.}, \end{cases}$ & $\frac{M}{2} \langle \bm{Y}, \bm{V}_{11} \rangle + \frac{M}{2} \langle \bm{I}_m, \bm{V}_{22} \rangle$ & $\frac{M}{2} V_{11,i,j}.$ \\
  & & $ \mbox{ s.t. } \begin{pmatrix} \bm{V}_{11} & \bm{\alpha} \\ \bm{\alpha}^\top & \bm{V}_{22}\end{pmatrix} \succeq \bm{0},$ & \\
  \midrule
  Spectral norm ($\bm{X}$ symmetric) & $\begin{cases} 0, \ & \text{if} \ \Vert \bm{X}\Vert_\sigma \leq M,\\ +\infty, \ & \text{o.w.}, \end{cases}$ & ${M} \langle \bm{Y}, \bm{V}_{11} + \bm{V}_{22} \rangle$ & $M(V_{11}+V_{22})_{i,j}.$ \\ \vspace{1mm}
  & & $\mbox{ s.t. } \bm{\alpha} = \bm{V}_{11} - \bm{V}_{22},$ & &\\
  & & $\bm{V}_{11}, \bm{V}_{2,2} \succeq \bm{0},$ &  \\ \midrule
Frobenius norm & $\frac{1}{2\gamma}\Vert \bm{X}\Vert_F^2$ & $\frac{\gamma}{2}\langle \bm{\alpha}, \bm{Y}\bm{\alpha}\rangle$ & $\frac{\gamma}{2}\langle \bm{\alpha}_i, \bm{\alpha}_j\rangle$ \\
  \bottomrule
    \end{tabular}
    \label{tab:reg_conj}
\end{table}
\subsection{A Saddle-Point Reformulation}\label{lowrankprojectionmatrix}
We now reformulate Problem \eqref{rankminproblem_proj} as a saddle-point problem. This reformulation is significant {\color{black} for two reasons. First, as shown in the proof of Theorem \ref{thm:innermax}, it leverages the nonlinear constraint $\bm{X}=\bm{Y}\bm{X}$ by introducing a new matrix of variables $\bm{V} \in \mathbb{R}^{n \times m}$ such that $\bm{V}=\bm{Y}\bm{X}$, giving:
\begin{align*}
    f(\bm{Y})=\min_{\bm{V}, \bm{X}}\{g(\bm{V})+\Omega(\bm{Y} \bm{X}): \bm{V}=\bm{Y}\bm{X}\},
\end{align*}a substitution reminiscent of the Douglas-Rachford splitting technique for composite convex optimization problems \citep[]{douglas1956numerical, eckstein1992douglas}—the proof of Theorem \ref{thm:innermax} shows that this substitution does not change the optimal objective value. Second, }it proves  that the regularizer $\Omega(\bm{X})$ drives the convexity and smoothness of $f(\bm{Y})$. To derive the problem's dual, we require:

\begin{assumption}\label{strongduality}
For each subproblem \eqref{prob:definesubproblem} generated by $f(\bm{Y})$ where $\bm{Y} \in \mathcal{Y}_n^k$, either the optimization problem is infeasible, or strong duality holds.
\end{assumption}

Assumption \ref{strongduality} holds under Slater's constraint qualification \citep[][Section 5.2.3]{boyd2004convex}. By invoking Assumption \ref{strongduality}, the following theorem reformulates \eqref{outerproblem} as a saddle-point problem:

\begin{theorem}\label{thm:innermax}
Suppose that Assumption \ref{strongduality} holds and $\Omega(\cdot)$ is either the spectral or Frobenius regularizer. Then, the following two optimization problems are equivalent:
\begin{align}
{\color{black} f(\bm{Y}) }&:= \min_{\bm{X} \in \mathbb{R}^{n \times m}} \quad  g(\bm{X})+\Omega(\bm{X}) \quad \text{\rm s.t.} \quad \bm{X}=\bm{Y}\bm{X}, \label{prob:dr1}\\
&= \max_{\bm{\alpha}, \bm{V}_{11}, \bm{V}_{22}} \quad h(\bm{\alpha})-\Omega^\star(\bm{\alpha}, \bm{Y}, \bm{V}_{11}, \bm{V}_{22}),\label{eqn:saddlepoint}
\end{align}
where $\displaystyle h(\bm{\alpha}) := \max_{\bm{\Pi} : \bm{C}-\bm{\alpha}-\bm{A}^\top\bm{\Pi} \in \mathcal{K}^\star} \langle \bm{b}, \bm{\Pi} \rangle$,
$\mathcal{K}^\star:=\{\bm{W}: \langle \bm{W}, \bm{X} \rangle \geq 0 \quad \forall \bm{X} \in \mathcal{K}\}$ denotes the dual cone to $\mathcal{K}$, and $\Omega^\star(\bm{\alpha}, \bm{Y}, \bm{V}_{11}, \bm{V}_{22})$ is defined in Table \ref{tab:reg_conj}.
\end{theorem}
\proof{Proof of Theorem \ref{thm:innermax}}
Let us fix $\bm{Y} \in \mathcal{Y}_n^k$, and suppose that strong duality holds for the inner minimization problem which defines $f(\bm{Y})$.
{\color{black}
To progress, we introduce a matrix $\bm{V} \in \mathbb{R}^{n \times m}$ such that $\bm{V}=\bm{Y}\bm{X}$ and obtain the relaxation:
\begin{align}\label{prob:dr2}
  \min_{\bm{X},\bm{V}} \: g(\bm{V}) +\Omega(\bm{Y} \bm{X}) \quad \mbox{  s.t.  } \bm{V} = \bm{Y}\bm{X}.
\end{align}
Let us verify that this relaxation is a valid substitution, i.e., that Problems \eqref{prob:dr1} and \eqref{prob:dr2} have the same optimal objective, $f(\bm{Y})$. If $\bm{X}$ is feasible for \eqref{prob:dr1}, then $(\bm{V} =\bm{X}, \bm{X})$ is obviously feasible for \eqref{prob:dr2} with same objective value. Similarly, let $(\bm{V}, \bm{X})$ be feasible for \eqref{prob:dr2}. $\bm{Y}\bm{V} = \bm{Y}^2 \bm{X} = \bm{Y} \bm{X} = \bm{V}$ since $\bm{Y}^2 = \bm{Y}$. Hence, $\bm{V}$ is feasible for \eqref{prob:dr1} with same objective value.}

{\color{black} Now, let} $\bm{\alpha}$ denote the dual variables associated with the coupling constraints $\bm{V} = \bm{Y}\bm{X}$. The minimization problem is then equivalent to its dual problem, which is given by:
\begin{align*}
  f(\bm{Y}) &= \max_{\bm{\alpha}} \: h(\bm{\alpha}) + \min_{\bm{X}} \left[ \Omega({\color{black}\bm{Y}} \bm{X}) + \langle \bm{\alpha}, \bm{Y}\bm{X} \rangle \right],
\end{align*}
where $h(\bm{\alpha}):=\inf_{\bm{V}} g(\bm{V})-\langle \bm{V}, \bm{\alpha}\rangle$ is, up to a sign, the Fenchel conjugate of $g$. By a standard application of Fenchel duality, it follows that
\begin{align*}
    h(\bm{\alpha})=\max_{\bm{\Pi}}\langle \bm{b}, \bm{\Pi}\rangle+\begin{cases} 0, \quad & \text{if} \ \bm{C}-\bm{\alpha}-\bm{A}^\top\bm{\Pi} \in \mathcal{K}^\star,\\ +\infty, \quad & \text{otherwise}. \end{cases}
\end{align*}

{\color{black} Finally, from Lemma \ref{lemma:conj.reg} we have $\displaystyle\min_{\bm{X}} \left\{ \Omega(\bm{Y} \bm{X}) + \langle \bm{\alpha},\bm{Y}\bm{X}  \rangle \right\} = \max_{\bm{V}_{11},\bm{V}_{22}} \: - \Omega^\star(\bm{\alpha}, \bm{Y}, \bm{V}_{11}, \bm{V}_{22})$, which concludes the proof. }

Alternatively, under either penalty, if the inner minimization problem defining $f(\bm{Y})$ is infeasible, then its dual problem is unbounded by weak duality.\footnote{Weak duality implies that the dual problem is either unfeasible or unbounded. Since the feasible set of the maximization problem does not depend on $\bm{Y}$, it is always feasible, unless the original problem is itself infeasible. Therefore, we assume without loss of generality that it is unbounded.}.
\hfill\Halmos \endproof
{\color{black}
\begin{remark}
In the unregularized case, i.e., $\Omega(\bm{X}) =0$, we can derive a similar reformulation:
\begin{align}
\min_{\bm{Y} \in \mathcal{Y}_n^k} \quad \max_{\bm{\alpha} \in \mathbb{R}^{n \times m}} \quad & h(\bm{\alpha})+{\lambda} \cdot \mathrm{tr}(\bm{Y}) \mbox{ s.t. } \quad \bm{Y} \bm{\alpha} = \bm{0}.
\end{align}
Under this lens, regularization of the primal problem is equivalent to a relaxation in the dual formulation: the hard constraint $\bm{Y} \bm{\alpha} = \bm{0}$ is penalized by $-\Omega^\star(\bm{\alpha}, \bm{Y}, \bm{V}_{11}, \bm{V}_{22})$.
\end{remark}

\begin{remark}
By Theorem \ref{thm:innermax} and Lemma \ref{lemma:conj.reg}, $f(\bm{Y})$ is convex as the point-wise maximum of functions which are linear in $\bm{Y}$.
\end{remark}

{

}
}
By Theorem \ref{thm:innermax}, when we evaluate $f(\hat{\bm{Y}})$, one of two alternatives occur. The first is that we have $f(\hat{\bm{Y}})< +\infty$ and there is some optimal $(\bm{\alpha}, \bm{V}_{11}, \bm{V}_{22})$. In this case, we construct the lower approximation $$f(\bm{Y}) \geq f(\hat{\bm{Y}})+\langle \bm{H}, \bm{Y}-\hat{\bm{Y}}\rangle,$$ where
{\color{black} $H_{i,j} = \frac{\partial}{\partial Y_{i,j}} \Omega^\star(\bm{\alpha}, \bm{Y}, \bm{V}_{11}, \bm{V}_{22})$ (see Table \ref{tab:reg_conj} for closed-form expression of the partial derivatives, which follow readily from Danskin's theorem \citep[see, e.g.,][Prop. B.22]{bertsekas1999nonlinear}).}
The second alternative is that $f(\hat{\bm{Y}})=+\infty$, in which case, by the conic duality theorem \citep[see][Chapter 2]{ben2001lectures} there exists a {\color{black}} $(\bm{\alpha}, \bm{\Pi})$ such that
\begin{align}
    \bm{C}-\bm{\alpha}-\bm{A}^\top \bm{\Pi} \in \mathcal{K}^\star, \ \text{and} \ \langle \bm{b}, \bm{\Pi}\rangle > \langle -\bm{H}, \hat{\bm{Y}} \rangle.
\end{align}
Under this alternative, we can separate $\hat{\bm{Y}}$ from the set of feasible $\bm{Y}$'s by imposing the cut $0 \geq \langle \bm{b}, \bm{\Pi}\rangle+\langle \bm{H}, \bm{Y} \rangle$. Under either alternative, we obtain a globally valid first-order underestimator of the form
\begin{align}
z f(\bm{Y})\geq h +\langle \bm{H}, \bm{Y}-\hat{\bm{Y}}\rangle,
\end{align}
where $z$, $h$ are defined as
\begin{align}
    z = \begin{cases} 1, & \text{ if } f(\hat{\bm{Y}}) < +\infty, \\
    0, & \text{ if } f(\hat{\bm{Y}}) = +\infty, \end{cases} \quad \mbox{ and } \quad h = \begin{cases} f(\hat{\bm{Y}}), & \text{ if } f(\hat{\bm{Y}}) < +\infty, \\
    \langle \bm{b}, \bm{\Pi}\rangle+\langle \bm{H}, \hat{\bm{Y}} \rangle, & \text{ if } f(\hat{\bm{Y}}) = +\infty. \end{cases}
\end{align}
This observation suggests that a valid numerical strategy for minimizing $f(\bm{Y})$ is to iteratively minimize and refine a piecewise linear underestimator of $f(\bm{Y})$ defined by the pointwise supremum of a finite number of underestimators of the form $z f(\bm{Y})\geq h +\langle \bm{H}, \bm{Y}-\bm{Y}\rangle$. Indeed, as we will see in Section \ref{sec:algorithmic}, this strategy gives rise to the global optimization algorithm known as outer-approximation.

{\color{black}
\paragraph{Smoothness} We now demonstrate that $f(\bm{Y})$ is smooth, in the sense of Lipschitz continuity, under a boundedness assumption on the size of the dual variables, which is a crucial property for ensuring the convergence of our global optimization methods and bounding the quality of our semidefinite relaxation and greedy rounding methods. Formally, the following result follows directly from Theorem \ref{thm:innermax}.
\begin{lemma}\label{lemma:lipschitz}
Let $\bm{Y}, \bm{Y}' \in \mathrm{Conv}(\mathcal{Y}^k_n)$ be on the convex hull of the orthogonal projection matrices. Then
\begin{align*}
    f(\bm{Y}) - f(\bm{Y}') \leq \Omega^\star(\bm{\alpha}^\star(\bm{Y}), \bm{Y}' - \bm{Y}, \bm{V}_{11}^{\star}(\bm{Y}), \bm{V}_{22}^{\star}(\bm{Y})).
\end{align*}
Moreover, suppose $\bm{\alpha}^\star(\bm{Y}), \bm{V}_{11}^{\star}(\bm{Y}),\bm{V}_{22}^{\star}(\bm{Y})$ can be bounded independently from $\bm{Y}$, i.e., $\Vert\bm{\alpha}^\star(\bm{Y})\Vert_\sigma \leq L_1$, $\Vert \bm{V}_{11}^{\star}(\bm{Y})\Vert_\sigma \leq L_2$, $\Vert \bm{V}_{22}^{\star}(\bm{Y})\Vert_\sigma \leq L_2$. Then, under spectral regularization we have
\begin{align}
    f(\bm{Y}) - f(\bm{Y}') \leq M \langle \bm{V}_{11}^{\star}(\bm{Y}), \bm{Y}'-\bm{Y}\rangle \leq M L_2 \Vert \bm{Y}'-\bm{Y}\Vert_*,
\end{align}
and under Frobenius regularization we have
\begin{align}
    f(\bm{Y}) - f(\bm{Y}') \leq \frac{\gamma}{2} \langle \bm{\alpha}^{\star\top}(\bm{Y}) \bm{\alpha}^{\star}(\bm{Y}), \bm{Y}'-\bm{Y}\rangle \leq \frac{\gamma}{2} L_1^2 \Vert \bm{Y}'-\bm{Y}\Vert_*,
\end{align}
where the bounds involving $L_1, L_2$ follow from Holder's inequality\footnote{\color{black}Namely, $\vert\langle \bm{X}, \bm{Y}\rangle\vert \leq \Vert \bm{X}\Vert_\sigma \Vert \bm{Y}\Vert_*$, since the $\Vert \cdot\Vert_\sigma$ and $\Vert \cdot\Vert_*$, as the matrix analogs of the $\ell_\infty$ and $\ell_1$ norms, are dual.}.
\end{lemma}
\begin{remark}
Section \ref{ssec:computingM} develops disciplined techniques for computing an $M$ such that the constraint $\Vert \bm{X}\Vert_\sigma \leq M$ in the primal does not alter the optimal objective. The same technique, applied to the dual, yields explicit bounds on $L_1$. Moreover, since there exists an optimal pair $(\bm{V}_{11}$, $\bm{V}_{22})$ which is an explicit functions of an optimal $\bm{\alpha}$, this translates into explicit bounds on $L_2$.
\end{remark}
}

\subsection{Semidefinite Relaxations}
To lower bound \eqref{outerproblem}'s objective, we {\color{black}invoke Lemma \ref{lemma:convhull} to }relax the non-convex constraint $\bm{Y} \in \mathcal{Y}_n^k$ to $$\bm{Y} \in \mathrm{Conv}\left(\mathcal{Y}_n^k\right)=\{\bm{Y} \in S^n: \bm{0} \preceq \bm{Y} \preceq \mathbb{I}, \mathrm{tr}(\bm{Y}) \leq k\}.$$ This yields the saddle-point problem
\begin{align}\label{prob:sdprelax}
    \min_{\bm{Y} \in \mathrm{Conv}\left(\mathcal{Y}_n^k\right)} \ \max_{\bm{\alpha}, \bm{V}_{11}, \bm{V}_{22} \in S^m} h(\bm{\alpha})-\Omega^\star \left(\bm{\alpha}, \bm{Y}, \bm{V}_{11},\bm{V}_{22}\right)+\lambda \cdot \mathrm{tr}(\bm{Y}).
\end{align}

Problem \eqref{prob:sdprelax} can in turn be reformulated as an SDO. Indeed, under Assumption \ref{strongduality}, we obtain a semidefinite formulation by taking Problem \eqref{prob:sdprelax}'s dual with respect to $\bm{\alpha}$. Formally, we have the following results (proofs deferred to Appendix \ref{subsec:perspectivecutproof} and \ref{proof:thm:spectralcovrelax} respectively):

\begin{lemma}\label{thm:perspectivereformulation}
Suppose that Assumption \ref{strongduality} holds. Then, strong duality holds between:
\begin{align}
    \min_{\bm{Y} \in \mathrm{Conv}\left(\mathcal{Y}_n^k\right)} \ \max_{\bm{\alpha} \in \mathbb{R}^{n \times m}} \quad & h(\bm{\alpha})-\frac{\gamma}{2}\langle \bm{\alpha}, \bm{Y}\bm{\alpha} \rangle+\lambda \cdot \mathrm{tr}(\bm{Y}),\\
    \min_{\bm{Y} \in \mathrm{Conv}\left(\mathcal{Y}_n^k\right)}\  \min_{\bm{X} \in \mathbb{R}^{n \times m}, \bm{\theta} \in S^n} \quad & g(\bm{X})+\frac{1}{2\gamma}\mathrm{tr}\left(\bm{\theta}\right)+\lambda \cdot \mathrm{tr}(\bm{Y})\quad
     \text{\rm s.t.} \quad \begin{pmatrix} \bm{\theta} & \bm{X}\\ \bm{X}^\top & \bm{Y}\end{pmatrix} \succeq \bm{0}. \label{bidualproblem}
\end{align}

\end{lemma}

\begin{lemma}\label{thm:spectralcovrelax}
Suppose that Assumption \ref{strongduality} holds. Then, strong duality holds between:
\begin{align}
    \min_{\bm{Y} \in \mathrm{Conv}\left(\mathcal{Y}_n^k\right)} \ \max_{\bm{\alpha} \in S^n, \bm{V}_{11}, \bm{V}_{22} \succeq \bm{0}}\quad & h(\bm{\alpha})-M\langle \bm{Y}, \bm{V}_{11}+ \bm{V}_{22}\rangle+\lambda \cdot \mathrm{tr}(\bm{Y}) \quad \text{\rm s.t.} \quad \bm{\alpha}=\bm{V}_{11}-\bm{V}_{22},\\
    \min_{\bm{Y} \in \mathrm{Conv}\left(\mathcal{Y}_n^k\right)} \ \min_{\bm{X} \in S^n} \quad & g(\bm{X})+\lambda \cdot \mathrm{tr}(\bm{Y}) \quad \text{\rm s.t.}\quad -M \bm{Y} \preceq \bm{X}\preceq M \bm{Y}. \label{bidualbigM}
\end{align}
\end{lemma}
We now offer some remarks on these bi-dual problems:
\begin{itemize}
    \item We can derive a more general version of {\color{black}Lemma} \ref{thm:spectralcovrelax} without the symmetry assumption on $\bm{X}$ in much the same manner, via the Schur complement lemma.
    \item Problem \eqref{bidualproblem}'s formulation generalizes the perspective relaxation from vectors to matrices. This suggests that \eqref{bidualproblem} is an efficient formulation for addressing rank constraints, as perspective formulations efficiently address cardinality constrained problems with conic quadratic \citep[][]{ gunluk2010perspective} or power cone \citep[][]{akturk2009strong} objectives{\color{black}, indeed, they provide a theoretical basis for scalable algorithms for sparse regression \citep{bertsimas2017sparse, hazimeh2020sparse}, sparse portfolio selection \citep{zheng2014improving,bertsimas2018scalable} and network design \citep{fischetti2016redesigning} problems among others.}.
\end{itemize}

{\subsection{Convex Penalty Interpretations of Relaxations}
In this section, we consider instances where rank is penalized in the objective only and interpret the above convex relaxations as penalty functions, in the tradition of \cite{fazel2002matrix, recht2010guaranteed}. In the presence of the Frobenius penalty, our first result generalizes the \textit{reverse Huber penalty} of \cite{pilanci2015sparse, dong2015regularization} from cardinality to rank objectives (proof deferred to Appendix \ref{ssec:appendgenrevhuber}).

\begin{lemma}\label{prop:generalizedrevhuber}
Suppose that Assumption \ref{strongduality} holds. Then, the following problems are equivalent:
\begin{align}
    & \min_{\bm{Y} \in \mathrm{Conv}(\mathcal{Y}_n)}  \ \min_{\bm{X} \in \mathbb{R}^{n \times m}, \bm{\theta} \in S^n} \quad g(\bm{X})+\frac{1}{2\gamma}\mathrm{tr}\left(\bm{\theta}\right)+\lambda \cdot \mathrm{tr}(\bm{Y})
\quad  \text{\rm s.t.} \quad  \begin{pmatrix} \bm{\theta} & \bm{X}\\ \bm{X}^\top & \bm{Y}\end{pmatrix} \succeq \bm{0},\label{prob:revhub1}\\
    & \min_{\bm{X} \in \mathbb{R}^{n \times m}} \quad  g(\bm{X})+\sum_{i=1}^n \min\left( \sqrt{\frac{2\lambda}{\gamma}}\sigma_i(\bm{X}),
     \lambda+\frac{\sigma_i(\bm{X})^2}{2\gamma} \right). \label{prob:revhub2}
\end{align}
\end{lemma}
{
\begin{remark}\color{black}

\begin{align*}
  \text{Since} \quad  \min_{0 \leq \theta \leq 1} \left[ \lambda \theta +\frac{t^2}{\theta}\right]=\begin{cases} 2\sqrt{\lambda}\vert t \vert, & \text{if} \ \vert t \vert \leq \sqrt{\lambda},\\
    t^2+\lambda, & \text{otherwise,}\end{cases}
\end{align*}
the proof of Lemma \ref{prop:generalizedrevhuber} reveals that Problems \eqref{prob:revhub1}-\eqref{prob:revhub2} are equivalent to minimizing
\begin{align}
\min_{\bm{X} \in \mathbb{R}^{n \times m}, \bm{\theta} \in \mathbb{R}^n: \ \bm{0} \leq \bm{\theta} \leq \bm{e}} \quad & g(\bm{X})+\sum_{i=1}^n \left(\lambda \theta_i +\frac{\sigma_i(\bm{X})^2}{2\gamma \theta_i}\right),
\end{align}
which applies the smooth penalty $t \rightarrow \lambda \theta+\frac{t^2}{2\gamma\theta}: 0 \leq \theta \leq 1$ to model the non-convex cost $t \rightarrow \lambda \Vert t\Vert_0+\frac{t^2}{2\gamma}$ incurred by each singular value of $\bm{X}$. Indeed, this smooth penalty is precisely the convex envelope of the non-convex cost function \citep[see, e.g.,][]{gunluk2010perspective}. 
Compared to other penalties for low-rank problems \citep{fan2001variable,zhang2010nearly}, this generalized Huber penalty is convex, amenable to efficient alternating minimization procedures (see Section \ref{ssec:altminderiv}) and could be of independent interest to the statistical learning community.
\end{remark}
}

{{\color{black}Lemma} \ref{prop:generalizedrevhuber} proposes an alternative to the nuclear norm penalty for approximately solving low-rank problems. This is significant, as many low-rank problems have constraints $\bm{X} \succeq \bm{0}, \mathrm{tr}(\bm{X})=k$ (e.g. sparse PCA {\color{black}\citep{d2007direct}}, $k$-means clustering {\color{black}\citep{peng2007approximating}}), and under these constraints a nuclear norm cannot encourage low-rank solutions \citep{zhang2013counterexample}, while {\color{black}Lemma} \ref{prop:generalizedrevhuber}'s penalty can.}

Our next results relate rank minimization problems with a spectral regularizer to the nuclear norm penalty, in both the square symmetric and the rectangular case (proofs deferred to Appendix \ref{ssec:proofofnucnormprops}):
\begin{lemma}\label{prop:rectcase1}
Suppose that Assumption \ref{strongduality} holds. Then, the following problems are equivalent:
\begin{align}
    \min_{\bm{Y} \in \mathrm{Conv}\left(\mathcal{Y}_n\right)} \ \min_{\bm{X} \in S^n} \quad & g(\bm{X})+\lambda\cdot \mathrm{tr}(\bm{Y}) \quad \text{\rm s.t.}\quad -M \bm{Y} \preceq \bm{X}\preceq M \bm{Y},\label{prob:specprob1}\\
    \min_{\bm{X} \in S^n} \quad & g(\bm{X})+\frac{\lambda}{M}\Vert \bm{X}\Vert_* \quad \text{\rm s.t.}\quad \Vert \bm{X}\Vert_\sigma \leq M. \label{prob:spectprob2}
\end{align}
\end{lemma}

\begin{lemma}\label{prop:rectcase2}
Suppose that Assumption \ref{strongduality} holds. Then, the following problems are equivalent:
\begin{align}
    \min_{\bm{Y} \in \mathrm{Conv}(\mathcal{Y}_n), \bm{Y}' \in \mathrm{Conv}(\mathcal{Y}_m)} \ \min_{\bm{X} \in \mathbb{R}^{n \times m}} \quad & g(\bm{X})+\frac{\lambda}{2}\mathrm{tr}(\bm{Y})+\frac{\lambda}{2}\mathrm{tr}(\bm{Y}') \quad \text{\rm s.t.}\quad \begin{pmatrix} M \bm{Y} & \bm{X} \\ \bm{X}^\top & M \bm{Y}' \end{pmatrix} \succeq \bm{0}, \label{prob:spectprob3}\\
    \min_{\bm{X} \in \mathbb{R}^{n \times m}} \quad & g(\bm{X})+\frac{\lambda}{M}\Vert \bm{X}\Vert_* \quad \text{\rm s.t.} \quad \Vert \bm{X}\Vert_\sigma \leq M. \label{prob:spectprob4}
\end{align}
\end{lemma}

{

}

}

{\color{black}
\subsection{Bounding the Spectral Penalty}\label{ssec:computingM}
We now present techniques for computing an $M$ such that the optimal values of \eqref{rankminproblem_proj_noreg} and \eqref{rankminproblem_proj} agree; these are essentially a generalization of similar techniques for logically constrained {\color{black}MICO}s \citep[][Section 2.3]{bertsimas2016best}. We first consider the positive semidefinite case, then develop the general case. Note that the $M$s obtained here are not, in general, the smallest possible—computing this quantity is NP-hard even for {\color{black}MICO}s \citep{kleinert2020there}.

\paragraph{Positive semidefinite case.}
Let $\bm{X} \in S^n_+$. Then, since $\mathrm{tr}(\bm{X})=\sum_{i=1}^n \lambda_i(\bm{X})\geq \lambda_1(\bm{X})$, the optimal value of the following problem gives a valid bound on $M$:
\begin{align}\label{prob:tr}
    M_{tr}:=\max_{\bm{X} \in S^n_+} \ \mathrm{tr}(\bm{X}) \ \text{s.t.} \ \bm{A}\bm{X}=\bm{B}, \bm{X} \in \mathcal{K}.
\end{align}
Alternatively, since the volume of $\bm{X}$ and its spectral radius are related, we can maximize
\begin{align}\label{prob:logdet}
    \bm{X}^\star=\arg\max_{\bm{X} \in S^n_+: \bm{A}\bm{X}=\bm{B}, \mathcal{X} \in \mathcal{K}} \log\det(\bm{X}).
\end{align}
We have $\det(\bm{X})\leq \left(\lambda_{\max}(\bm{X})\right)^n$, which implies we can set $M_{det}=\sqrt[n]{\det(\bm{X}^\star)}$. Alternatively, we could solve both \eqref{prob:tr} and \eqref{prob:logdet} and set $M=\min(M_{tr}, M_{det})$, which also gives a valid bound.

\paragraph{General case} Let $\bm{X} \in \mathbb{R}^{n \times m}$. Then, computing a valid $M$ is more expensive, because $\Vert \bm{X}\Vert_\sigma$, $\Vert \bm{X}\Vert_F$, and $\Vert\bm{X}\Vert_*$ are not mixed-integer convex representable when maximizing \citep[Corollary 4.1]{lubin2017regularity}. To progress, let us evaluate $n\cdot m$ values $M_{i,j}$, each computed by solving two conic problems:
\begin{align}
    M_{i,j}:=\max_{\bm{X} \in \mathbb{R^{n \times m}}: \bm{A}\bm{X}=\bm{B}, \mathcal{X} \in \mathcal{K}} \ \vert X_{i,j} \vert
\end{align}
Then, a valid $M$ is given by $\sum_{i,j}M_{i,j}$, since $\sum_{i,j} M_{i,j} \geq \Vert \bm{X}\Vert_1 \geq \Vert \bm{X}\Vert_\sigma$ \citep{boyd2004convex}.

\paragraph{Unbounded interpretation.} We remind the reader that interpreting the case where $M=+\infty$ requires caution. Indeed, when a feasible $\bm{X}$ and extreme ray $\bm{W}$ give an unbounded direction such that $$\mathrm{rank}(\bm{X}+\lambda \bm{W}) \leq k\quad \forall \lambda \geq 0$$ we have a certificate that no valid bound on $M$ exists. Alternatively, when $\mathrm{rank}(\bm{X}+\lambda \bm{W})\geq k+1$, we don't actually know whether a valid $M$ exists, since the set $\{\bm{X} \in S^n_+: \bm{A}\bm{X}=\bm{B}, \mathcal{X} \in \mathcal{K}\}$ could be unbounded, even while its intersection with a low-rank set is bounded. This difficulty also arises in the sparsity-constrained {\color{black}(MICO)} case however, and therefore should not be unexpected; it can be dealt with by cross-validating $M$/$\gamma$, which is usually acceptable since $\bm{A}, \bm{B}$ are usually estimated from data.
}

\section{Efficient Algorithmic Approaches}\label{sec:algorithmic}
In this section, we present an efficient numerical approach to solve Problem \eqref{rankminproblem} and its convex relaxations. The backbone is an outer-approximation strategy, embedded within a non-convex {\color{black}QCQO} branch-and-bound procedure to solve the problem exactly. We also propose rounding heuristics to find good feasible solutions, and semidefinite free methods for optimizing over \eqref{rankminproblem}'s convex relaxations.

{\color{black}The primary motivations for developing an outer-approximation procedure and solving mixed-projection problem as saddle-point problems are twofold. First, we are not aware of any exact solvers which address mixed-projection problems with semidefinite constraints. Instead, 
a decomposition strategy like outer-approximation can be readily implemented using a conjunction of \verb|Gurobi| (to solve non-convex quadratically constrained master problems) and \verb|Mosek| (to solve conic subproblems). Second, decomposition schemes for mixed-integer semidefinite problems typically outperform one-shot strategies \citep{belotti2013mixed}, so we expect - and observe in Section \ref{ssec:benchmarkalg1} - a similar comparison for mixed-projection optimization, \color{black} hence connecting the frameworks in both theory (see Table \ref{tab:analogybetweenbinaryandprojection}) and practice. 
}
\subsection{A Globally Optimal Cutting-Plane Method}\label{cuttingplanemethod}
The analysis in the previous section reveals that evaluating $f(\bm{Y})$ yields a globally valid first-order underestimator of $f(\cdot)$. Therefore, a numerically efficient strategy for minimizing $f(\bm{Y})$ is to iteratively minimize and refine a piecewise linear underestimator of $f(\bm{Y})$. This strategy is known as outer-approximation (OA), and was originally proposed by \citet{duran1986outer}. OA iteratively constructs underestimators of the following form at each iterate $t+1$:
\begin{align}
    f_{t+1}(\bm{Y})=\max_{1 \leq i \leq t} \left\{f(\bm{Y}_i)+\langle \bm{H}_i, \bm{Y}-\bm{Y}_i\rangle\right\}.
\end{align}
By iteratively minimizing $f_{t+1}(\bm{Y})$ and imposing the resulting cuts when constructing the next underestimator, we obtain a non-decreasing sequence of underestimators $f_t(\bm{Y}_t)$ and non-increasing sequence of overestimators $\min_{i \in [t]} f(\bm{Y}_i)$ which converge to an $\epsilon$-optimal solution within a finite number of iterations; see also Section \ref{lowrankprojectionmatrix} for details on cut generation. Indeed, since $\mathrm{Conv}\left(\mathcal{Y}_n^k\right)$ is a compact set and $f(\cdot)$ is an $L$-Lipschitz continuous function in $\bm{Y}$, OA never visits a ball of radius $\frac{\epsilon}{L}$ twice.

We now formalize this numerical procedure in Algorithm \ref{alg:cuttingPlaneMethod}, and state its convergence properties (proof of convergence deferred to Appendix \ref{ssec:proofofconv}):

\begin{algorithm*}
\caption{An outer-approximation method for Problem \eqref{outerproblem}}
\label{alg:cuttingPlaneMethod}
\begin{algorithmic}
\REQUIRE Initial solution $\bm{Y}_1$
\STATE $t \leftarrow 1 $
\REPEAT
\STATE Compute $\bm{Y}_{t+1}, \theta_{t+1}$ solution of
{\vspace{-3mm}
\begin{align*}
\min_{\bm{Y} \in \mathcal{Y}_n^k, \theta} \: \theta+\lambda\cdot \mathrm{tr}(\bm{Y}) \quad \mbox{ s.t. } z_i \theta \geq h_i +\langle \bm{H}_i, \bm{Y}-\bm{Y}_i\rangle \quad \forall i \in [t].
\end{align*}}\vspace{-5mm}
\STATE Compute $f(\bm{Y}_{t+1})$, $\bm{H}_{t+1}$, $z_{t+1}$, $d_{t+1}$
\UNTIL{$ f(\bm{Y}_t)-\theta_t \leq \varepsilon$}
\RETURN $\bm{Y}_t$
\end{algorithmic}
\end{algorithm*}

\begin{theorem}\label{proofofconvergence}
Suppose that Assumptions \ref{assumption:regularizer}-\ref{strongduality} hold, and that there exists some Lipschitz constant $L$ such that for any feasible $\bm{Y}, \bm{Y}' \in \mathrm{Conv}(\mathcal{Y}_n^k)$ we have:
$
    \vert f(\bm{Y})-f(\bm{Y}')\vert \leq L \Vert \bm{Y}-\bm{Y}'\Vert_F,
$
and for any feasibility cut $\langle \bm{H}_i, \bm{Y}-\bm{Y}_i\rangle +h_i \leq 0$ we have
$
    \vert \langle \bm{H}_i, \bm{Y}-\bm{Y}'\rangle\vert \leq L \Vert \bm{Y}-\bm{Y}'\Vert_F.
$
Let $\bm{Y}_t \in {\color{black}\mathcal{Y}_n^k}$ be a feasible solution returned by the $t^{\text{th}}$ iterate of Algorithm \ref{alg:cuttingPlaneMethod}, where
\begin{align*}
    t \geq \left(\frac{L k}{\epsilon}+1\right)^{n^2}.
\end{align*}
Then, $\bm{Y}_t$ is an $\epsilon$-optimal and $\epsilon$-feasible solution to Problem \eqref{rankminproblem_proj}.  Moreover, suppose that we set $\epsilon \rightarrow 0$. Then, any limit point of $\{\bm{Y}_t\}_{t=1}^\infty$ solves \eqref{rankminproblem_proj}.
\end{theorem}

\subsubsection{Optimizing Over Orthogonal Projection Matrices}
To successfully implement Algorithm \ref{alg:cuttingPlaneMethod}, we need to repeatedly solve optimization problems of the form \begin{align}
\min_{\bm{Y} \in \mathcal{Y}_n^k, \theta} \: \theta+\lambda\cdot \mathrm{tr}(\bm{Y}) \quad \mbox{ s.t. } z_i \theta \geq h_i +\langle \bm{H}_i, \bm{Y}-\bm{Y}_i\rangle \quad \forall i \in [t],
\end{align}
which requires a tractable representation of  $\mathcal{Y}^k_n$. {\color{black}Fortunately,} \verb|Gurobi| $9.0$ contains a globally optimal {\color{black}spatial} branch-and-bound method for general {\color{black}QCQO}s {\color{black}which recursively partitions the feasible region into boxes and invokes the ubiquitous McCormick inequalities to obtain valid upper and lower bounds on each box—see \citet[][]{gurobi2020} for a discussion of Gurobi's bilinear solver, \citet{belotti2013mixed} for a general theory of spatial branch-and-bound. Therefore,} we represent $\bm{Y}$ by introducing a matrix $\bm{U} \in \mathbb{R}^{n \times k}$ and requiring that $\bm{Y}=\bm{U}\bm{U}^\top$ and $\bm{U}^\top \bm{U}=\mathbb{I}$. This allows Algorithm \ref{alg:cuttingPlaneMethod} to be implemented by iteratively solving a sequence of {\color{black}QCQO}s and conic optimization problems. Moreover, to decrease the amount of branching required in each iteration of Algorithm \ref{alg:cuttingPlaneMethod}, we impose an outer-approximation of the valid constraint $\bm{Y} \succeq \bm{U}\bm{U}^\top$
. Specifically, we strengthen the formulation by imposing second-order cone relaxations of the PSD constraint. First, we require that the $2 \times 2$ minors in $\bm{Y}$ are non-negative {\color{black}, i.e., $Y_{i,j}^2 \leq Y_{i,i}Y_{j,j} \ \forall i,j \in [n]$}, as proposed in \cite{ahmadi2019dsos,bertsimas2019polyhedral}. Second, we require that the on-diagonal entries of $\bm{Y}\succeq \bm{U}\bm{U}^\top$ are non-negative {\color{black}i.e., $Y_{i,i}\geq \sum_{i=1}^k U_{i,t}^2 \ \forall i \in [n]$}. Finally, we follow \citet[Proposition 5]{atamturk2019rank} in taking a second-order cone approximation of the $2 \times 2$ minors in $\bm{Y}\succeq \bm{U}\bm{U}^\top$ {\color{black}i.e., $0 \geq \Vert \bm{U}_i\pm\bm{U}_j\Vert_2^2\pm 2Y_{i,j}- Y_{i,i}- Y_{j,j},\ \forall i,j \in [n]$}. All told, we have\footnote{It should be noted that this formulation is rather complicated because non-convex {\color{black}QCQO} solvers such as Gurobi currently do not model PSD constraints. If they did, we would supplant the second-order cone constraints with $\bm{Y}\succeq \bm{U}\bm{U}^\top$ and thereby obtain a simpler master problem.}:
\begin{equation}
    \begin{aligned}
    \min_{\bm{Y} \in S^n, \bm{U} \in \mathbb{R}^{n \times k}, \theta} \quad & \theta+\lambda\cdot \mathrm{tr}(\bm{Y}) \quad \mbox{ s.t. } z_i \theta \geq h_i +\langle \bm{H}_i, \bm{Y}-\bm{Y}_i\rangle\quad \forall i \in [t], \\
    & \bm{Y}=\bm{U}\bm{U}^\top, \bm{U}^\top \bm{U}=\mathbb{I}, Y_{i,i}Y_{j,j}\geq Y_{i,j}^2\ \forall i,j \in [n], Y_{i,i} \geq \sum_{t=1}^k U_{i,t}^2\ \forall i \in [n], \mathrm{tr}(\bm{Y})=k,\\
    & 0 \geq \Vert \bm{U}_i+\bm{U}_j\Vert_2^2-2Y_{i,j}-Y_{i,i}-Y_{j,j},\ 0 \geq  \Vert\bm{U}_i-\bm{U}_j\Vert_2^2+2Y_{i,j}- Y_{i,i}-Y_{j,j}\quad \forall i,j \in [n].
    \end{aligned}
\end{equation}

Finally, for a given $\bm{Y}, \bm{U}$, we strengthen this formulation by imposing second-order cone cuts of the form $\langle \bm{Y}-\bm{U}\bm{U}^\top, \bm{u}\bm{u}^\top \rangle \geq 0$, where $\bm{u}$ is the most negative eigenvector of $\bm{Y}-\bm{U}\bm{U}^\top$, as proposed by \cite{sherali2002enhancing}.

As described, a linear optimization problem over the set of orthogonal projection matrices is solved at each iteration, hence building a new branch-and-bound tree each time. We refer to this implementation as a ``multi-tree'' method. Although inefficient if implemented naively, multi-tree methods benefit from gradually tightening the numerical tolerance of the solver as the number of cuts increases.

To improve the efficiency of Algorithm \ref{alg:cuttingPlaneMethod}, one can integrate the entire procedure within a single branch-and-cut tree using lazy callbacks, as originally proposed in the context of MICO by \cite{quesada1992lp}. Henceforth, we refer to this implementation as a ``single-tree'' method. However, the benefit from using multi-tree over single-tree is not straightforward for it depends on how the method is engineered. We benchmark both implementations in Section \ref{ssec:benchmarkalg1}. 

{\color{black}
}

\subsubsection{A Simple Benchmark}\label{ssec:qcqopgurobi}
We now lay out a simple approach for solving low-rank problems exactly, which we will compare against in our numerical experiments. Rather than introducing an orthogonal projection matrix $\bm{Y}$, we let $\bm{X}=\bm{U}\bm{V}^\top$ where $\bm{U} \in \mathbb{R}^{n \times k}$ and $\bm{V} \in \mathbb{R}^{m \times k}$, and $\bm{U}$, $\bm{V}$ are both bounded in absolute value by big-M constraints of the form $\vert U_{i,j} \vert \leq 1, \vert V_{i,j} \vert \leq M$. Assuming that the objective and constraints are {\color{black}QCQO} representable, as occurs for all of the examples mentioned in the introduction, this formulation can then be optimized over using \verb|Gurobi|'s piecewise linear reformulation technique for general {\color{black}QCQO}s. Formally, a rank constraint $\mathrm{Rank}(\bm{X}) \leq k$ leads to:
\begin{align*}
    \min_{\bm{X} \in \mathbb{R}^{n \times m}, \bm{U} \in \mathbb{R}^{n \times k}, \bm{V} \in \mathbb{R}^{m \times k}} \quad & \Omega(\bm{X})+g(\bm{X})\\
    \text{s.t.} \quad & X \geq \bm{U}\bm{V}^\top-\epsilon \bm{E}, X \leq \bm{U}\bm{V}^\top+\epsilon \bm{E}, \ \Vert \bm{U}_i\Vert_2 \leq 1\quad \forall i \in [n], \Vert \bm{V}\Vert_\infty \leq M,
    \end{align*}
where $\bm{E}$ is a matrix of all ones. 
Note however that, as we observe in Section \ref{sec:numeric}, this approach is significantly less efficient than the previously described cutting-plane approaches.

\subsection{Lower bounds via Semidefinite Relaxations}\label{sec:sdr:rr}
To certify optimality, high-quality lower bounds are of interest and can be obtained by relaxing the non-convex constraint $\bm{Y} \in \mathcal{Y}_n^k$ to $\bm{Y} \in \mathrm{Conv}\left(\mathcal{Y}_n^k\right)$ {\color{black} to obtain a semidefinite relaxation as discussed in Lemma \ref{lemma:convhull}}. In addition to a valid lower bound on \eqref{outerproblem}'s objective, the optimal solution to the relaxation $\bm{Y}^\star$ is a natural candidate for a random rounding strategy, for {\color{black}stronger} convex relaxations lead to superior random rounding strategies. We will explore such rounding strategies in detail in the next section.

The convex relaxation yields the optimization problem \eqref{prob:sdprelax}
which can be solved using a cutting-plane method (see Section \ref{ssec:kelleyrelax}), an alternating minimization method (see Section \ref{ssec:altminderiv}) or reformulated as an SDO and solved as such.
Since Algorithm \ref{alg:cuttingPlaneMethod} is also an outer-approximation scheme, solving the convex relaxation via a cutting-plane method has the additional benefit of producing valid linear lower-approximations of $f(\bm{Y})$ to initialize Algorithm \ref{alg:cuttingPlaneMethod} with.


\subsubsection{Cutting-Plane Methods for Improving the Root Node Bound}\label{ssec:kelleyrelax} As mentioned previously, Problem \eqref{prob:sdprelax} can be solved by a cutting-plane method such as Kelley's algorithm \citep[see][]{kelley1960cutting}, which is a continuous analog of Algorithm \ref{alg:cuttingPlaneMethod} that solves Problem \eqref{outerproblem} over $\mathrm{Conv}(\mathcal{Y}^k_n)$, rather than $\mathcal{Y}^k_n$. 
The main benefit of such a cutting-plane method is that the cuts generated are valid for both $\mathrm{Conv}(\mathcal{Y}^k_n)$ and $\mathcal{Y}^k_n$, and therefore can be used to initialize Algorithm \ref{alg:cuttingPlaneMethod} and ensure that its \textit{initial} lower bound is equal to the semidefinite relaxation. As demonstrated by \citet{fischetti2016redesigning} in the context of MICO and facility location problems, this approach often accelerates the convergence of decomposition schemes by orders of magnitude. We present pseudocode in Appendix \ref{ssec:inout}

Figure \ref{fig:setupvalidation}'s left panel illustrates the convergence of Kelley's method and the \verb|in-out| method for solving the semidefinite relaxation of a noiseless matrix completion problem\footnote{The data generation process is detailed in Section \ref{ssec:compwithheurmethods}. Here, $n=100$, $p=0.25$, $r=1$, and  $\gamma=\frac{20}{p}$.}. Note that in our plot of the \verb|in-out| method on the continuous relaxation we omit the time required to first solve the SDO relaxation; this is negligible ($38.4$s) compared to the time required for either approach to solve the relaxation using cutting planes. Observe that the \verb|in-out| method's lower bound is both initially better and converges substantially faster to the optimal solution than Kelley's method. This 
justifies our use of the \verb|in-out| method over Kelley's method for a stabilizing cut loop in numerical experiments.

Once the relaxation is solved, the generated cuts are used to initialize Algorithm \ref{alg:cuttingPlaneMethod}. Figure \ref{fig:setupvalidation}'s right panel displays the convergence profile of the lower bound of Algorithm \ref{alg:cuttingPlaneMethod} initialized with cuts from Kelley's or the \verb|in-out| method (with a limit of $100$ cuts). We use a single-tree implementation of Algorithm \ref{alg:cuttingPlaneMethod}\footnote{We warm-start the upper bound with greedy rounding and the Burer-Monterio local improvement heuristic described in Section \ref{ssec:ubound}. To mitigate against numerical instability, we opted to be conservative with our parameters, and therefore turned Gurobi's heuristics off, set FuncPieceError and FuncPieceLength to their minimum possible values ($10^{-5}$ and $10^{-6}$), set the MIP gap to $1\%$ and the time limit for each solve to one hour.} and again a noiseless matrix completion setting\footnote{Here, $n=10$, $p=0.25$, $r=1$, and  $\gamma=\frac{5}{p}$.}.
We also consider the impact of using the SOC inequalities  $Y_{i,j}^2 \leq Y_{i,i}Y_{j,j}$ in the master problem formulation. Using the \verb|in-out| method and imposing the SOC inequalities are both vitally important for obtaining high-quality lower bounds from Algorithm \ref{alg:cuttingPlaneMethod}. Accordingly, we make use of both ingredients in our numerical experiments.
\begin{figure}[]
    \centering
    \begin{subfigure}[t]{.45\linewidth}
    \centering
            \includegraphics[scale=0.35]{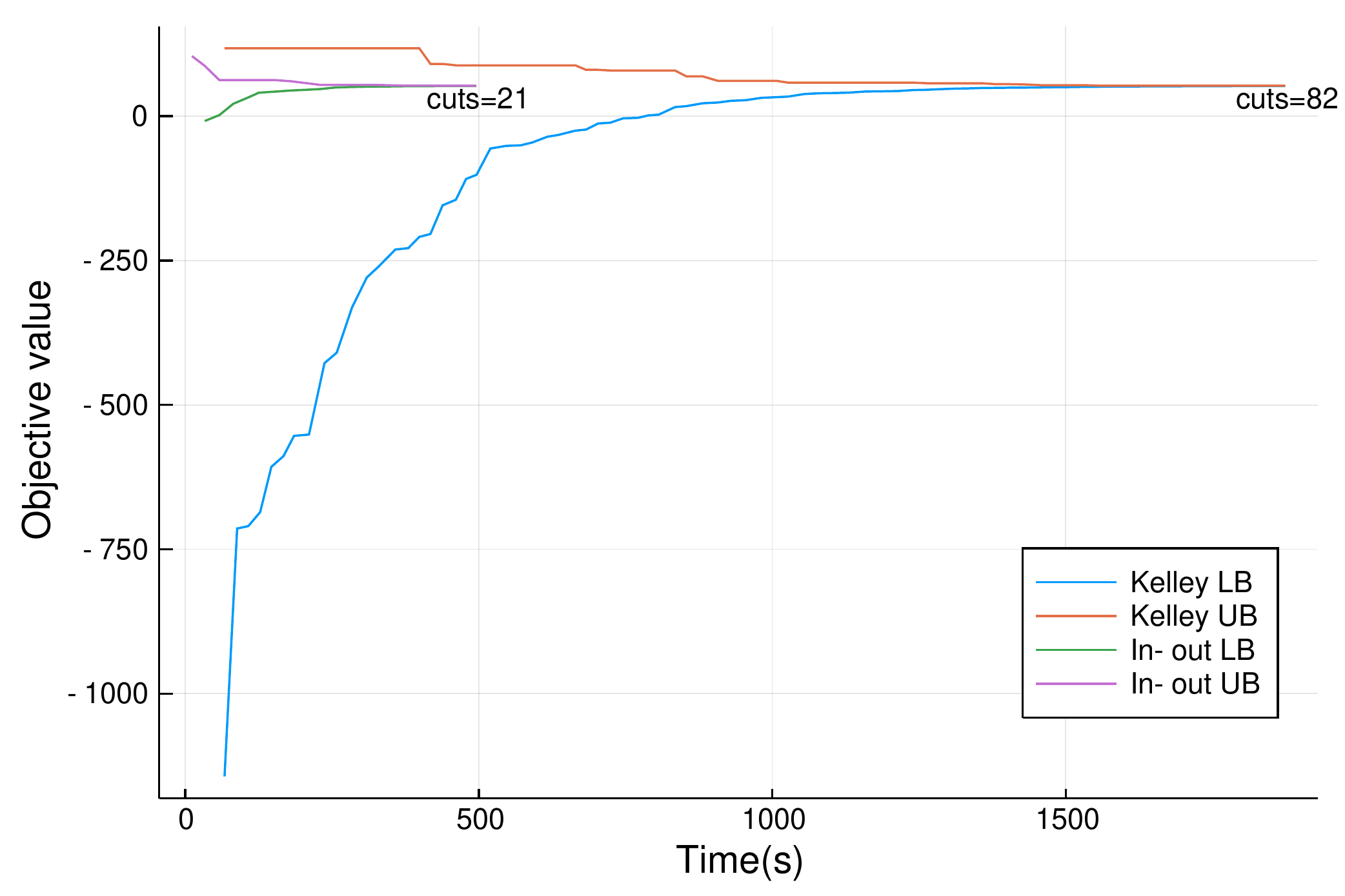}
    \end{subfigure}
        \begin{subfigure}[t]{.45\linewidth}
            \includegraphics[scale=0.35]{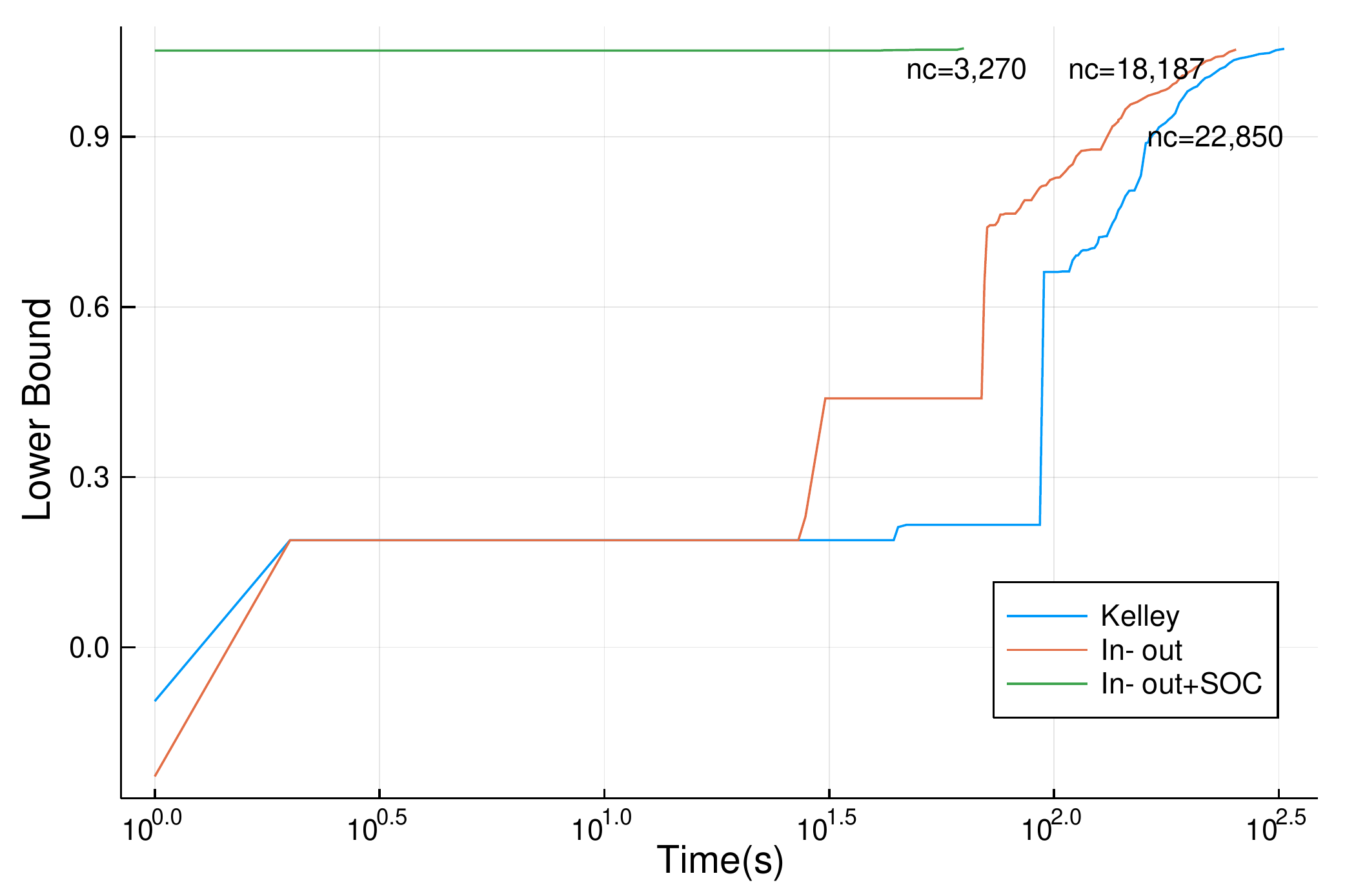}
    \end{subfigure}
   \caption{Convergence behavior of Kelley's method and the in-out method for solving the semidefinite relaxation of a synthetic matrix completion instance where $n=100$ (left), and lower bounds generated by a single-tree implementation of Algorithm \ref{alg:cuttingPlaneMethod} for a synthetic matrix completion instance where $n=10$ (right).}
   \label{fig:setupvalidation}
\end{figure}

\subsubsection{Solving the Semidefinite Relaxation at Scale via Alternating Minimization}\label{ssec:altminderiv}
In preliminary numerical experiments, we found that modern IPM codes such as \verb|Mosek| $9.0$ cannot optimize over the Frobenius/nuclear norm penalties when $n > 200$ on a standard laptop. As real-world low-rank problems are often large-scale, we now explore more scalable alternatives for optimizing over these penalties. As scalable alternatives for the nuclear norm penalty have been studied, we focus on the Frobenius penalty, and refer to \citep[]{recht2010guaranteed} for nuclear norm minimization.
We begin our analysis with the following result (proof deferred to Appendix \ref{ssec:proofoflemmaaltmin1}):

\begin{lemma}\label{lemma:optygivenx}
For any fixed $\bm{X}_t$ in Problem \eqref{bidualproblem}, an optimal choice of $\bm{\theta}$ is given by $\bm{\theta}^\star=\bm{X}_t^\top (\bm{Y}^\star)^\dag \bm{X}_t$, where $\bm{Y}^\star=\sum_{i=1}^n \rho_i^\star \bm{u}_i \bm{u}_i^\top$, $\bm{X}_t=\bm{U}\bm{\Sigma}\bm{V}^\top$ is an SVD of $\bm{X}_t$, and $\bm{\rho}^\star$ is an optimal solution to the following second order cone problem:
\begin{align}
    \min_{\bm{\rho} \in [0, 1]^n:\ \bm{e}^\top \bm{\rho} \leq k} \ \lambda \cdot \bm{e}^\top \bm{\rho} +\sum_{i=1}^n \frac{\sigma(\bm{X}_t)^2}{2\gamma \rho_i}.
\end{align}
\end{lemma}

As optimizing over $\bm{X}$ for a fixed $\bm{Y}_t$ is straightforward, {\color{black}Lemma} \ref{lemma:optygivenx} suggests a viable approach for optimize over the Frobenius norm penalty is alternating minimization \citep[AM; see][for a modern implementation]{beck2009fast}. By specializing \cite{beck2009fast}'s implementation of AM to the Frobenius norm penalty, we obtain an efficient numerical strategy for obtaining an optimal solution to \eqref{bidualproblem}, which we present in Algorithm \ref{alg:altmin}; we {\color{black}note that since $\langle \bm{X}\bm{X}, \bm{Y}^\dag\rangle$ is jointly convex in $\bm{X}, \bm{Y}$ (this follows directly from Lemma \ref{thm:perspectivereformulation}), alternating minimization converges to an optimal solution to the semidefinite relaxation under standard convergence conditions for block coordinate descent techniques for convex programs \citep[see, e.g.,][Section 3.7]{bertsekas1999nonlinear} such as the introduction of a proximal term.}

{\color{black}
We now discuss some enhancements to Algorithm \ref{alg:altmin} which improve its rate of convergence in practice.
\begin{itemize}
    \item Imposing a proximal regularization term in the objective, namely $+\frac{\tau}{2}\Vert \bm{X}-\bm{X}_t\Vert_F^2$, improves the rate of convergence of the method by stabilizing the iterates; we make use of this in our experiments.
    \item The method stalls when the eigenvalues of $\bm{Y}_t$ are near zero (a) due to numerical instability and (b) because $\bm{Y}_t$ is near the boundary of $\mathrm{Conv}(\mathcal{Y}^k_n)$. Therefore, to accelerate convergence, we require that $\lambda_{\min}(\bm{Y}) \geq \frac{K}{t}$ at the $t$th iterate, where $K\approxeq 10^{-2}$. In practice, this introduces very little error.
    \item Selecting an optimal $\bm{W}^{t+1}$ is generally much cheaper than selecting an optimal $\bm{V}^{t+1}$, since the former problem involves optimizing over $n$ eigenvalues, rather than $n^2$ variables. Therefore, efficient implementations of Algorithm \ref{alg:altmin} necessarily require efficient methods for obtaining $\bm{V}^{t+1}$. In the case of matrix completion, $g(\bm{X})$ is a quadratic form, which implies that obtaining $\bm{V}^{t+1}$ is equivalent to solving a linear system, which we do iteratively in our numerical experiments.
    \item We solve for $\bm{V}^{t+1}$ by solving the first-order optimality condition using a successive over-relaxation linear technique, or in rare instances where the linear system solver fails to converge we use \verb|Ipopt| to solve the QP's first-order optimality condition.
\end{itemize}}

\begin{algorithm*}[h!]
\caption{An Accelerated Alternating Minimization Algorithm \citep[c.f.][]{beck2009fast}}
\label{alg:altmin}
\begin{algorithmic}
\REQUIRE Initial solution $\bm{X}_1, \tau_1\gets 1$
\STATE $t \leftarrow 1, T_{\max} $
\REPEAT
\STATE Compute $\bm{W}^{t+1}$ solution of $\argminE_{\bm{Y} \in \mathrm{Conv}(\mathcal{Y}_n^k)} \ g(\bm{X}_t)+\frac{1}{2\gamma}\langle \bm{X}_t\bm{X}_t^\top, \bm{Y}^\dag \rangle$
\STATE Set $\bm{Y}^{t+1}=\bm{W}^t+\frac{\tau_{t-1}}{\tau_{t+1}}(\bm{W}_t-\bm{W}_{t-1})$
\STATE Compute $\bm{V}^{t+1}$ solution of $\argminE_{\bm{X} \in \mathbb{R}^{n \times m}} \ g(\bm{X})+\frac{1}{2\gamma}\langle \bm{X}\bm{X}^\top, \bm{Y}_t^\dag \rangle$
\STATE Set $\bm{X}^{t+1}=\bm{V}^t+\frac{\tau_{t-1}}{\tau_{t+1}}(\bm{V}_t-\bm{V}_{t-1})$
\STATE Set $\tau_{t+1}=\frac{1+\sqrt{1+4 \tau_t^2}}{2}$
\STATE If $t \mod 20=0$ compute dual bound at $\bm{Y}^{t+1}$ via Equation \eqref{max_conv_rankk}.
\STATE $t \leftarrow t+1 $
\UNTIL{$t > T_{\max}$ or duality gap $\leq \epsilon$}
\RETURN $\bm{X}_t, \bm{Y}_t$
\end{algorithmic}
\end{algorithm*}

{
To confirm that Algorithm \ref{alg:altmin} has indeed converged (at least approximately) to an optimal solution, we require a dual certificate. As optimizing over the set of dual variables $\bm{\alpha}$ for a fixed $\bm{Y}_t$ does not supply such a bound, we now invoke strong duality to derive a globally valid lower bound. Formally, we have the following result (proof deferred to Appendix \ref{sec:proofofpropositionconvminimax}):

\begin{lemma}\label{prop:conv_fullmax}
Suppose that Assumption \ref{strongduality} holds. Then, strong duality holds between:
\begin{align}\label{minmaxconv_rankkcase}
    \min_{\bm{Y} \in \mathrm{Conv}\left(\mathcal{Y}^k_n\right)} \ \max_{\bm{\alpha} \in \mathbb{R}^{n \times m}} \quad & h(\bm{\alpha})-\frac{\gamma}{2}\sum_{i=1}^n \sum_{j=1}^n Y_{i,j}\langle \bm{\alpha}_i, \bm{\alpha}_j \rangle, \\
     \max_{\substack{\bm{\alpha} \in \mathbb{R}^{n \times m}, \\\bm{U} \succeq \bm{0}, t \geq 0}} \quad & h(\bm{\alpha})-\mathrm{tr}(\bm{U})-kt \quad \text{\rm s.t.} \quad\bm{U}+\mathbb{I}t \succeq \frac{\gamma}{2}\bm{\alpha}\bm{\alpha}^\top. \label{max_conv_rankk}
\end{align}
\end{lemma}
{\color{black}Lemma} \ref{prop:conv_fullmax} demonstrates that Problem \eqref{rankminproblem}'s semidefinite relaxation is equivalent to maximizing the dual conjugate $h(\bm{\alpha})$, minus the $k$ largest eigenvalues of $\frac{\gamma}{2}\bm{\alpha}\bm{\alpha}^\top$. Moreover, as proven in the special case of sparse regression by \citet{bertsimas2019sparse}, one can show that if the $k$th and $k+1$th largest eigenvalues of $\bm{\alpha}\bm{\alpha}^\top$ in a solution to \eqref{minmaxconv_rankkcase} are distinct then Problem \eqref{minmaxconv_rankkcase}'s lower bound is tight. 

}

\subsection{Upper Bounds via Greedy Rounding}\label{ssec:ubound}
We now propose a greedy rounding method for rounding $\bm{Y}^\star$, an optimal $\bm{Y}$ in a semidefinite relaxation of Problem \eqref{rankminproblem_proj}, to obtain certifiably near-optimal solutions to Problem \eqref{rankminproblem_proj} quickly. Rounding schemes for approximately solving low-rank optimization problems by rounding their SDO relaxations have received a great deal of attention since they were first proposed by \citet{goemans1995improved}. Our analysis is, however, more general than typically conducted when solving low-rank problems, as it involves rounding a projection matrix $\bm{Y}$, rather than rounding $\bm{X}$, and therefore is able to generalize to the rank-$k$ case for $k>1$, which has historically been challenging.

Observe that for any feasible $\bm{Y} \in \mathrm{Conv}(\mathcal{Y}_n)$, $0 \leq \lambda_i(\bm{Y}) \leq 1$ for each eigenvalue of $\bm{Y}$, and $\bm{Y}$ is a projection matrix if and only if its eigenvalues are binary. Combining this observation with the Lipschitz continuity of $f(\bm{Y})$ in $\bm{Y}$ suggests that high-quality feasible projection matrices can be found in the neighborhood of a solution to the semidefinite relaxation, and a good method for obtaining them is to greedily round the eigenvalues of $\bm{Y}$. Namely, let $\bm{Y}^\star$ denote a solution to the semidefinite relaxation \eqref{prob:sdprelax}, $\bm{Y}^\star=\bm{U}\bm{\Lambda}^\star\bm{U}^\top$ be a singular value decomposition of $\bm{Y}^\star$ {\color{black}such that $\bm{\Lambda}$ is a diagonal matrix with on-diagonal entries $\Lambda_{i,i}$}, and $\bm{\Lambda}_{greedy}$ be a diagonal matrix obtained from rounding up (to 1) $k$ of the highest diagonal coefficients of $\bm{\Lambda}^\star$, and rounding down (to 0) the $n-k$ others{\color{black}, with diagonal entries $\Lambda_{i,i}:=(\Lambda_{greedy})_{i,i}$}. We then let $\bm{Y}_{greedy} = \bm{U} \bm{\Lambda}_{greedy} \bm{U}^\top$.
We now provide guarantees on the quality of the greedily rounded solution (proof deferred to Appendix \ref{ssec:proofoftheoremrandround}):
\begin{theorem}\label{thm:randrounding2}
Let $\bm{Y}^\star$ denote a solution to the semidefinite relaxation \eqref{prob:sdprelax}, $\bm{Y}^\star=\bm{U}\bm{\Lambda}\bm{U}^\top$ be a singular value decomposition of $\bm{Y}^\star$, $\mathcal{R}$ denote the indices of strictly fractional diagonal entries in $\bm{\Lambda}$, and $\bm{\alpha}^\star(\bm{Y})$ denote {\color{black}an optimal} choice of $\bm{\alpha}$ for a given $\bm{Y}${\color{black}, i.e., $$\alpha^\star(\bm{Y})\in \arg\max_{\bm{\alpha}}\left\{\max_{\bm{V}_{11}\bm{V}_{22}}h(\bm{\alpha})-\Omega^\star \left(\bm{\alpha}, \bm{Y}, \bm{V}_{11},\bm{V}_{22}\right)\right\}.$$}Suppose that for any $\bm{Y} \in \mathcal{Y}_n^k$, we have $\sigma_{\max}(\bm{\alpha}^\star(\bm{Y})) \leq L$. Then, {\color{black} any valid rounding of $\bm{Y}^\star$ which preserves the relaxation's eigenbasis, i.e., $\bm{Y}_{rounded}=\bm{U}\bm{\Lambda}_{rounded}\bm{U}^\top$ where $\bm{Y}^\star=\bm{U}\bm{\Lambda}\bm{U}^\top$ and $\bm{\Lambda}_{rounded}$ is a diagonal matrix with binary diagonal entries $\Lambda_{i,i}^{rounded}$ such that $\mathrm{tr}(\bm{\Lambda}_{rounded}) \leq k$, satisfies}
\begin{align}\label{eqn:concentrationbound1}
    {\color{black}f(\bm{Y}_{rounded})}-f(\bm{Y}^\star) \leq \frac{\gamma}{2}L^2 \vert\mathcal{R}\vert \max_{\bm{\beta} \geq \bm{0}: \Vert \bm{\beta} \|_1 \leq 1} \sum_{i \in \mathcal{R}} (\Lambda_{i,i}^\star-{\color{black}\Lambda_{i,i}^{rounded}})\beta_i,
\end{align}
under the Frobenius penalty and
\begin{align}\label{eqn:eqn:concentrationbound2}
    {\color{black}f(\bm{Y}_{rounded})}-f(\bm{Y}^\star) \leq M L \vert\mathcal{R}\vert \max_{\bm{\beta} \geq \bm{0}: \Vert \bm{\beta} \|_1 \leq 1} \sum_{i \in \mathcal{R}} (\Lambda_{i,i}^\star-{\color{black}\Lambda_{i,i}^{rounded}})\beta_i,
\end{align}
for the spectral penalty. {\color{black}Moreover, let $\bm{Y}_{greedy} = \bm{U} \bm{\Lambda}_{greedy} \bm{U}^\top$ be an instance of $\bm{Y}_{rounded}$ obtained by setting
$\Lambda_{i,i}=1$ for $k$ of the highest diagonal coefficients in $\bm{\Lambda}^\star$. Then, the above bounds imply that }$0 \leq f(\bm{Y}_{greedy})-f(\bm{Y}^\star) \leq \epsilon$, where $\epsilon=M L \min(\vert \mathcal{R}\vert, n-k)$ for the spectral penalty and $\epsilon=\frac{\gamma}{2} \min(\vert \mathcal{R}\vert, n-k) L^2$ for the Frobenius penalty.
\end{theorem}
{\color{black} This result calls for multiple remarks:
\begin{itemize}
    \item When the relaxation gap ${\color{black}f(\bm{Y}_{\text{greedy}})}-f(\bm{Y}^\star)=0$, and the optimal solution to the relaxation, $\bm{Y}^\star$, is unique, $\vert \mathcal{R}\vert=0$. This justifies retaining $\mathcal{R}$ in the bound, rather than replacing it with $n$.
    \item The techniques introduced in Section \ref{ssec:computingM} for computing an $M$ so that an optimal solution $\bm{X}^\star$ obeys $\Vert \bm{X}^\star\Vert_\sigma \leq M$, also apply to computing an explicit $L$ such that $\Vert \bm{\alpha}^\star\Vert_\sigma \leq L$ in the above bound.
\item The rounding technique is robust, because it minimizes the worst-case Lipschitz upper bound, under the assumption $\sigma_{\max}(\bm{\alpha}^\star) \leq L$ (i.e., we have no information about which coordinate\footnote{If we had this information then, as the proof of Theorem \ref{thm:randrounding2} suggests, we would greedily round to one $k$ of the indices with the largest values of $L_i \Lambda_{i,i}^\star$.} has the largest Lipschitz upper bound). For instance, under Frobenius regularization the bound is
    \begin{align}
        f(\bm{Y}_{rounded})-f(\bm{Y}^\star) \leq \frac{\gamma}{2}L^2 \vert\mathcal{R}\vert \max_{\bm{\beta} \geq \bm{0}: \Vert \bm{\beta}{\color{black} \Vert}_1 \leq 1} \sum_{i \in \mathcal{R}} (\Lambda_{i,i}^\star-{\color{black}\Lambda_{i,i}^{rounded}})\beta_i,
    \end{align}
     which is minimized over $\bm{\Lambda}^{rounded}: \mathrm{tr}(\bm{\Lambda}^{rounded})\leq k$ by solving:
     \begin{align}
      \min_{\bm{\lambda} \in \mathcal{S}^k_n}\max_{\bm{\beta} \geq \bm{0}: \Vert \bm{\beta} {\color{black} \Vert}_1 \leq 1}\frac{\gamma}{2}L^2 \vert\mathcal{R}\vert  \sum_{i \in \mathcal{R}} (\Lambda_{i,i}^\star-\lambda_{i})\beta_i,
    \end{align}
    i.e., rounding greedily. This interpretation suggests that greedy rounding never performs too badly.
\end{itemize}
}

To improve the greedily rounded solution, we implement a local search strategy which obtains even higher quality warm-starts. Namely, a variant of the popular Burer-Monterio (BM) heuristic \citep[][]{burer2003nonlinear}, which seeks low-rank solutions $\bm{X}$ by applying a non-linear factorization $\bm{X}=\bm{U}\bm{V}^\top$, where $\bm{U} \in \mathbb{R}^{n \times l}, \bm{V} \in \mathbb{R}^{m \times k}$ and iteratively optimizing over $\bm{U}$ for a fixed $\bm{V}$ (resp. $\bm{V}$ for a fixed $\bm{U}$) until convergence to a local optima occurs. This strategy improves our greedily rounded solution because we initially set $\bm{U}$ to be the square root of $\bm{Y}_{greedy}$ and optimize over $\bm{V}$; recall that if $\bm{Y}$ is a projection matrix we have $\bm{Y}=\bm{U}\bm{U}^\top$ and $\bm{X}=\bm{U}\bm{\Sigma}\bm{V}^\top$ for some singular value decomposition $\bm{U}, \bm{\Sigma}, \bm{V}^\top$.


\section{Numerical Experiments} \label{sec:numeric}
In this section, we evaluate the algorithmic strategies derived in the previous section, implemented in \verb|Julia| 1.3 using \verb|JuMP.jl| $0.20.1$, \verb|Gurobi| $9.0.1$ to solve the non-convex {\color{black}QCQO} master problems\footnote{We remark that Gurobi solves the non-convex {\color{black}QCQO} master problems by translating them to piecewise linear optimization problems. Since rank constraints are not MICO representable, this introduces some error. To mitigate against this error, we set the Gurobi parameters FuncPieceError and FuncPieceLength to their minimum possible values ($10^{-6}$ and $10^{-5}$ respectively). Additionally, we set NonConvex to $2$, and otherwise use default Gurobi/Mosek parameters.}, and Mosek $9.1$ to solve the conic subproblems/continuous relaxations. Except where indicated otherwise, all experiments were performed on a Intel Xeon E5---2690 v4 2.6GHz CPU core using 32 GB RAM. {\color{black}To bridge the gap between theory and practice, we have made our code freely available on \verb|Github| at \verb|github.com/ryancorywright/MixedProjectionSoftware|.}

We evaluate the different ingredients of our numerical strategy on a matrix completion example: First, we solve the semidefinite relaxation by implementing Algorithm \ref{alg:altmin} and demonstrate its increased scalability over \verb|Mosek|'s IPM in Section \ref{ssec:altminexper}. From the solution of the relaxation, our rounding and local search heuristics then provide near-optimal solutions that outperform state-of-the-art heuristic methods, as discussed in Section \ref{ssec:compwithheurmethods}. We implement Algorithm \ref{alg:cuttingPlaneMethod}, benchmark its performance and, for the first time, solve low-rank matrix completion to certifiable optimality in Section \ref{ssec:benchmarkalg1}. {\color{black}In Section \ref{ssec:impactofreg}, we explore the role which regularization plays in our numerical strategy, by showing that increasing the amount of regularization in Problem \eqref{rankminproblem} decreases the relative gap, the problem's complexity, and the amount of time required to solve the problem to optimality. Finally, in Section \ref{ssec:syntheticcoordrec} we solve sensor location problems to certifiable optimality.}

\subsection{Exploring the Scalability of the Convex Relaxations}\label{ssec:altminexper}
In this section, we explore the relative scalability of \verb|Mosek|'s interior point method and Algorithm \ref{alg:altmin}.

We consider convex relaxations of matrix completion problems. Similarly to \cite{candes2010matrix}, we generate two low-rank matrices $\bm{M}_L, \bm{M}_R \in \mathbb{R}^{n \times r}$ with i.i.d. $\mathcal{N}(0,1)$ entries, and attempt to recover the matrix $\bm{M}=\bm{M}_L \bm{M}_R^\top$ given a proportion $p$ of its observations. Here, we fix $p=0.25$ and $k=r=5$, vary $n$, and set $\gamma=\frac{20}{p}$ {\color{black} where we scale $\gamma$ proportionally to $1/p$ so that the relative importance of $\Vert \bm{X}\Vert_F^2$ and $\sum_{(i,j) \in \Omega}(X_{i,j}-A_{i,j})^2$ remains constant with $p$}.

We solve the continuous relaxation
\begin{align}
    \min_{\bm{X} \in \mathbb{R}^{n \times n}, \bm{Y} \in \mathrm{Conv}(\mathcal{Y}_n^k), \bm{\theta} \in S^n} \ \frac{1}{2 \gamma}\mathrm{tr}(\bm{\theta})+\sum_{(i,j) \in \Omega}(X_{i,j}-A_{i,j})^2 \quad \text{s.t.} \quad \begin{pmatrix} \bm{\theta} & \bm{X}\\ \bm{X}^\top & \bm{Y} \end{pmatrix} \succeq \bm{0}.
\end{align}
Table \ref{tab:comparison} reports the time required by Algorithm \ref{alg:altmin} to obtain a solution with a relative duality gap of $0.1\%$. To evaluate numerical stability, we also report the relative MSE of the greedily rounded solution; experiments {\color{black}where} $n \leq 250$ were run on a standard MacBook pro with $16$GB RAM, while larger experiments were run on the previously described cluster with $100$GB RAM.

\begin{table}[h!]
\centering\footnotesize
\caption{Scalability of convex relaxations, averaged over $5$ matrices. Problem is regularized with Frobenius norm and $\gamma=\frac{20}{p}$. ``-'' indicates an instance could not be solved with the supplied memory budget.}
\begin{tabular}{@{}l r r r r r | l r r @{}} \toprule
$n$  & \multicolumn{2}{c@{\hspace{0mm}}}{Mosek} &\multicolumn{2}{c@{\hspace{0mm}}}{Algorithm \ref{alg:altmin}} & & $n$ & \multicolumn{2}{c@{\hspace{0mm}}}{Algorithm \ref{alg:altmin}} \\
\cmidrule(l){2-3} \cmidrule(l){4-5} \cmidrule(l){8-9} & Relative MSE & Time (s) & Relative MSE  & Time (s) & & & Relative MSE  & Time (s) \\\midrule
$50$ & $0.429$ & $2.28$ & $0.438$ & $17.28$ & & $350$ & $0.058$ & $6,970$\\
$100$ & $0.138$ & $47.20$ & $0.139$ & $79.01$
 & & $400$ & $0.056$ & $8,096$ \\
$150$ & $0.082$ & $336.1$ & $0.081$ & $228.7$
 & & $450$ & $0.055$ & $26,350$ \\
 $200$ & $0.0675$ & $1,906$ & $0.067$ & $841.7$ & & $500$ & $0.054$ & $28,920$ \\
 $250$ & - & - & $0.062$ & $1,419$ & & $550$ & $0.0536$ & $39,060$\\
 $300$ & - & - & $0.059$ & $2,897$ & & $600$ & $0.0525$ & $38,470$\\
\bottomrule
\end{tabular}
\label{tab:comparison}
\end{table}

Our results demonstrate the efficiency of Algorithm \ref{alg:altmin}: the relative MSE is comparable to \verb|Mosek|'s, but computational time does not explode with $n$. Since it does not require solving any SDOs and avoids the computational burden of performing the Newton step in an IPM, Algorithm \ref{alg:altmin} scales beyond $n=600$ ($1,440,000$ decision variables), compared to $n=200$ for IPMs ($80,000$ decision variables).

\subsection{Numerical Evaluation of Greedy Rounding on Matrix Completion Problems}\label{ssec:compwithheurmethods}
In this section, we compare the greedy rounding method with state-of-the-art heuristic methods, and demonstrate that, by combining greedy rounding with the local search heuristic of \citep{burer2003nonlinear}, our approach outperforms state-of-the-art heuristic methods and therefore should be considered as a viable and efficient warm-start for Algorithm \ref{alg:cuttingPlaneMethod}.

We consider the previous matrix completion problems and assess the ability to recover the low-rank matrix $\bm{M}$ (up to a relative MSE of $1\%$), for varying fraction of observed entries $p$ and {\color{black} rank $r$, with $n=100$ fixed. Note that, other than the inclusion of a Frobenius regularization term, this is the same experimental setup considered by \cite{candes2009exact, recht2010guaranteed} among others.}

We compare the performance of four methods: the greedy rounding method, both with and without the local improvement heuristic from \cite{burer2003nonlinear}, against the local improvement heuristic alone (with a thresholded-SVD initialization point) and the nuclear norm approach. Specifically, the greedy rounding method takes the solution of the previous convex relaxation with $\gamma=\frac{500}{p}$ and rounds its singular values to generate a feasible solution $\bm{Y}_{greedy}$. For the local improvement heuristic, we solve:
\begin{align*}
    \min_{\bm{X} \in \mathbb{R}^{n \times n}, \bm{U}, \bm{V} \in \mathbb{R}^{n \times k}} \ \frac{1}{2 \gamma}\Vert \bm{X}\Vert_2^2+\sum_{(i,j) \in \Omega}(X_{i,j}-A_{i,j})^2 \quad \text{s.t.}\quad \bm{X}=\bm{U}\bm{V}^\top,
\end{align*}
for $\gamma=\frac{500}{p}$ and $k=r$, and iteratively optimize over $\bm{U}$ and $\bm{V}$ using \verb|Mosek|. We provide an initial value for $\bm{U}$ by either taking the first $k$ left-singular vectors of a matrix $\bm{A}$ where unobserved entries are replaced by 0, or taking the square root of $\bm{Y}_{greedy}$. For the nuclear norm regularization strategy, since our observations are noiseless, we solve:
$
    \min_{\bm{X} \in \mathbb{R}^{n \times n}} \Vert \bm{X}\Vert_* \ \text{s.t.}\ X_{i,j}=A_{i,j} \quad \forall (i,j) \in \Omega.
$

\begin{figure}
    \centering
\begin{subfigure}[t]{.42\linewidth}
    	\centering
    	\includegraphics[width=.9\linewidth]{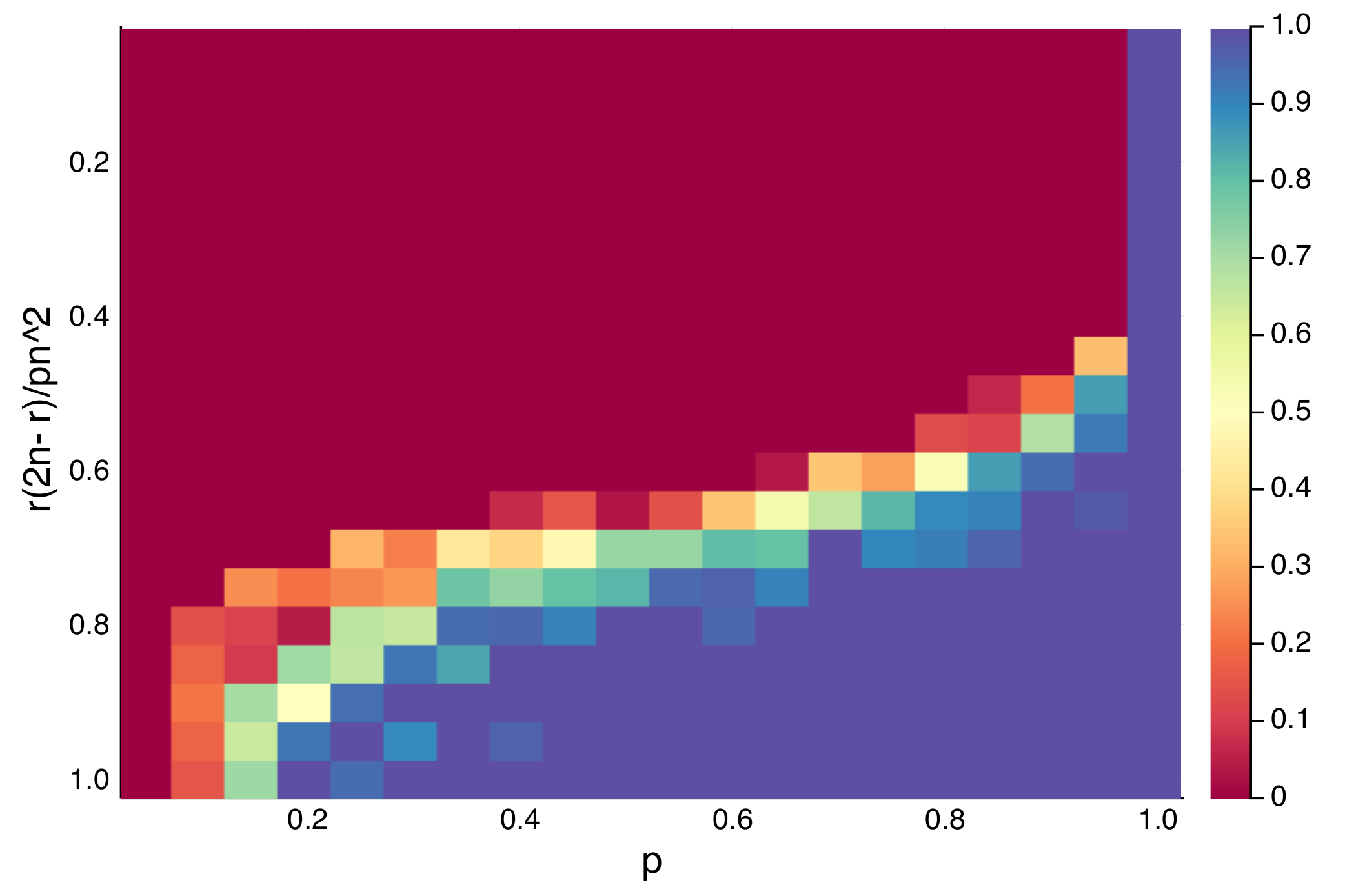}
        \caption{\footnotesize Greedy Rounding}
    \end{subfigure} %
    \begin{subfigure}[t]{.42\linewidth}
    	\centering
    	\includegraphics[width=.9\linewidth]{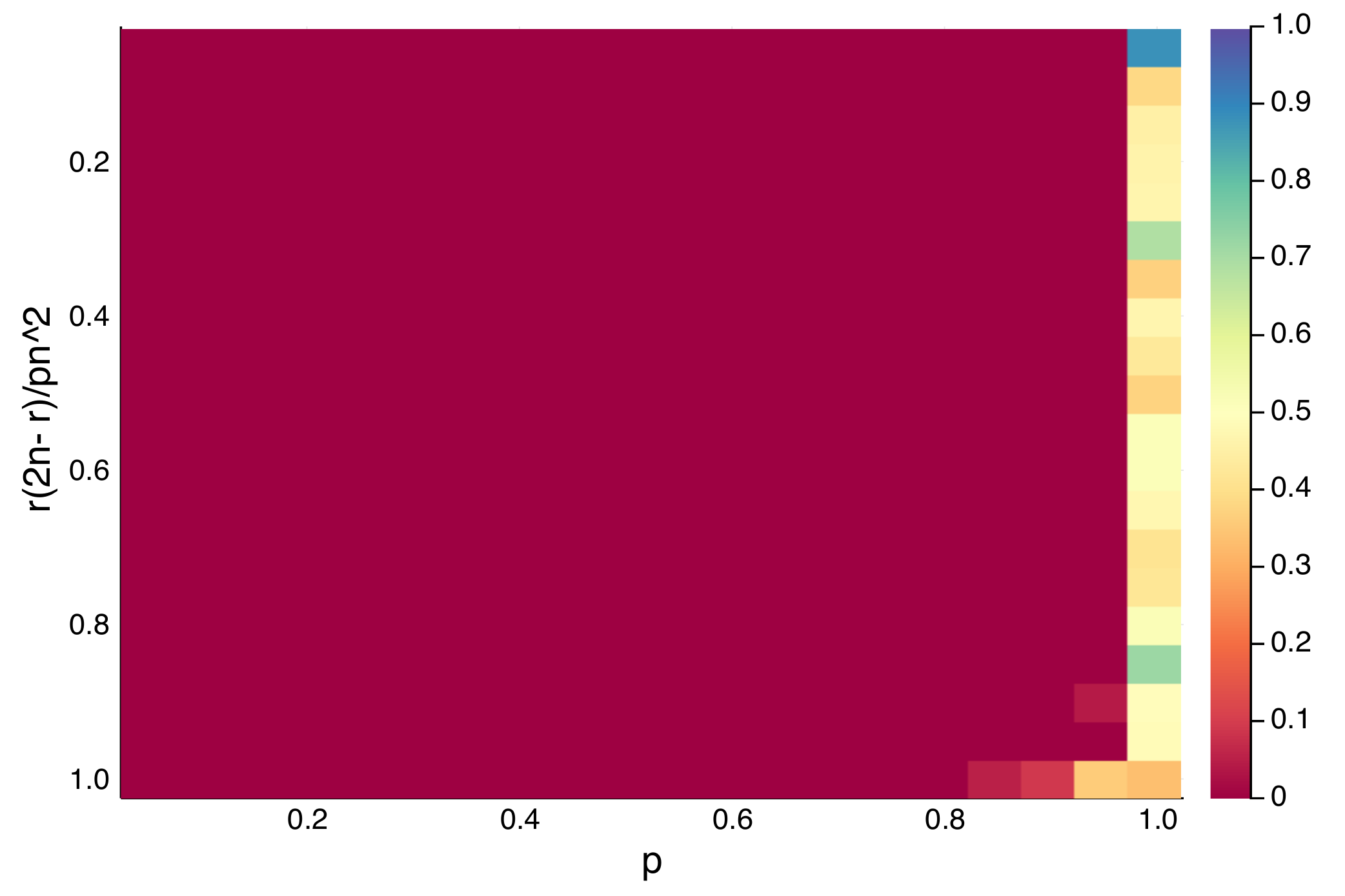}
    \caption{\footnotesize Nuclear Norm}
\end{subfigure}\\
\begin{subfigure}[t]{.42\linewidth}
    	\centering
    	\includegraphics[width=.9\linewidth]{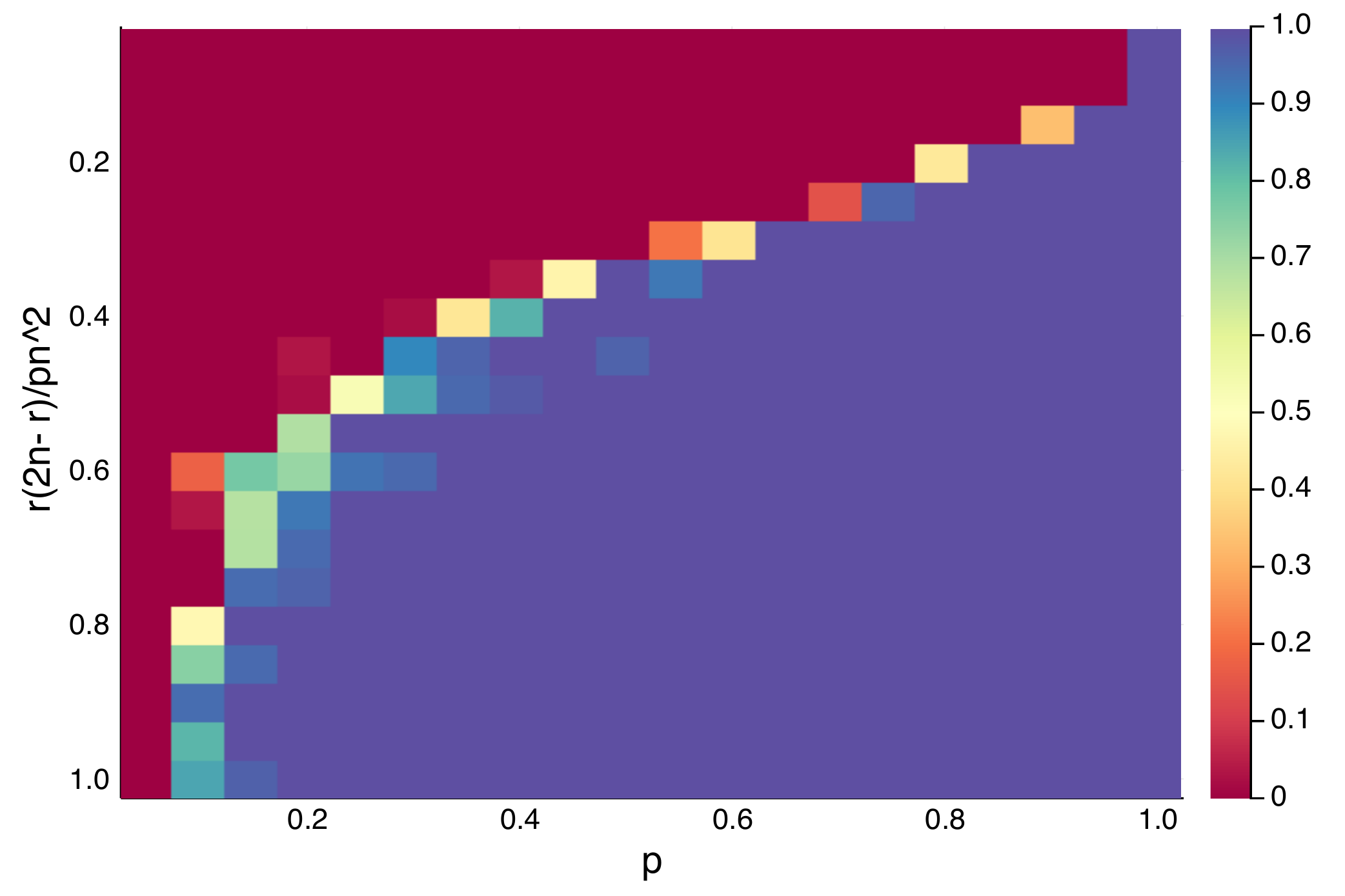}
        \caption{\footnotesize SVD+Local Improvement}
    \end{subfigure} %
    \begin{subfigure}[t]{.42\linewidth}
    	\centering
    	\includegraphics[width=.9\linewidth]{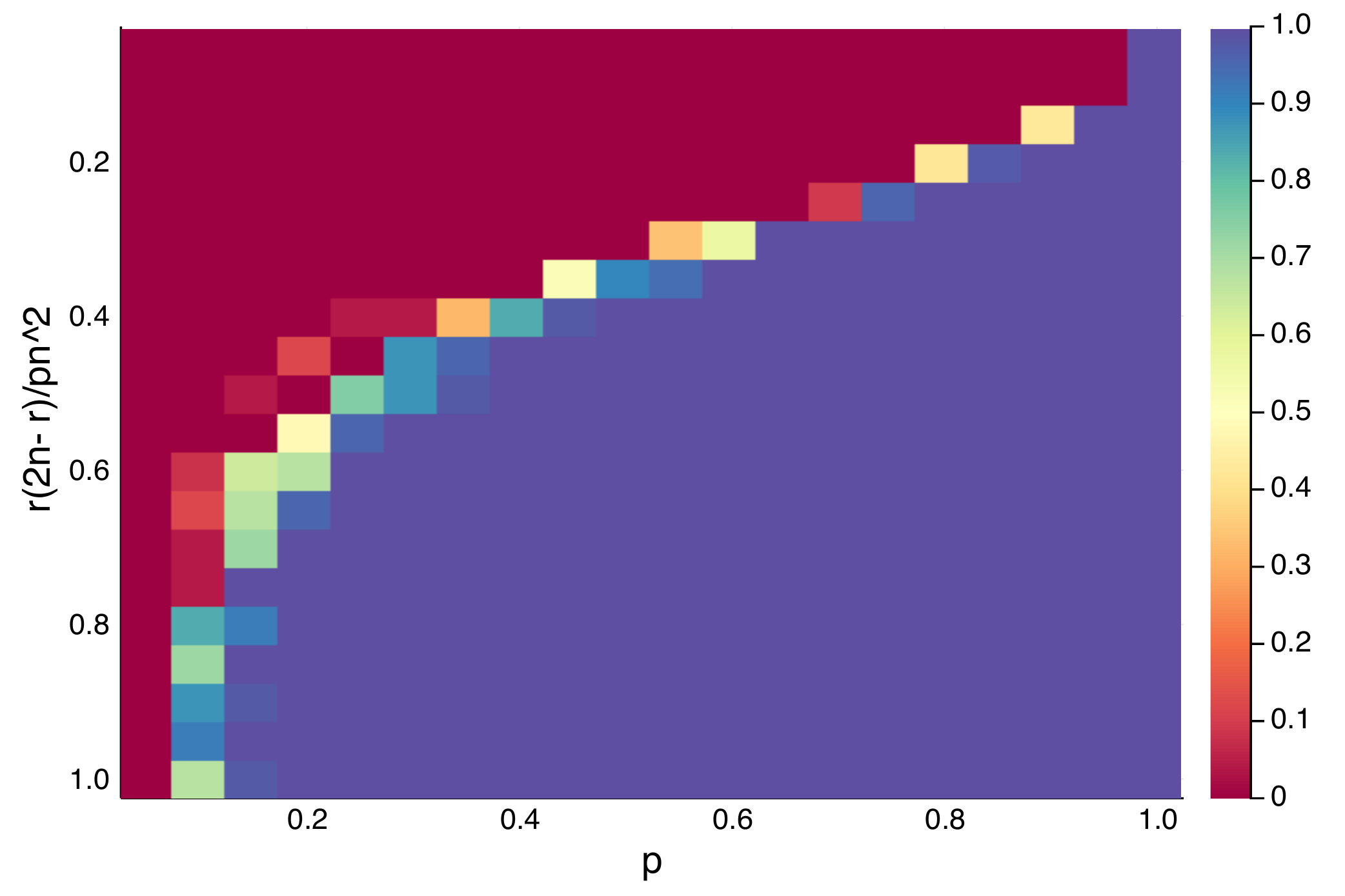}
    \caption{\footnotesize Greedy+Local Improvement}
    \end{subfigure}
\caption{Prop. matrices recovered with $\leq 1\%$ relative MSE (higher is better), for different values of $p$ (x-axis) and ${r(2n-r)}{/pn^2} \propto 1/n$ (y-axis), averaged over $25$ rank-$r$ matrices.}
\label{fig:proprec}
\end{figure}

\begin{figure}
    \centering
\begin{subfigure}[t]{.42\linewidth}
    	\centering
    	\includegraphics[width=.9\linewidth]{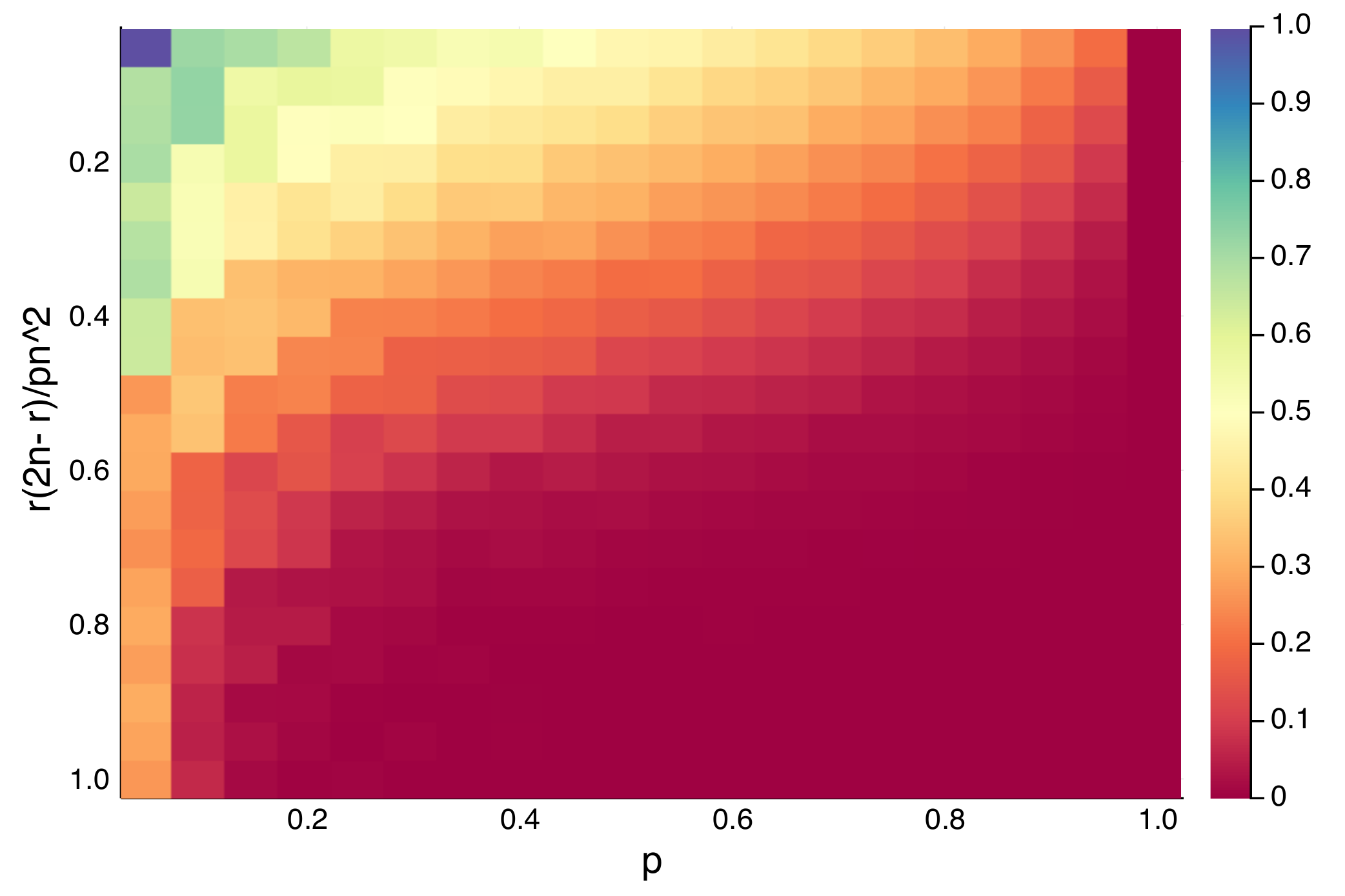}
        \caption{\footnotesize Greedy Rounding}
    \end{subfigure} %
    \begin{subfigure}[t]{.42\linewidth}
    	\centering
    	\includegraphics[width=.9\linewidth]{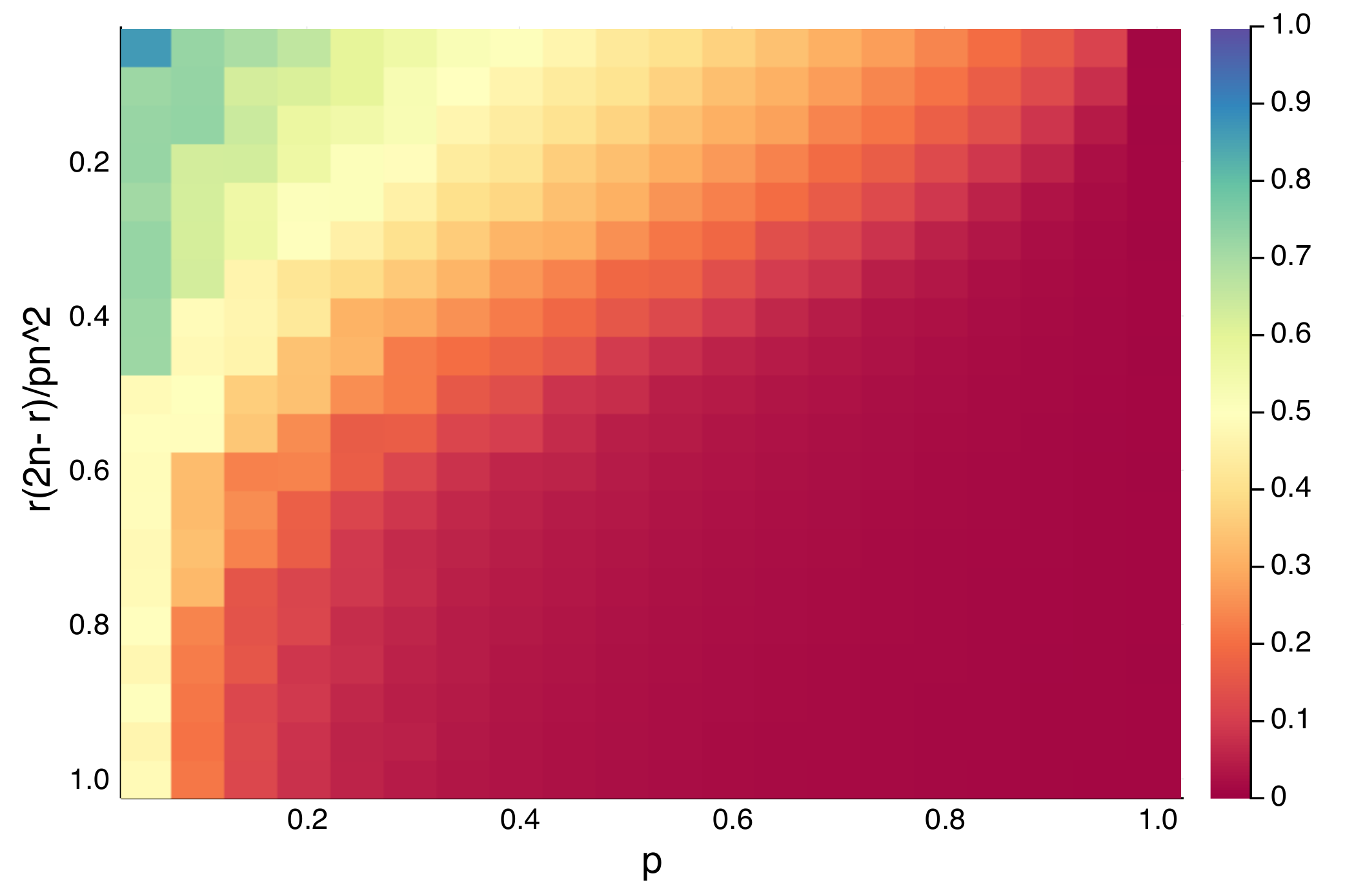}
    \caption{\footnotesize Nuclear Norm}
\end{subfigure} \\
\begin{subfigure}[t]{.42\linewidth}
    	\centering
    	\includegraphics[width=.9\linewidth]{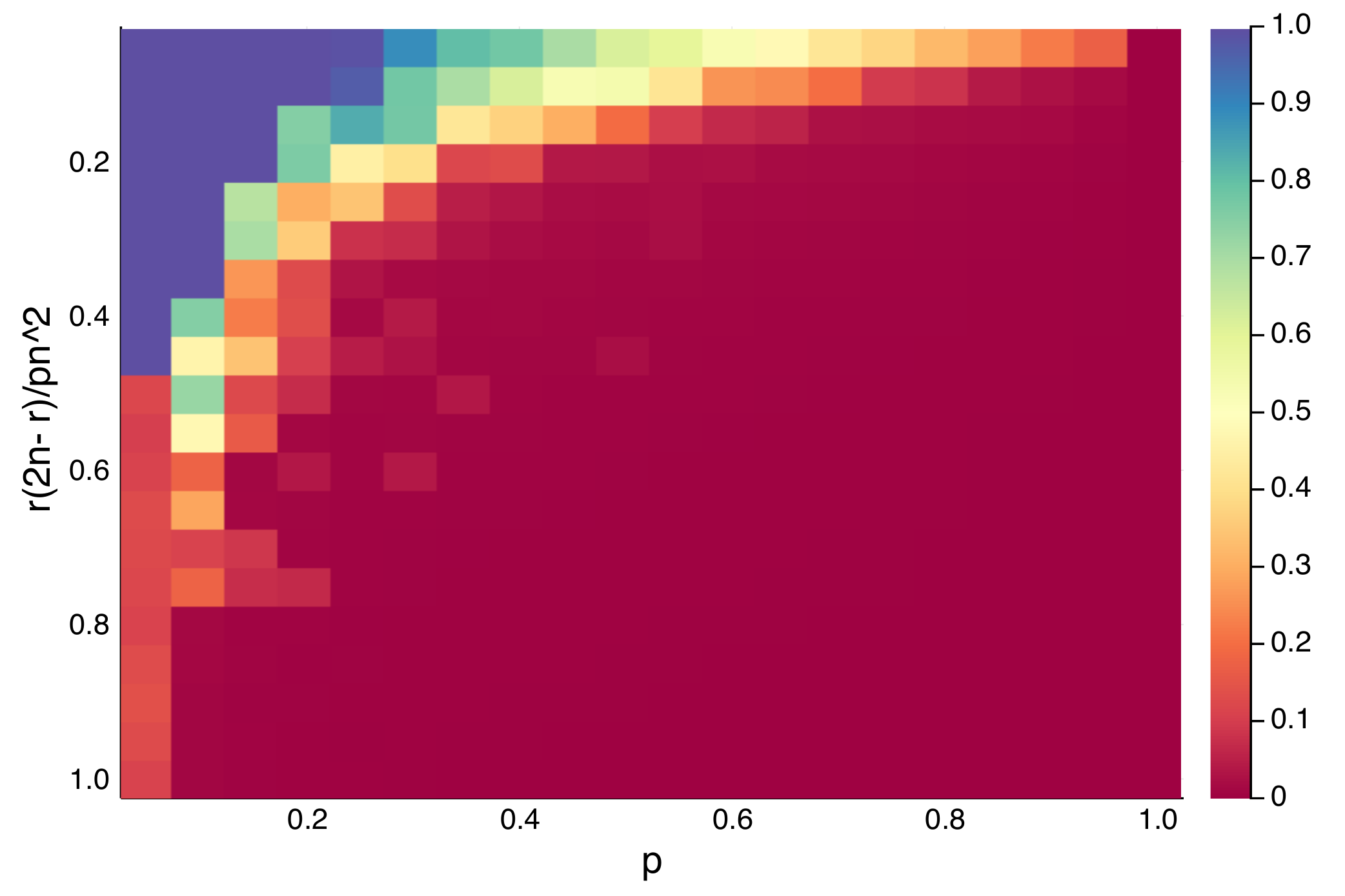}
        \caption{\footnotesize SVD+Local Improvement}
    \end{subfigure} %
    \begin{subfigure}[t]{.42\linewidth}
    	\centering
    	\includegraphics[width=.9\linewidth]{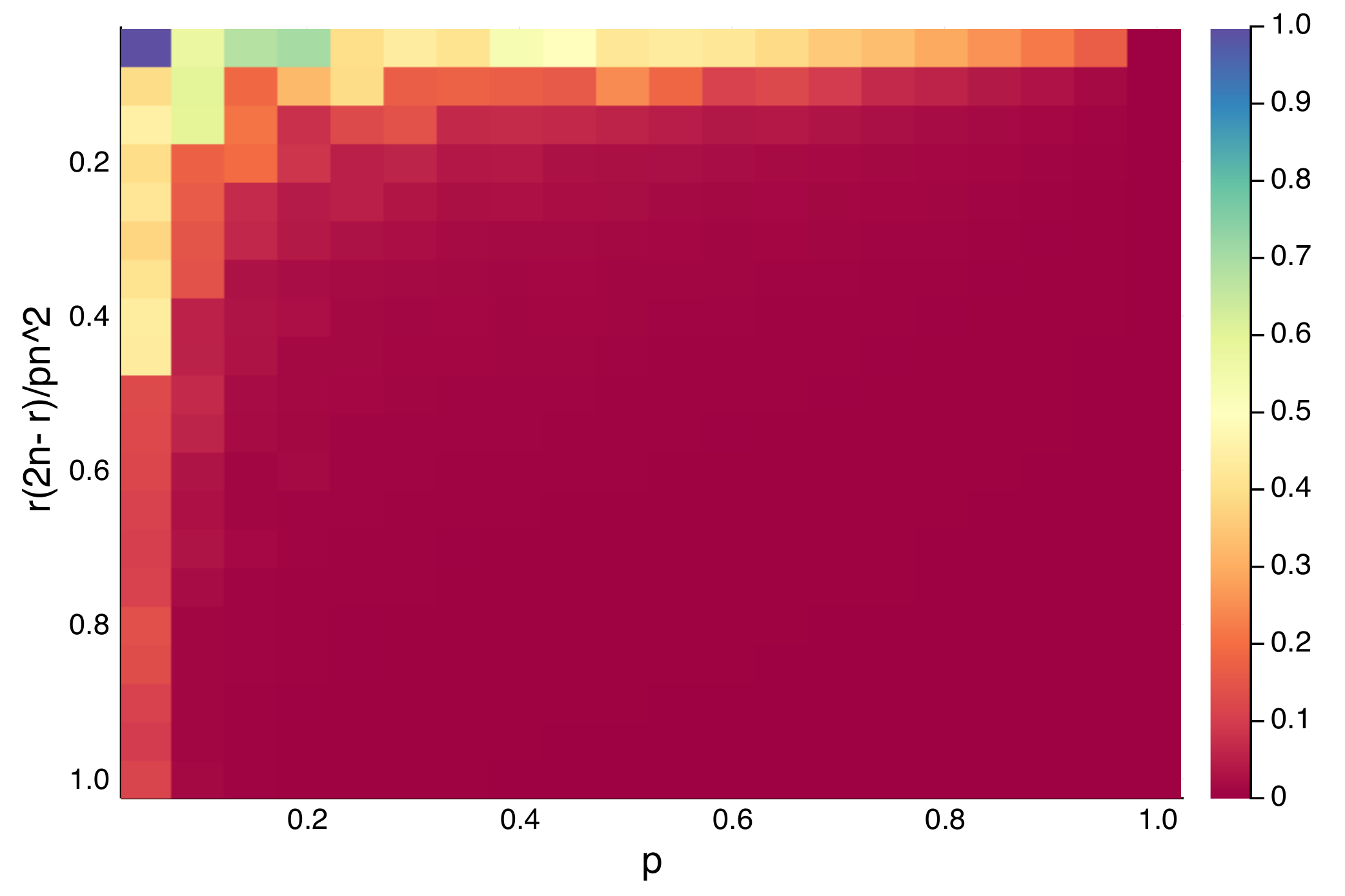}
    \caption{\footnotesize Greedy+Local Improvement}
    \end{subfigure}
\caption{Average relative MSE (lower is better), averaged over $25$ rank-$r$ matrices. We cap the relative MSE at $1.0$.}
\label{fig:mse_rr20}
\end{figure}

Figure \ref{fig:proprec} depicts the proportion of times the matrix was recovered exactly (averaged over $25$ samples per {\color{black}tuple} of $(n,p,r)$), while Figure \ref{fig:mse_rr20} depicts the relative average MSE over all instances. {\color{black}As in \cite{candes2009exact, recht2010guaranteed}, we vary $p$ between $0$ and $1$ and consider all possible ranks $r$ such that $r(2n-r) \leq pn^2$. 
}
From this set of experiments, we make several observations: First, greedy rounding and the local improvement heuristic outperform nuclear norm minimization both in terms of average relative MSE and amount of data required to recover the matrix. Second, the local improvement heuristic improves upon greedy rounding. In terms of its ability to recover the underlying matrix exactly, it performs equally well with either initialization strategy. However, initialization with the greedy rounding supplies dramatically lower average MSEs in instances where no approach recovers the true matrix exactly. This suggests that initialization strategies for the Burer-Monterio heuristic should be revisited and greedy rounding considered as a viable and more accurate alternative than selecting a random feasible point.

\subsection{Benchmarking Algorithm \ref{alg:cuttingPlaneMethod} on Synthetic Matrix Completion Problems}\label{ssec:benchmarkalg1}
We now benchmark Algorithm \ref{alg:cuttingPlaneMethod} on matrix completion problems {\color{black}where $n \in \{10, 20, 30\}$}.

We first compare the two different implementations of Algorithm  \ref{alg:cuttingPlaneMethod}, single- and multi-tree, with solving the problem directly as a {\color{black}QCQO} in \verb|Gurobi| $9.0$ (Section \ref{ssec:qcqopgurobi}). In Algorithm  \ref{alg:cuttingPlaneMethod}, the lower bounds are warm-started with $200$ cuts from the in-out method, and greedy rounding with local search improvement is used for the upper bounds{\color{black}; if a single-tree instance fails to find a feasible solution (due to numerical instability in \verb|Gurobi|) we return the gap between the warm-start and the semidefinite relaxation.} At the $t$th iteration, we impose a time limit of $10t$ seconds for generating the new cut so as to increase numerical precision as the solver progresses. We also impose a limit of $20$ cuts for the multi-tree approach, a time limit of {\color{black}$30,000$s} for the single-tree approach{\footnote{\color{black}We require a larger time limit than $3,600$s, since Gurobi often fails to find any feasible solutions within this time limit due to the numerical difficulties inherent in integrating lazy constraint callbacks and a non-convex master problem.}}, {\color{black}a time limit of $3,600$s for} \verb|Gurobi|, and an optimality gap of $1\%$\footnote{We report the {\color{black}absolute} gap between the better of Gurobi's lower bound and the semidefinite lower bound, compare to the objective value which we evaluate directly; this is sometimes $1-2\%$ even when Gurobi reports that it has found an optimal solution, due to numerical instability in Gurobi. {\color{black}Note that we report the absolute, rather than relative, gap since the relative gap depends on the quality of Gurobi's approximation of $\mathcal{Y}^k_n$, which is controlled by the parameter FuncPieceError and cannot be set lower than $10^{-6}$;
also note that the objective values are on the order of $0.5$-$5.0$ for the problems reported in Table $4$.}}. Average runtime, number of nodes, and optimality gap are reported in Table \ref{tab:comparison2}. {\color{black}Note that the same random instances were solved by all three approaches (by fixing the random seeds), to facilitate a less noisy comparison}.

We observe that Algorithm  \ref{alg:cuttingPlaneMethod} drastically improves upon \verb|Gurobi| both in terms of computational time (reduced by up to an order of magnitude) and accuracy ({\color{black}absolute gap reduced by around an order of magnitude}). {\color{black} Multi-tree dominates single-tree and Gurobi in terms of runtime and the quality of the solution found, although single-tree occasionally has a smaller gap at termination. Moreover, multi-tree consistently finds high-quality feasible solutions earlier than single tree and accepts our warm-start more consistently, which suggests it may scale better to high-dimensional settings. }

\begin{table}[h]
\footnotesize
\caption{{\color{black}Scalability of Algorithm \ref{alg:cuttingPlaneMethod} vs. Gurobi for solving rank-$1$ matrix completion problems to certifiable optimality, averaged over $20$ random matrices per row. In multi-tree, Nodes denotes the number of nodes expanded in the final branch-and-cut tree, while Nodes (t) denotes the number of nodes expanded over all trees for the multi-tree implementation.}}
\begin{tabular}{@{}l l l r r r r r r r r r r r r r r r  @{}} \toprule
  &  & & \multicolumn{3}{c@{\hspace{0mm}}}{Gurobi (direct)} &\multicolumn{4}{c@{\hspace{0mm}}}{Algorithm \ref{alg:cuttingPlaneMethod} (single-tree)} &\multicolumn{4}{c@{\hspace{0mm}}}{Algorithm \ref{alg:cuttingPlaneMethod} (multi-tree)}  \\
\cmidrule(l){4-6} \cmidrule(l){7-10} \cmidrule(l){11-14}
$n$ & $p$ & $\gamma$ & Time(s) & Nodes & Gap & Time(s) & Nodes & Gap & Cuts & Time(s) & Nodes & Nodes (t)& Gap & Cuts \\\midrule
$10$ & $0.1$ & $20/p$ & {\color{black}$>3,600$} & {\color{black}$313,700$} & {\color{black}$0.0301$} & {\color{black}$10,310$} & {\color{black}$40,060$} & {\color{black}$0.0004$} & {\color{black}$23,460$} & {\color{black}$252.3$} & {\color{black}$630.2$} & {\color{black}$10,099$} & {\color{black}$0.0019$} & {\color{black}$2.95$}\\
$10$ & $0.2$ & $20/p$ & {\color{black}$>3,600$} & {\color{black}$299,200$} & {\color{black}$0.0854$} & {\color{black}$19,440$} & {\color{black}$28,430$} & {\color{black}$0.0229$} & {\color{black}$19,370$} & {\color{black}$1,672$} & {\color{black}$2,277$} & {\color{black}$28,895$} & {\color{black}$0.0104$} & {\color{black}$11.0$}\\
$10$ & $0.3$ & $20/p$ & {\color{black}$>3,600$} & {\color{black}$274,500$} & {\color{black}$0.1167$} & {\color{black}$20,368$} & {\color{black}$25,480$} & {\color{black}$0.0433$} & {\color{black}$20,290$} & {\color{black}$2,319$} & {\color{black}$2,684$} & {\color{black}$30,906$} & {\color{black}$0.0317$} & {\color{black}$15.4$}\\
\midrule
 $10$ & $0.1$ & $100/p$ & {\color{black}$>3,600$} & {\color{black}$281,800$} & {\color{black}$0.0068$} & {\color{black}$18,580$} & {\color{black}$62,680$} & {\color{black}$0.0015$} & {\color{black}$42,200$} & {\color{black}$239.5$} & {\color{black}$405$} & {\color{black}$7,499$} & {\color{black}$0.0003$} & {\color{black}$3.20$}\\
$10$ & $0.2$ & $100/p$ & {\color{black}$>3,600$} & {\color{black}$271,300$} & {\color{black}$0.0178$} & {\color{black}$27,990$} & {\color{black}$39,330$} & {\color{black}$0.0492$} & {\color{black}$31,060$} & {\color{black}$1,269$} & {\color{black}$1,931$} & {\color{black}$21,420$} & {\color{black}$0.0042$} & {\color{black}$8.40$}\\
$10$ & $0.3$ & $100/p$ & {\color{black}$3,239$} & {\color{black}$237,000$} & {\color{black}$0.0178$} & {\color{black}$25,750$} & {\color{black}$29,350$} & {\color{black}$0.0434$} & {\color{black}$23,390$} & {\color{black}$2,472$} & {\color{black}$2,134$} & {\color{black}$32,196$} & {\color{black}$0.0098$} & {\color{black}$19.6$}\\
\midrule
$20$ & $0.1$ & $20/p$ & {\color{black}$>3,600$} & {\color{black}$80,760$} & {\color{black}$0.8915$} & {\color{black}$>30,000$} & {\color{black}$13,110$} & {\color{black}$0.741$} & {\color{black}$13,070$} & {\color{black}$2,917$} & {\color{black}$413.2$} & {\color{black}$4,835$} & {\color{black}$0.0166$} &  {\color{black}$18.5$}\\
$20$ & $0.2$ & $20/p$ & {\color{black}$>3,600$} & {\color{black}$65,310$} & {\color{black}$4.094$} & {\color{black}$>30,000$} & {\color{black}$7,023$} & {\color{black}$0.1816$} & {\color{black}$7,008$} & {\color{black}$3,512$} & {\color{black}$143.0$} & {\color{black}$1,746$} & {\color{black}$0.247$} & {\color{black}$20.0$} \\
$20$ & $0.3$ & $20/p$ & {\color{black}$>3,600$} & {\color{black}$64,850$} & {\color{black}$4.745$} & {\color{black}$28,700$} & {\color{black}$6,914$} & {\color{black}$0.1066$} & {\color{black}$6,828$} & {\color{black}$3,287$} & {\color{black}$204.9$} & {\color{black}$2,081$} & {\color{black}$0.253$} & {\color{black} $19.6$}\\
\midrule
$20$ & $0.1$ & $100/p$ & {\color{black}$>3,600$} & {\color{black}$60,830$} & {\color{black}$0.428$} & {\color{black}$>30,000$} & {\color{black}$13,790$} & {\color{black}$0.7714$} & {\color{black}$13,799$} & {\color{black}$3,013$} & {\color{black}$508.0$} & {\color{black}$6,152$} & {\color{black}$0.0072$} & {\color{black}$17.6$}\\
$20$ & $0.2$ & $100/p$ & {\color{black}$>3,600$} & {\color{black}$43,850$} & {\color{black}$1.421$} & {\color{black}$>30,000$} & {\color{black}$6,395$} & {\color{black}$0.0543$} & {\color{black}$6,395$} & {\color{black}$3,106$} & {\color{black}$80.9$} & {\color{black}$1,027$} & {\color{black}$0.0903$} & {\color{black}$17.6$}\\
$20$ & $0.3$ & $100/p$ & {\color{black}$>3,600$} & {\color{black}$55,150$} & {\color{black}$2.810$} & {\color{black}$29,530$} & {\color{black}$6,538$} & {\color{black}$0.0271$} & {\color{black}$6,510$} & {\color{black}$2,910$} & {\color{black}$62.7$} & {\color{black}$744.2$} & {\color{black}$0.1368$} & {\color{black}$17.0$}\\
\bottomrule
\end{tabular}
\label{tab:comparison2}
\end{table}

{Next, we evaluate the performance of the multi-tree implementation of Algorithm \ref{alg:cuttingPlaneMethod} on a more extensive test-set, including instances where $\mathrm{Rank}(\bm{M})>1$, in Table \ref{tab:comparison3}. 
Note that when $r=1$ we use the same experimental setup (although we impose a time limit of $30t$ seconds, or $7200$ seconds if there has been no improvement for two consecutive iterations, a cut limit of $50$ cuts when $n>20$), and when $r>1$ we increase the time limit per iteration to $300$t seconds (or $7200$ seconds if there has been no improvement for two consecutive iterations), and allow up to $100$ PSD cuts per iteration to be added at the root node via a user cut callback, in order to strengthen the approximation of the PSD constraint $\bm{Y}\succeq \bm{0}$. We observe that the problem's complexity increases with the rank, although not too excessively. Moreover, when $r>1$ the bound gap is actually smaller when $\gamma=\frac{100}{p}$ than when $\gamma=\frac{20}{p}$. We believe this is because Gurobi cannot represent the SDO constraint $\bm{Y} \succeq \bm{0}$ and its SOC approximation is inexact (even with PSD cuts), and in some cases refining this approximation is actually harder than refining our approximation of $g(\bm{X})$.

Note that the main bottleneck inhibiting solving matrix completion problems {\color{black}where $n \geq 50$} is \verb|Gurobi| itself, as the non-convex solver takes increasing amounts of time to process warm-starts (sometimes in the $100$s or $1000$s of seconds) when $n$ increases. We believe this may be because of the way \verb|Gurobi| translates orthogonal projection matrices to a piecewise linear formulation. Encouragingly, this suggests that our approach may successfully scale to $100 \times 100$ matrices as \verb|Gurobi| improves their solver.

\begin{table}
\footnotesize
\caption{Scalability of Algorithm \ref{alg:cuttingPlaneMethod} (multi-tree) for solving low-rank matrix completion problems to certifiable optimality, averaged over $20$ random matrices per row.}
\begin{tabular}{@{}l l l r r l r r r r r r r r r r r  r@{}} \toprule
  &  & & \multicolumn{4}{c@{\hspace{0mm}}}{Rank-$1$} &\multicolumn{4}{c@{\hspace{0mm}}}{Rank-$2$} &\multicolumn{4}{c@{\hspace{0mm}}}{Rank-$3$}  \\
\cmidrule(l){4-7} \cmidrule(l){8-11} \cmidrule(l){12-15}
$n$ & $p$ & $\gamma$ & Time(s) & Nodes & Gap & Cuts & Time(s) & Nodes & Gap & Cuts & Time(s) & Nodes & Gap  & Cuts \\\midrule
$10$ & $0.1$ & $20/p$ & $182.1$ & $9,755$ & $0.0005$ & $2.56$ & $24,220$ & $35,670$ & $0.0034$ & $5.78$ & $37,780$ & $39,870$ & $0.0071$ & $9.28$\\
$10$ & $0.2$ & $20/p$ & $3,508$ & $21,060$ & $0.0026$ & $10.8$ & $209,900$ & $108,000$ & $0.0252$ & $35.3$ & $135,260$ & $35,870$ & $0.031$ & $26.2$\\
$10$ & $0.3$ & $20/p$ & $5,488$ & $30,970$ & $0.0039$ & $13.1$ & $302,200$ & $70,500$ & $0.0866$ & $50.0$ & $302,100$ & $31,870$ & $0.0197$ & $50.0$\\
\midrule
$10$ & $0.1$ & $100/p$ & $656.5$ & $28,870$ & $0.0001$ & $2.14$ & $676.1$ & $25,493$ & $0.0009$ & $1.83$ & $842.7$ & $20,700$ & $0.0024$ & $1.79$\\
$10$ & $0.2$ & $100/p$ & $1,107$ & $10,010$ & $0.0009$ & $4.29$ & $2,065$ & $42,490$ & $0.0019$ & $5.61$ & $57,530$ & $36,910$ & $0.0124$ & $10.7$\\
$10$ & $0.3$ & $100/p$ & $3,364$ & $48,730$ & $0.0022$ & $6.30$ & $272,300$ & $33,150$ & $0.0195$ & $44.7$ & $249,700$ & $35,530$ & $0.0499$ & $42.2$\\
\midrule
$20$ & $0.1$ & $20/p$ & $2,017$	 & $4,756$ &	$0.0061$ & $8.20$ & $253,900$	& $8,030$	& $0.0279$	& $42.7$ & $255,400$ & $3,015$ & $0.0309$ & $43.2$ \\
$20$ & $0.2$ & $20/p$ & $6,369$ & $6,636$ & $0.0136$ & $15.0$ & $298,700$	& $3,342$	& $0.549$	& $50.0$ & $295,500$ & $236.5$ & $0.879$	& $50.0$\\
$20$ & $0.3$ & $20/p$ & $6,687$ & $4,187$ & $0.0082$ & $18.4$ & $296,500$ & $3,175$ & $1.123$ & $50.0$ & $291,100$ & $41.35$ &	$2.147$ & $50.0$ \\
\midrule
$20$ & $0.1$ & $100/p$ & $1,266$ & $8,792$ & $0.0087$ & $8.35$ & $211,700$	& $6,860$ &	$0.0073$ &	$34.24$ & $171,900$ &	$2,350$	& $0.0131$ & $29.8$\\
$20$ & $0.2$ & $100/p$ & $1,220$ & $2,710$ & $0.0104$ & $7.80$ & $302,800$ & $2,426$ & $0.123$ & $50.0$ & $298,800$ & $221.4$ & $0.123$ & $50.0$\\
$20$ & $0.3$ & $100/p$ & $1,272$ & $1,837$ & $0.0064$ & $3.14$ & $299,000$	& $2,518$	& $0.264$ & $50.0$ & $293,500$ & $43.0$ & $0.659$ & $50.0$\\
\midrule
$30$ & $0.1$ & $20/p$ & $300,300$ & $2,735$ & $0.0905$ & $50.0$  & $304,300$ & $164.0$ & $0.790$ & $50.0$ & $303,100$ & $1.10$ & $0.365$ & $50.0$\\
$30$ & $0.2$ & $20/p$ & $298,700$ & $1,511$ &	$0.136$ & $50.0$ & $301,700$ & $9.62$ & $3.105$ & $50.0$ & $302,600$ & $1.00$ & $5.581$ & $50.0$\\
$30$ & $0.3$ & $20/p$ & $183,800$ &	$1,743$	& $0.0476$	& $36.9$ & $303,000$ & $1.63$ & $5.232$ & $50.0$ & $305,000$ & $0.70$ & $14.60$	& $50.0$\\
\midrule
$30$ & $0.1$ & $100/p$ & $305,600$ & $2,262$ & $0.0273$ & $50.0$ & $302,800$ & $97.40$ & $0.0973$ & $50.0$ & $305,000$ & $1.90$ & $0.0967$ & $50.0$\\
$30$ & $0.2$ & $100/p$ & $246,300$ & $3,285$ & $0.0315$ & $43.6$ & $304,300$ & $6.17$ & $0.697$ & $50.0$ & $302,600$ & $1.00$ & $1.419$ & $50.0$\\
$30$ & $0.3$ & $100/p$ & $25,970$ & $11,020$ & $0.0089$ & $17.1$ & $304,000$ & $1.00$ & $0.923$ & $50.0$ & $304,700$ & $1.00$ & $3.221$ & $50.0$\\
\bottomrule
\end{tabular}
\label{tab:comparison3}
\end{table}
}

\begin{figure}
    \centering
\begin{subfigure}[t]{.47\linewidth}
    	\centering
    	\includegraphics[width=.85\linewidth]{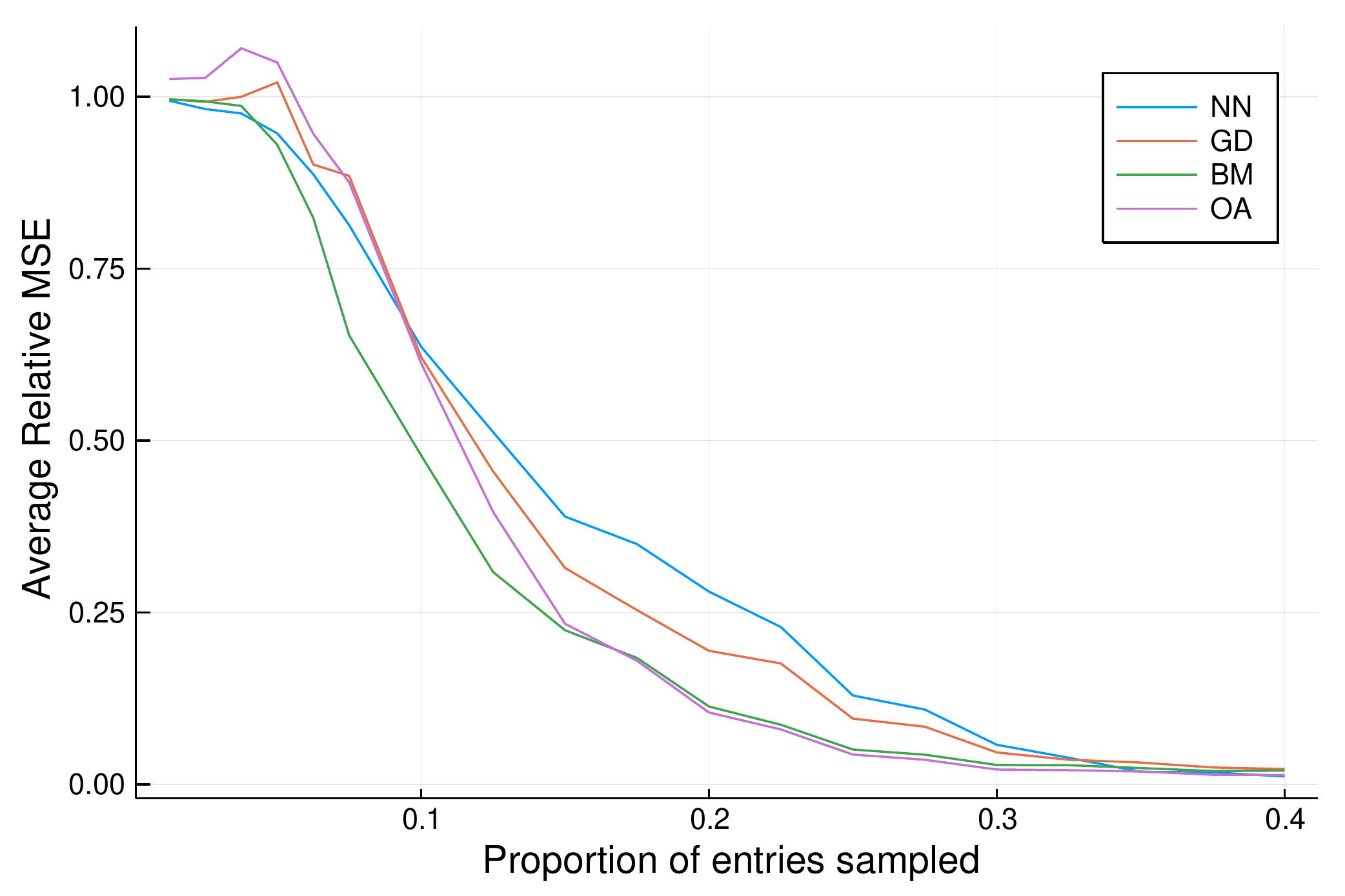}
     	\caption{\footnotesize $n=25$}
    \end{subfigure} %
    \begin{subfigure}[t]{.47\linewidth}
    	\centering
    	\includegraphics[width=.85\linewidth]{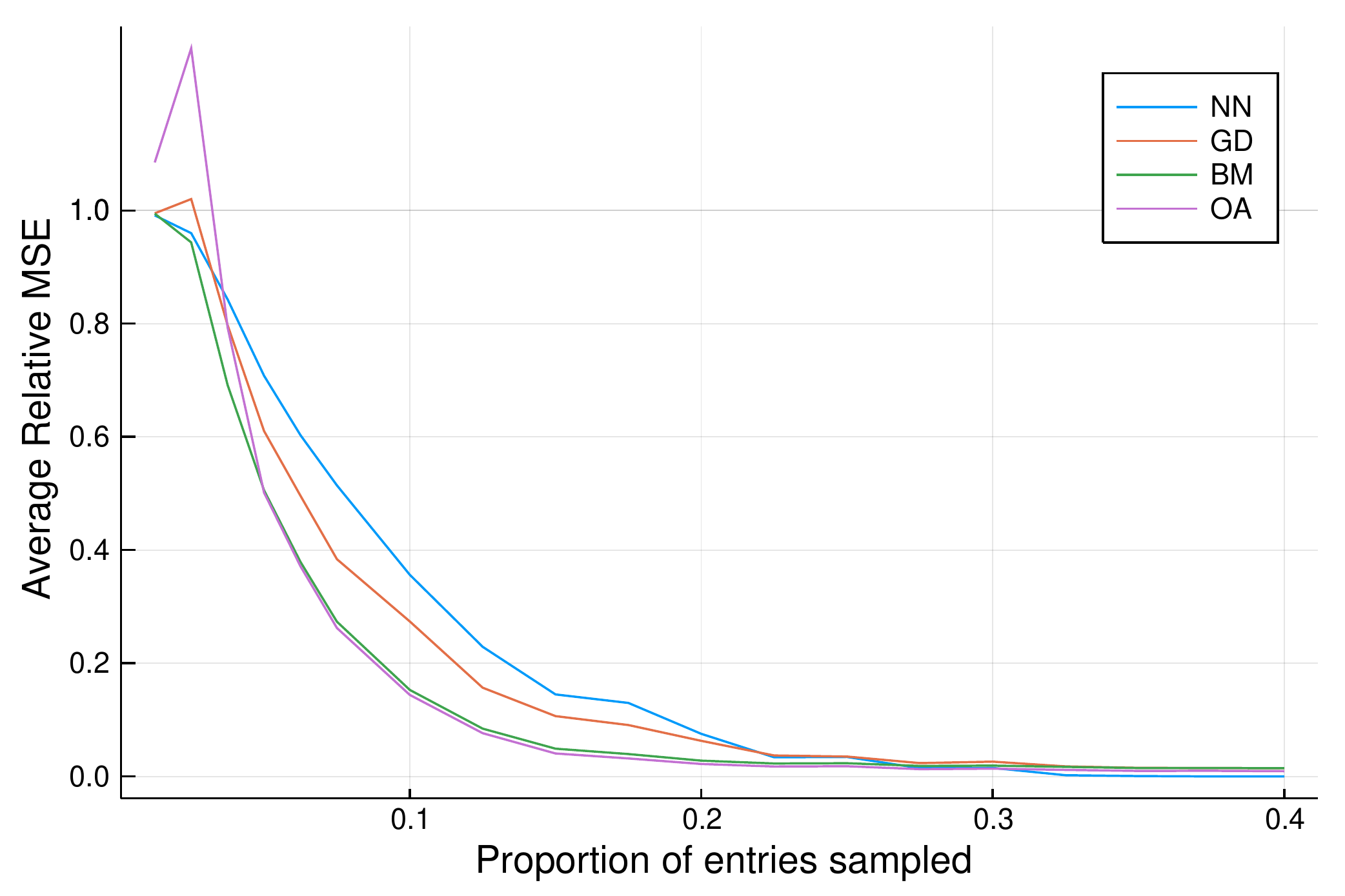}
     	\caption{\footnotesize $n=50$}
    \end{subfigure}
\caption{Average relative MSE for nuclear norm (NN), greedy rounding (GD), Burer-Monterio (BM), and outer-approximation (OA) when imputing a rank-$1$ $n \times n$ matrix. All results are averaged over $25$ matrices.}
\label{fig:mse_rr2}
\end{figure}

Finally, we compare the solution from the exact formulation \eqref{rankminproblem_proj} solved using Algorithm \ref{alg:cuttingPlaneMethod} (multi-tree) with the initial warm-start we proposed and two state-of-the-art heuristics, namely nuclear norm minimization and the Burer-Monterio approach, as in Section \ref{ssec:compwithheurmethods}. Here, we take $n \in \{25, 50\}$, $r=1$, $p$ ranging from $0$ to $0.4$, and $\gamma=\frac{100}{p}$.
Figure \ref{fig:mse_rr2} depicts the average relative MSE over the entire matrix, averaged over $25$ random instances per value of $p$. 
When $p \geq 0.2$, the exact method supplies an out-of-sample relative MSE around $0.6\%$ lower than Burer-Monterio\footnote{Because we ran all methods on the same random instances, this difference is statistically significant, with a p-value of $2\times 10^{-51}$ (resp. $2\times 10^{-129}$) that the relative MSE is lower for the exact method when $n=25$ (resp. $n=50$).}.

{\color{black}
\subsection{Exploring the Impact of Regularization on Problem Complexity}\label{ssec:impactofreg}
We now examine the impact of the regularization term $\frac{1}{2\gamma}\Vert \bm{X}\Vert_F^2$ on the problem complexity, as captured by the relative in-sample duality gap between the semidefinite relaxation and the objective value of the greedy solution with a BM local improvement heuristic. We generate the problem data in the same manner as the previous experiment, and display results for four values of $\gamma$
in Figure \ref{fig:bgapcomp}.
\begin{figure}
    \centering
\begin{subfigure}[t]{.42\linewidth}
    	\centering
    	\includegraphics[width=.85\linewidth]{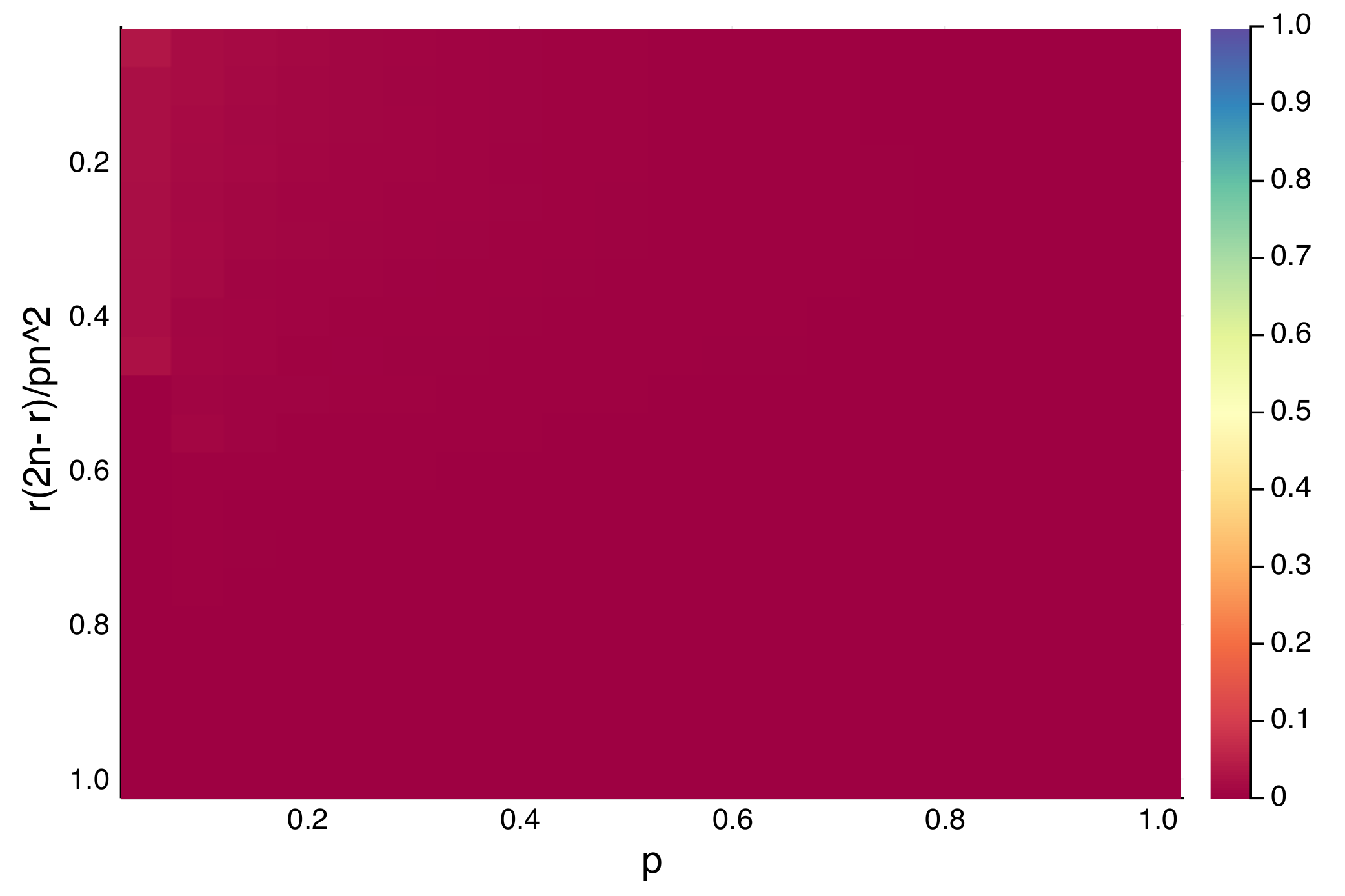}
        \caption{\footnotesize Relative in-sample bound gap ($\gamma=\frac{0.5}{p}$)}
    \end{subfigure} %
    \begin{subfigure}[t]{.42\linewidth}
    	\centering
    	\includegraphics[width=.85\linewidth]{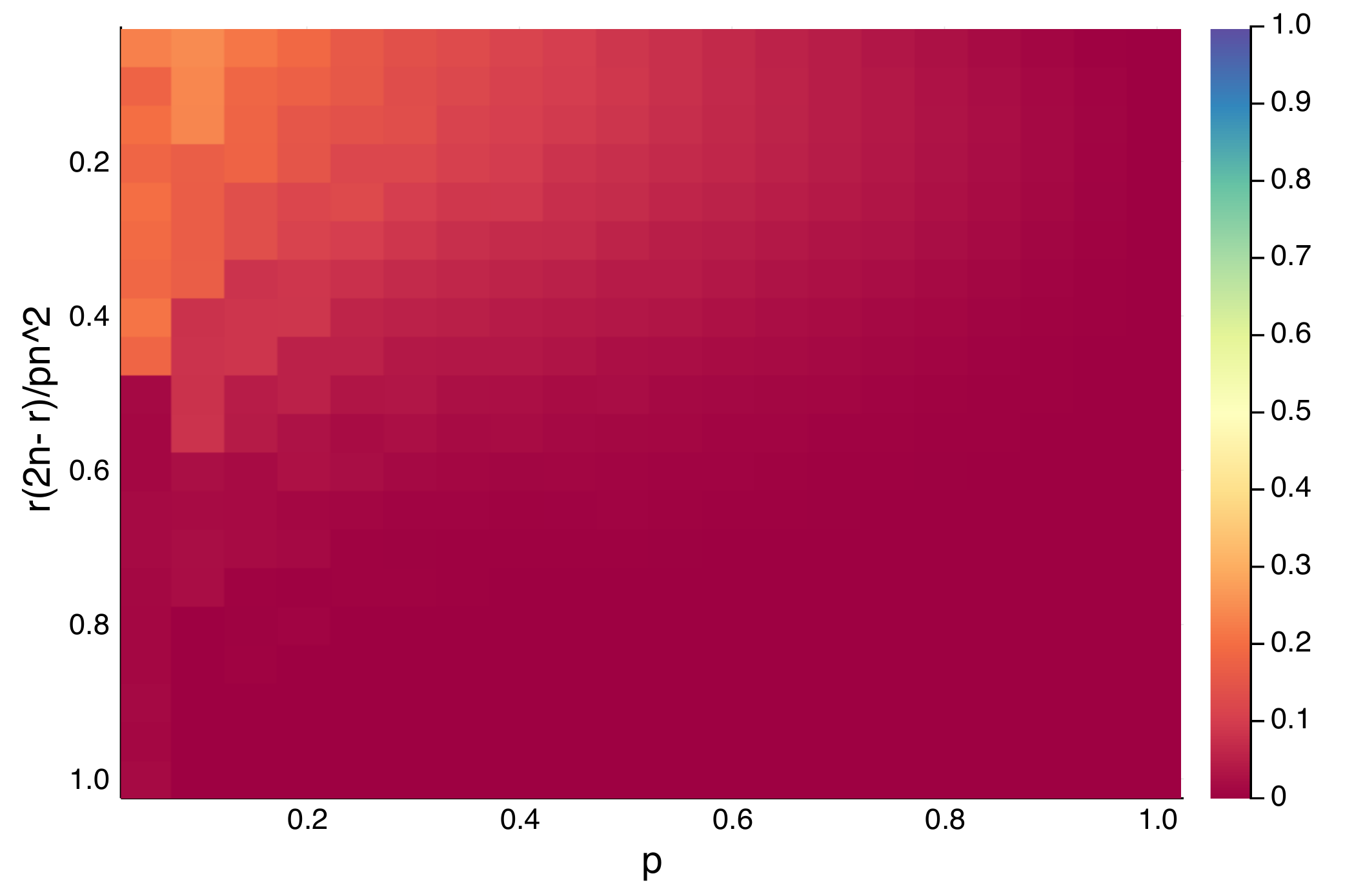}
    \caption{\footnotesize Relative in-sample bound gap ($\gamma=\frac{5}{p}$)}
\end{subfigure} \\
\begin{subfigure}[t]{.42\linewidth}
    	\centering
    	\includegraphics[width=.85\linewidth]{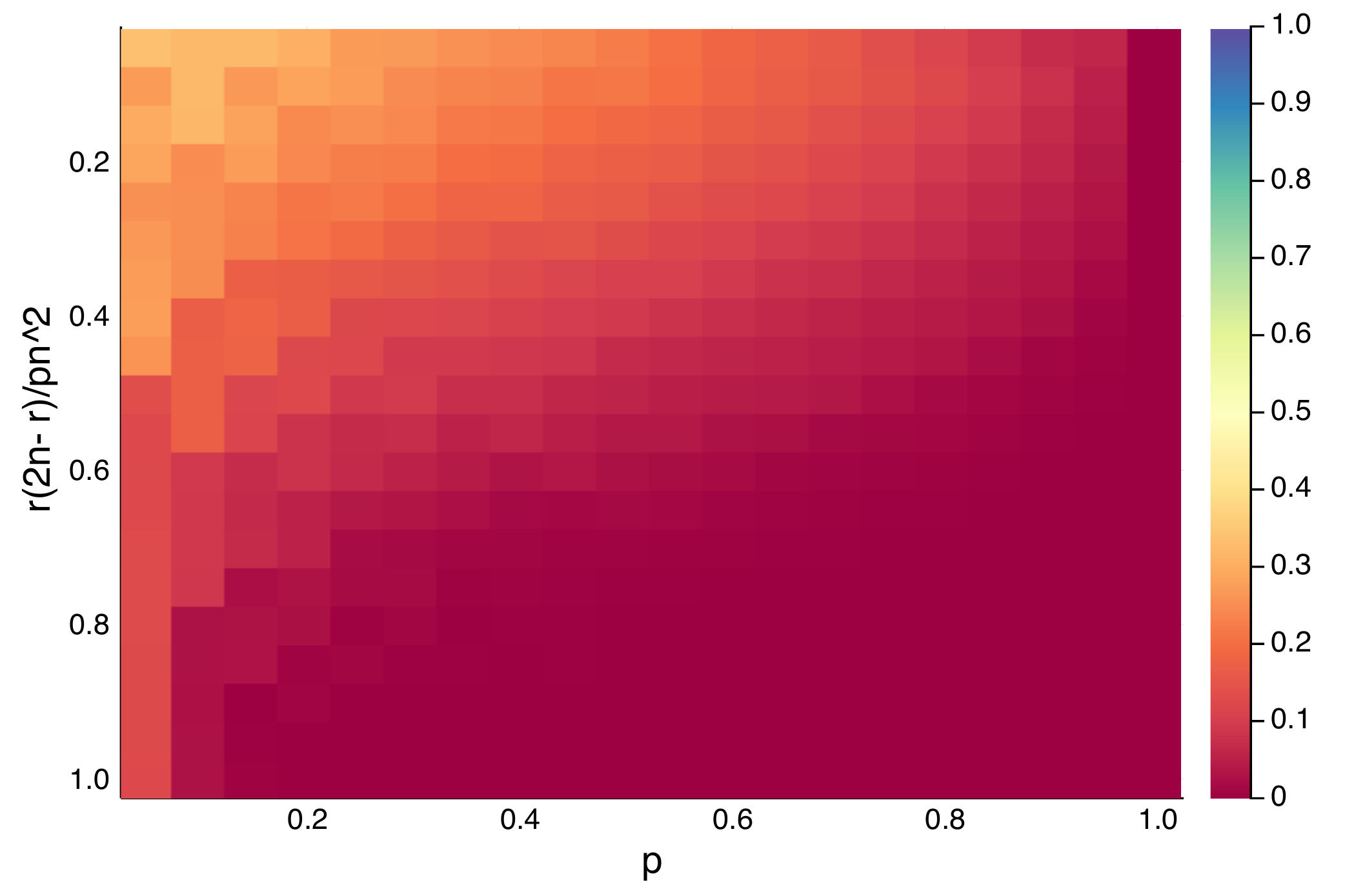}
        \caption{\footnotesize Relative in-sample bound gap ($\gamma=\frac{50}{p}$)}
    \end{subfigure} %
    \begin{subfigure}[t]{.42\linewidth}
    	\centering
    	\includegraphics[width=.85\linewidth]{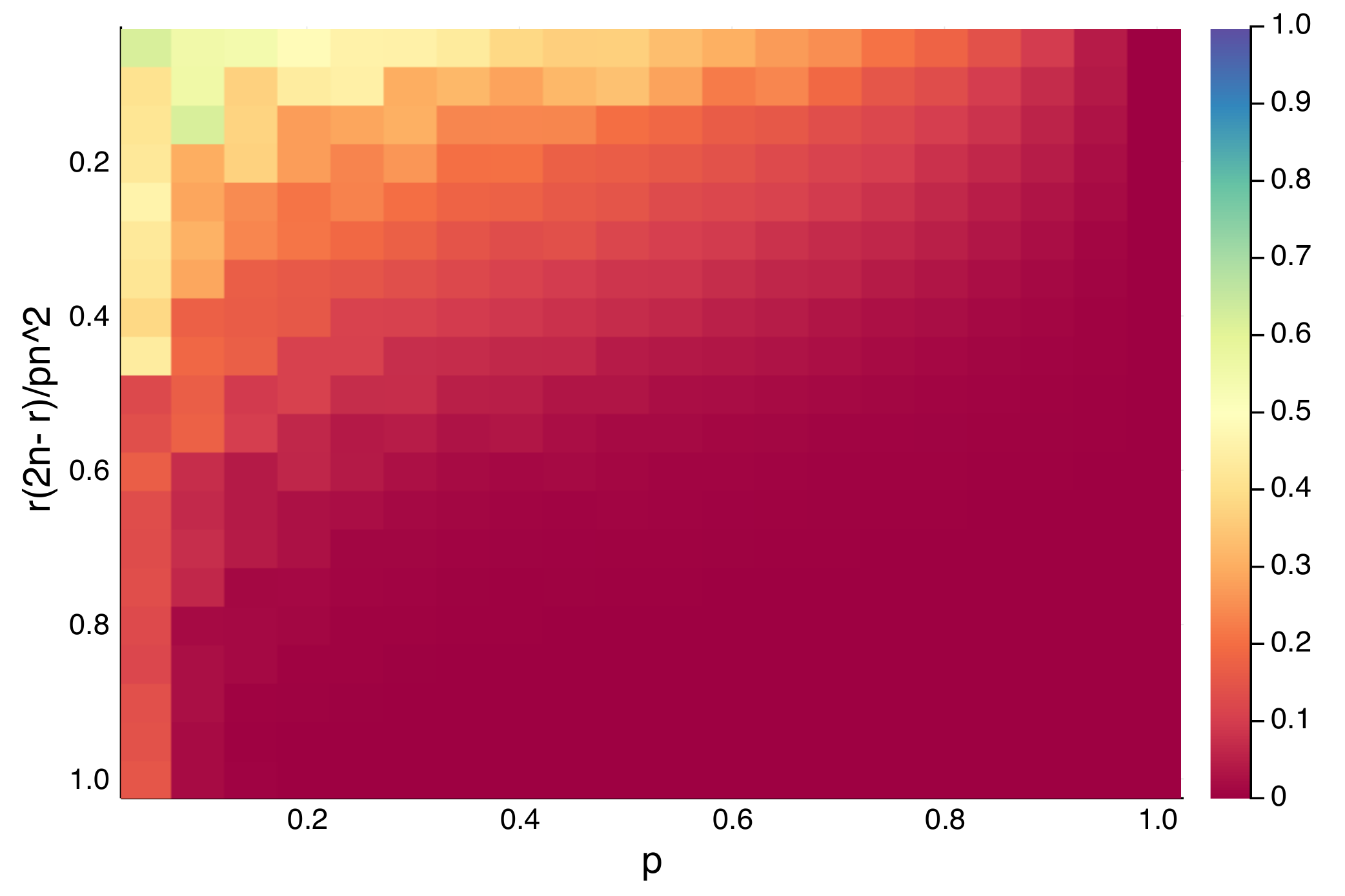}
    \caption{\footnotesize Relative in-sample bound gap ($\gamma=\frac{500}{p}$)}
    \end{subfigure}
\caption{Average relative in-sample bound gap $(\%)$, averaged over $25$ rank-$r$ matrices.}
\label{fig:bgapcomp}
\end{figure}
Observe that as $\gamma$ increases, both the duality gap and the problem's complexity increase.
This observation confirms similar results on the impact of regularization in mixed-integer conic optimization problems \citep[c.f.][]{bertsimas2018scalable,bertsimas2019unified}. Additionally, when $\gamma=\frac{500}{p}$ in Figure \ref{fig:bgapcomp}(d), the region where the in-sample duality gap is zero corresponds to exactly recovering the underlying matrix with high probability, while a strictly positive duality gap corresponds to instances with partial recovery only (see Figure \ref{fig:proprec}). This suggests a deep connection between relaxation tightness and statistical recovery.

While the relative in-sample semidefinite relaxation gap is a theoretical measure of problem difficulty, it does not indicate how fast Algorithm \ref{alg:cuttingPlaneMethod} converged in practice. In this direction, we solve the $20$ synthetic matrix completion problems considered in Table \ref{tab:comparison2} where $n\in \{10, 20\}$, $r=1$, $p \in \{0.2, 0.3\}$ for $20$ different values of $\gamma \in [10^0, 10^4]$ (distributed uniformly on a log-scale), and compare the relative in-sample semidefinite gap (greedily rounded solution vs. semidefinite bound) with Algorithm \ref{alg:cuttingPlaneMethod}'s runtimes in Figure \ref{fig:sensitivitytogamma3}, for the single-tree (left panel) and multi-tree (right panel) implementation. Results are averaged over $20$ random synthetic instances per value of $\gamma$.  We observe that the relaxation gap does correlate with runtime for single-tree. Yet, the relationship between the relaxation gap and runtime is less straightforward for multi-tree, as it depends on how Gurobi balances cut generation and node expansion, and the conditioning of the problem.
\begin{figure}[h!]
    \centering
    \begin{subfigure}[t]{.45\linewidth}
    \centering
            \includegraphics[scale=0.4]{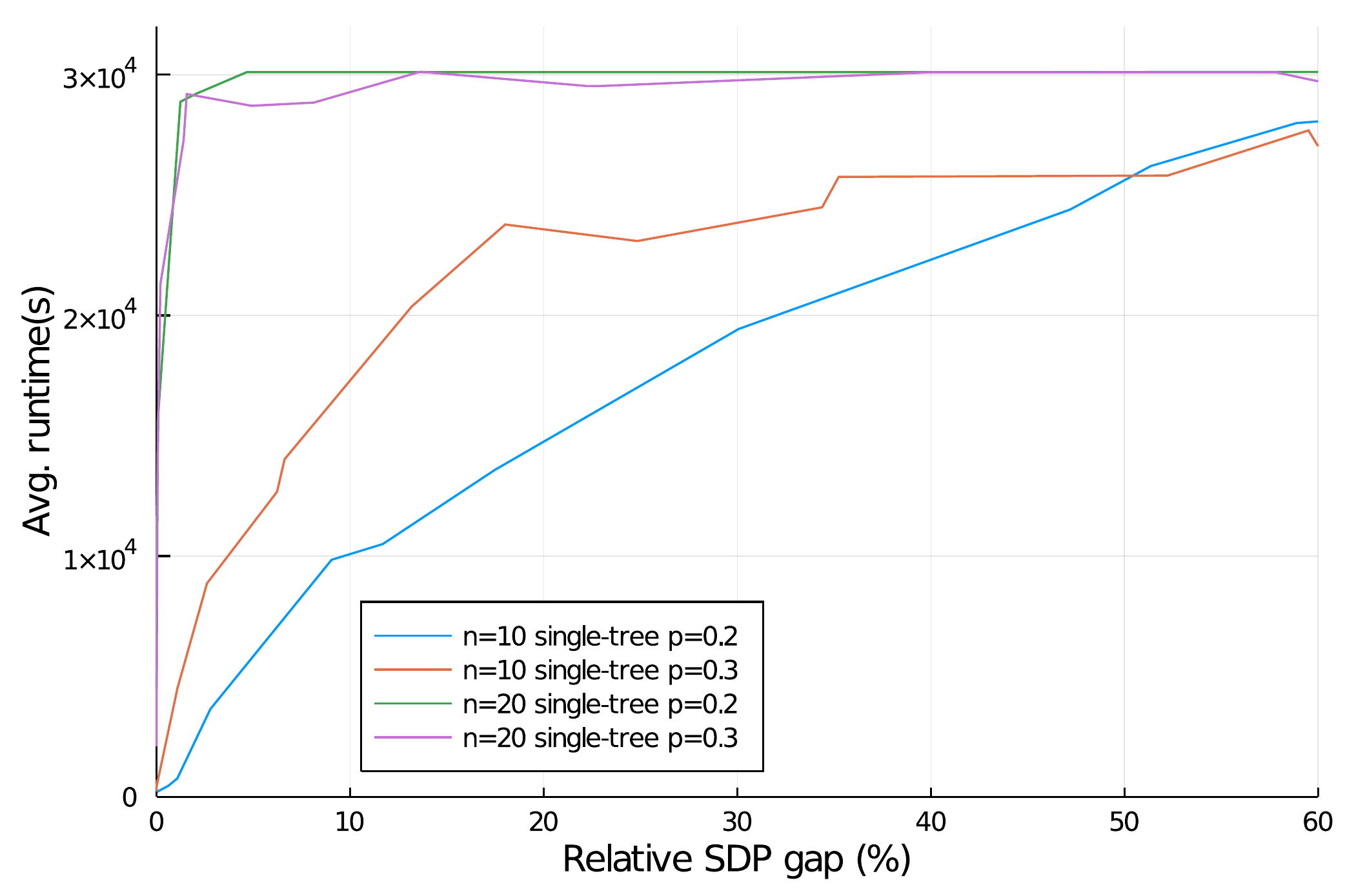}
    \end{subfigure}
        \begin{subfigure}[t]{.45\linewidth}
            \includegraphics[scale=0.4]{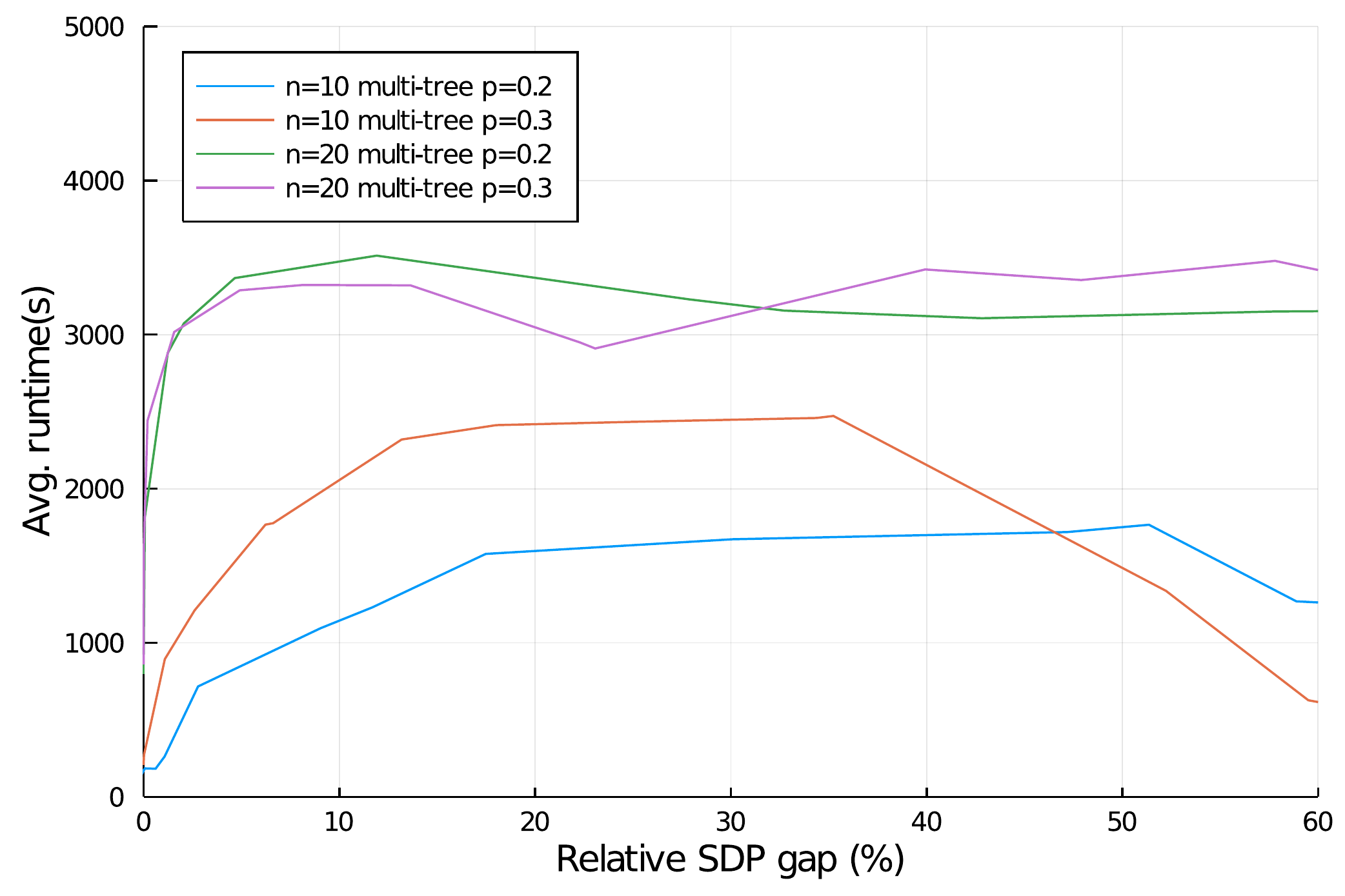}
    \end{subfigure}
   \caption{\color{black}Average runtime against relative semidefinite relaxation gap for Algorithm \ref{alg:cuttingPlaneMethod} single-tree (left) and multi-tree (right) over $20$ synthetic matrix completion instances per data point, where $p \in \{0.2, 0.3\}$, $r=1$, $n \in \{10, 20\}$.}
   \label{fig:sensitivitytogamma3}
\end{figure}

The regularizer $\gamma$ also impact the bias term $\tfrac{1}{2\gamma} \| \bm{X} \|_F^2$ added to the objective function, hence the suboptimality of the solution. To further illustrate the impact of the regularizer $\gamma$ on solve times and the trade-off between tractability and sub-optimality, Figure \ref{fig:sensitivitytogamma} reports the average runtime and MSE for the previously solved instances, as a function of $\gamma$. Figure \ref{fig:sensitivitytogamma} illustrates how
$\gamma$ balances tractability (runtime, top row) and optimality of the solution (MSE, bottom row). Also, single-tree (left panel) is one order of magnitude slower than multi-tree (right panel), and is also more numerically instable when $\gamma$ increases, largely because of the difficulty of combining a non-convex master problem and lazy constraint callbacks (which imposes many cuts, without processing the implications of these cuts as quickly). Echoing our findings in the previous section, this suggests that, while in MICO single-tree typically outperforms multi-tree, at the current state of technology multi-tree should be considered as a viable and potentially more efficient alternative for matrix completion problems which have non-convex master problems. However, as the algorithmic implementations of non-convex QCQOP solvers mature, this finding should be revisited.
\begin{figure}[h!]
    \centering
    \begin{subfigure}[t]{.45\linewidth}
    \centering
            \includegraphics[scale=0.4]{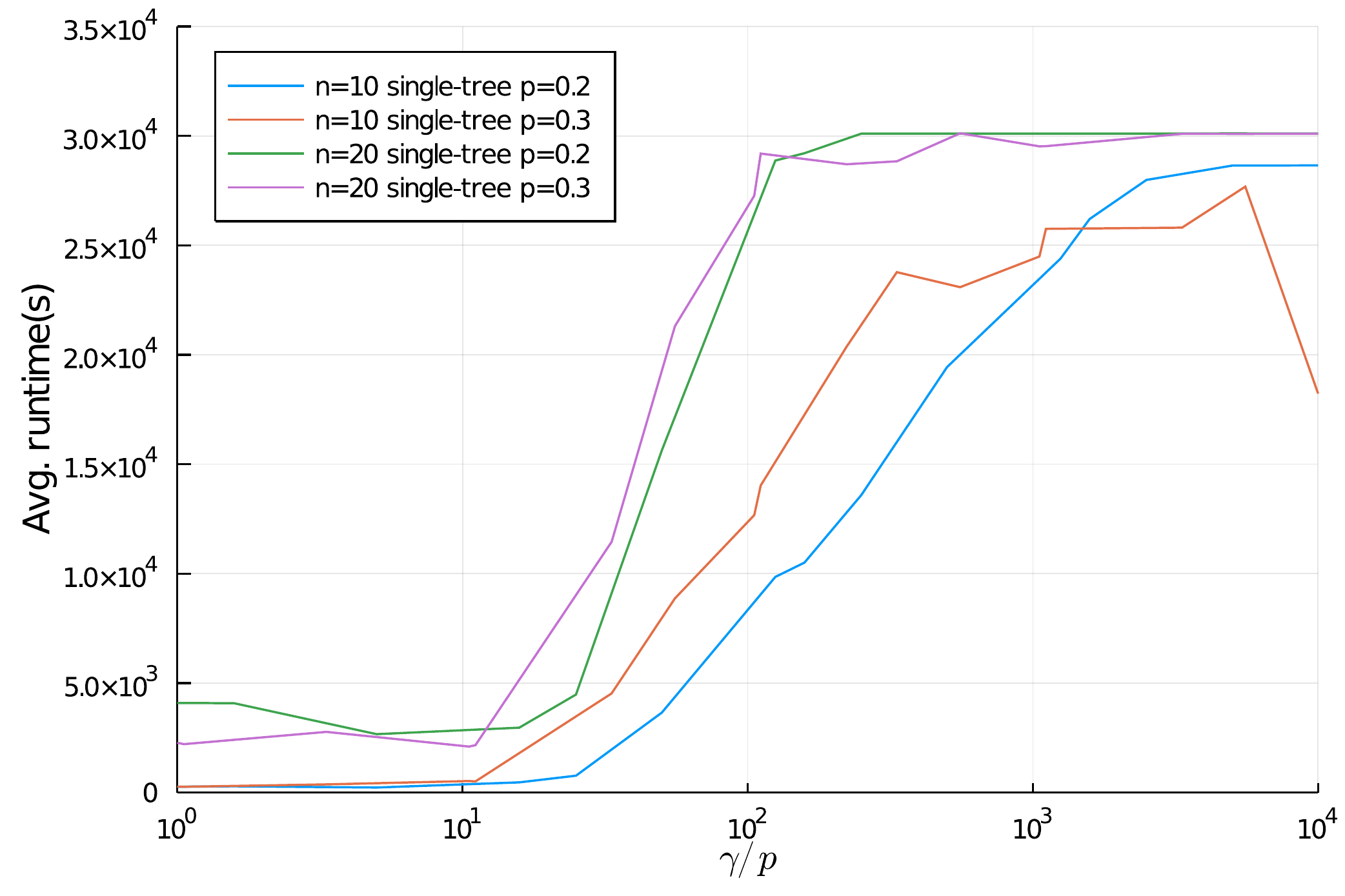}
    \end{subfigure}
        \begin{subfigure}[t]{.45\linewidth}
            \includegraphics[scale=0.4]{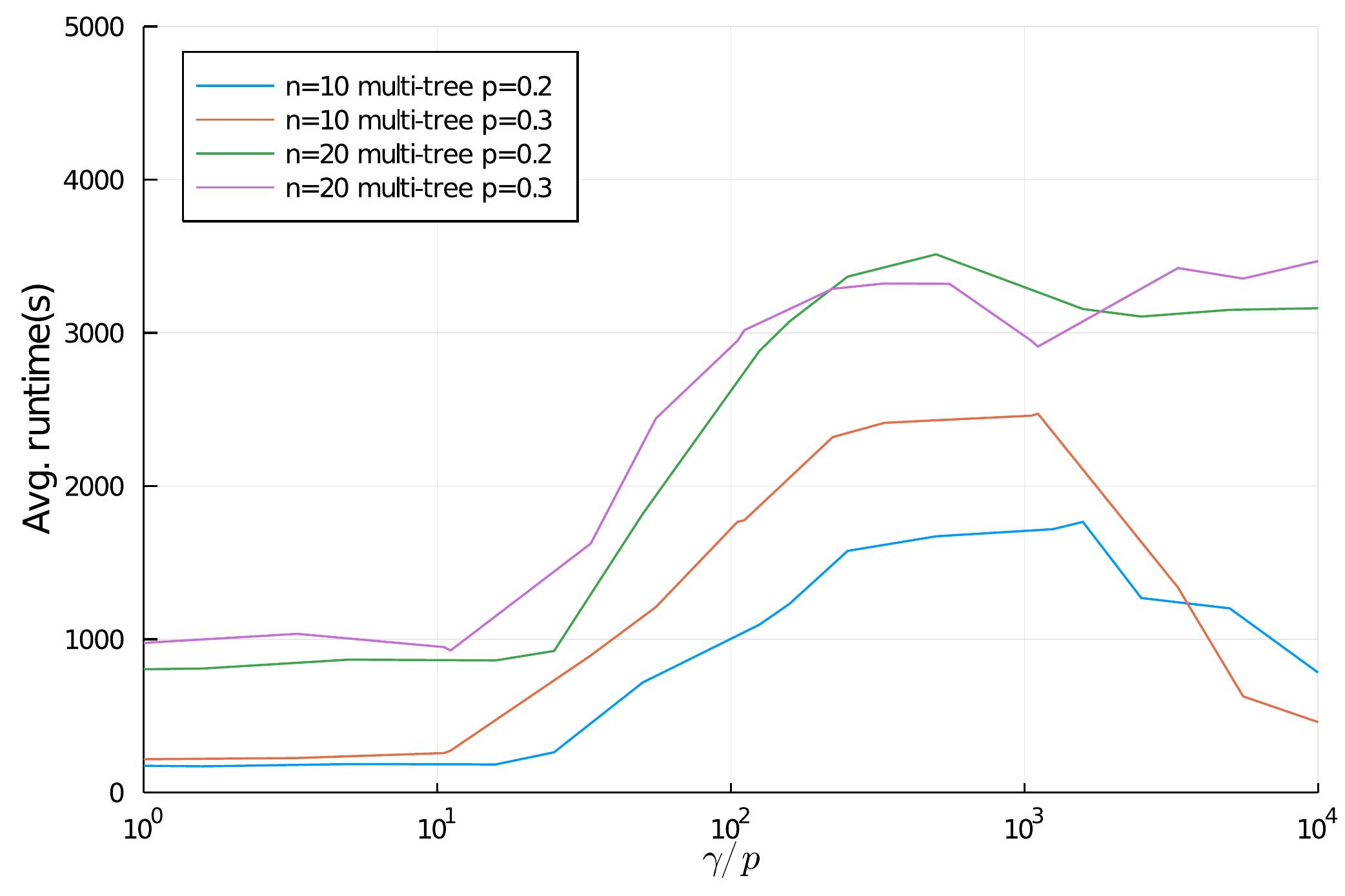}
    \end{subfigure}\\
    \begin{subfigure}[t]{.45\linewidth}
    \centering
            \includegraphics[scale=0.4]{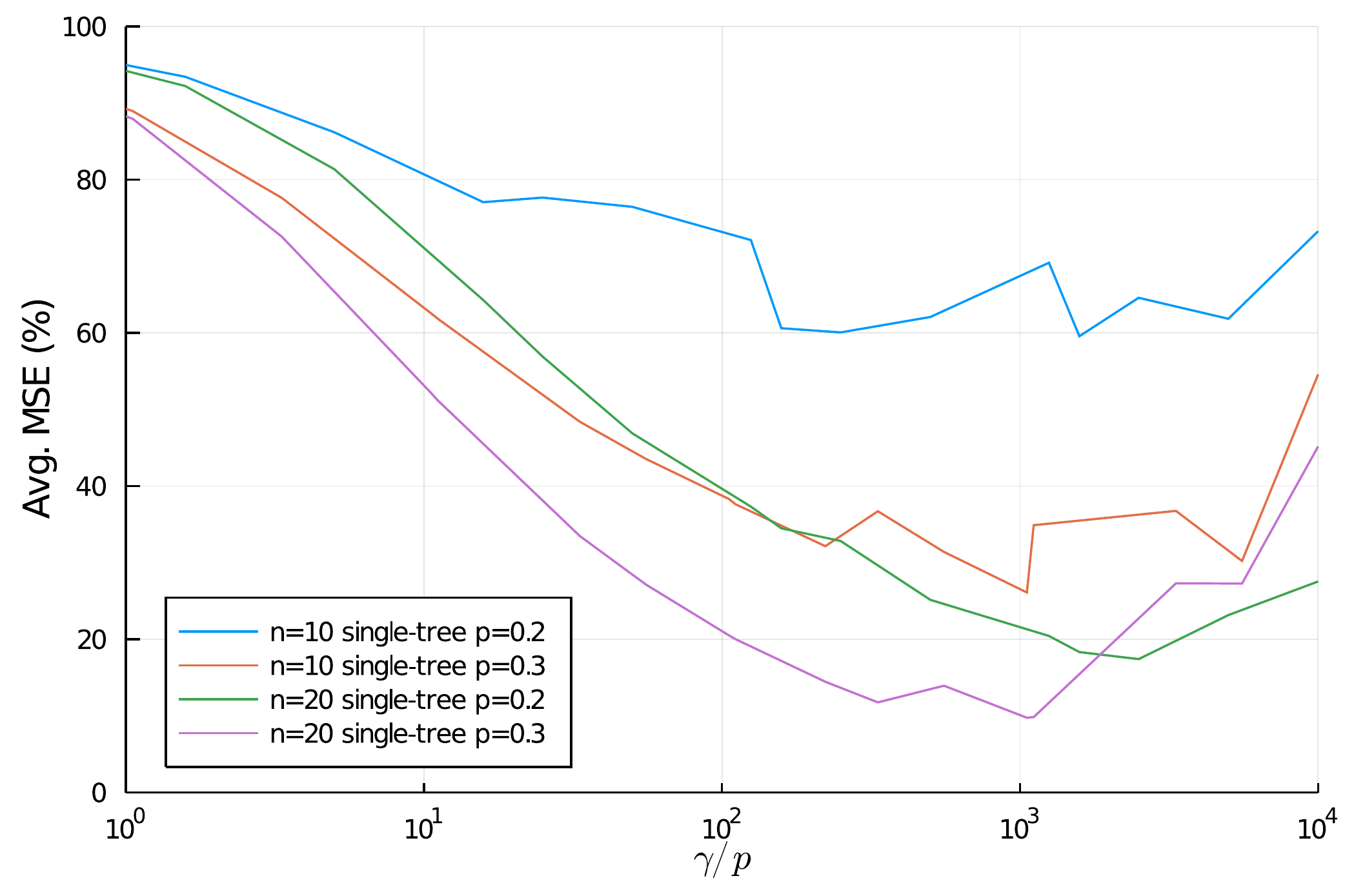}
    \end{subfigure}
        \begin{subfigure}[t]{.45\linewidth}
            \includegraphics[scale=0.4]{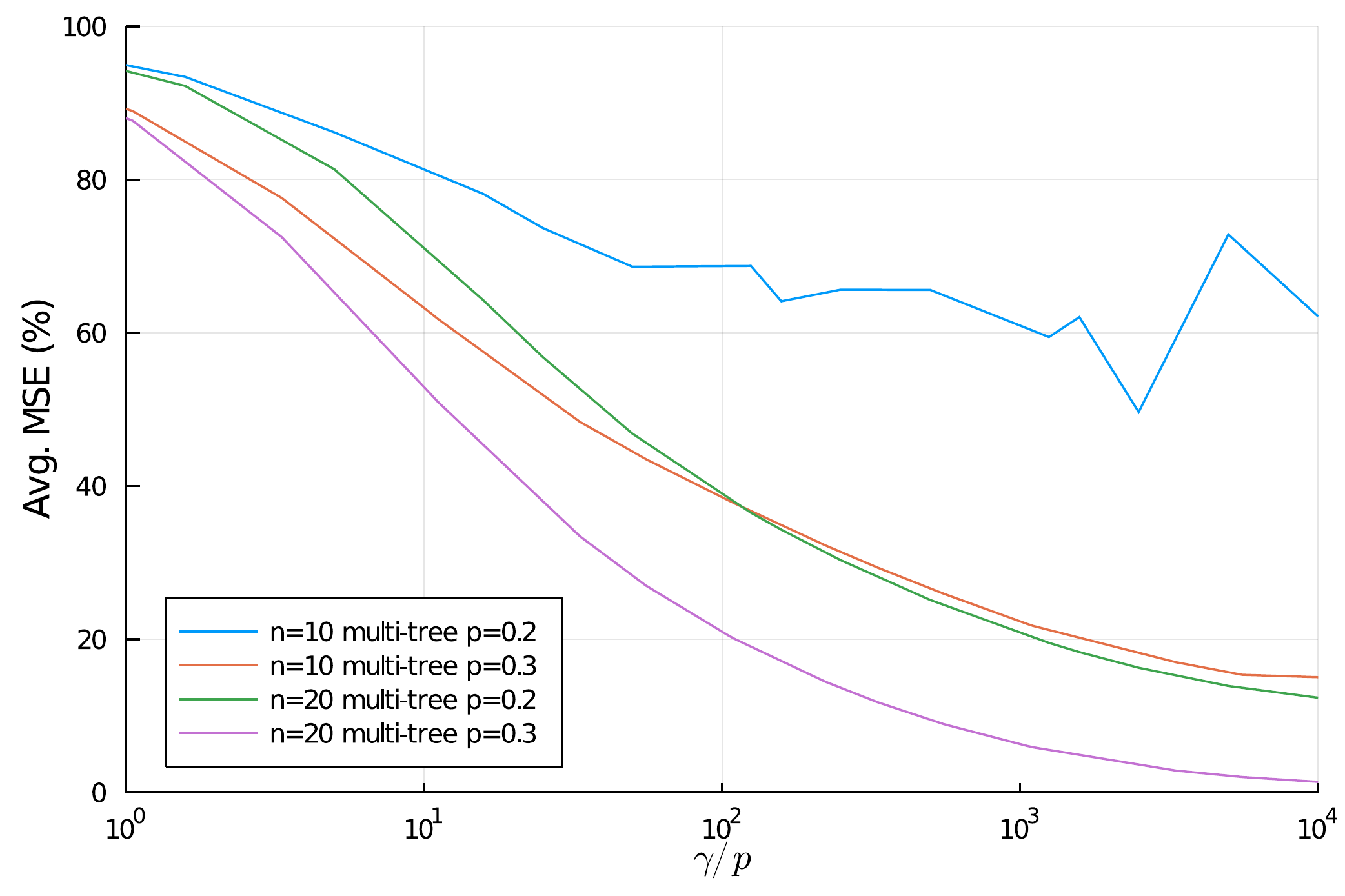}
    \end{subfigure}
   \caption{\color{black}Average runtime (top) and MSE (bottom) vs. $\gamma$ for Algorithm \ref{alg:cuttingPlaneMethod} single-tree (left) and multi-tree (right) implementations over $20$ synthetic matrix completion instances where $p \in \{0.2, 0.3\}$, $r=1$ and $n \in \{10, 20\}$. The same random seeds were used to generate random matrices completed by single-tree and multi-tree.}
   \label{fig:sensitivitytogamma}
\end{figure}
}

{
\subsection{Benchmarking Algorithm \ref{alg:cuttingPlaneMethod} on Synthetic Coordinate Recovery Problems}\label{ssec:syntheticcoordrec}
We now benchmark the performance of Algorithm \ref{alg:cuttingPlaneMethod} on anchor-free synthetic coordinate recovery problems, as previously studied by \cite{biswas2004semidefinite, luo2010semidefinite} among others.

Specifically, we sample $n$ coordinates $\bm{x}_i$ uniformly over $[-0.5, 0.5]^k$ for $k \in \{2, 3\}$, and attempt to recover a noisy Gram matrix $\bm{G} \in S^n_+$ of the $\bm{x}_i$'s, given a subset of observations of the underlying matrix. Similarly to \citet{biswas2004semidefinite}, we supply the distance between the points $D_{i,j}=\Vert \bm{x}_i-\bm{x}_j\Vert_2^2+z$, where $z\sim \mathcal{N}(0, 0.01)$, if and only if the radio range between the two points is such that $D_{i,j} \leq d_{radio}^2$. Note that we solve these problems in precisely the same fashion as the largest matrix completion problems solved in the previous section (multi-tree, with a limit of $50$ cut passes etc.)

Formally, in order to account for noise in the observed entries, we solve the following problem:
\begin{align*}
    \min_{\bm{Y} \in \mathcal{Y}^k_n}\ \min_{\bm{G} \in S^n_+} \quad & \frac{1}{2\gamma}\Vert \bm{G}\Vert_F^2+ \mathrm{tr}(\bm{G})+\lambda \cdot\Vert \bm{\xi}\Vert_1 \quad
    \text{s.t.} \quad G_{i,i}+G_{j,j}-2G_{i,j}+\xi_{i,j}=D_{i,j}\ \forall (i,j) \in \Omega_,\ \bm{G}=\bm{Y}\bm{G},
\end{align*}
    where $\lambda>0$ is a penalty term which encourages robustness, and the Frobenius norm objective likewise encourages robustness against noise in $\bm{G}$. The performance of Algorithm \ref{alg:cuttingPlaneMethod} (multi-tree) on various synthetic instances is reported in Table \ref{tab:comparison5}, for $\gamma, n, d_{radio}, k$ varying.

\begin{table}[h]
\centering
\footnotesize
\caption{Scalability of Algorithm \ref{alg:cuttingPlaneMethod} (multi-tree) for solving sensor location problems to certifiable optimality, averaged over $20$ random instances per row. A ``-'' denotes an instance that cannot be solved within the time budget, because Gurobi fails to accept our warm-start and cannot find a feasible solution. We let $\lambda=n^2$ for all instances.}
\begin{tabular}{@{}l l l r r r r r r r r r r  r@{}} \toprule
  &  & &\multicolumn{4}{c@{\hspace{0mm}}}{Rank-$2$} &\multicolumn{4}{c@{\hspace{0mm}}}{Rank-$3$}  \\
\cmidrule(l){4-7} \cmidrule(l){8-11}
$n$ & $d_{radio}$ & $\gamma$ & Time(s) & Nodes & Gap & Cuts & Time(s) & Nodes & Gap & Cuts \\\midrule
$10$ & $0.1$ & $1/p$ & $135.3$ & $6,926$ & $0.0001$ & $1.00$ & $45.14$ & $0.02$ & $0.0000$ & $1.00$\\
$10$ & $0.2$ & $1/p$ & $3,189$ & $5,249$ & $0.0024$ & $11.5$ & $216.8$ & $7,819$ & $0.0022$ & $1.00$\\
\midrule
$10$ & $0.1$ & $100/p$ & $76.2$ & $1,155$ & $0.0000$& $1.00$ & $140.6$ & $950$ & $0.0000$ & $1.00$ \\
$10$ & $0.2$ & $100/p$ & $480.6$ & $0.05$ & $0.0001$ & $21.7$ & $92.6$ & $139$ & $0.0000$ & $1.14$\\
\midrule
$20$ & $0.1$ & $1/p$ & $3,475$ & $4,548$ & $0.0007$ & $13.0$ & $3,090$ & $9,740$ & $0.0001$ & $1.00$\\
$20$ & $0.2$ & $1/p$ & $73,000$ & $0.50$ & $0.0149$ & $50.0$ & $7,173$ & $5,313$ & $0.0038$ & $1.20$\\
\midrule
$20$ & $0.1$ & $100/p$ & $1,878$ & $0.00$ & $0.0000$ & $3.91$ & $64.9$ & $0.00$ & $0.0000$ & $1.07$\\
$20$ & $0.2$ & $100/p$ & $67,530$ & $0.20$ & $0.0044$ & $50.0$ & $55.7$ & $0.00$ & $0.0002$ & $1.00$\\
\bottomrule
\end{tabular}
\label{tab:comparison5}
\end{table}
We observe that the problem’s complexity increases with the rank and with the dimensionality of the Gram matrix, although not too excessively. Indeed, Algorithm \ref{alg:cuttingPlaneMethod} can solve coordinate recovery problems with tens of data points to certifiable optimality in hours.
}

\subsection{Summary of Findings from Numerical Experiments}
Our main findings from the numerical experiments in this section are as follows:
\begin{itemize}
    \item As demonstrated in Section \ref{ssec:altminexper}, Algorithm \ref{alg:altmin} successfully solves convex relaxations of low-rank problems where $n=100$s, in a faster and more scalable fashion than state-of-the-art interior point codes such as \verb|Mosek|.
    \item As demonstrated in Section \ref{ssec:compwithheurmethods}, the following strategy is almost as effective as solving a low-rank problem exactly: solving the convex relaxation using Algorithm \ref{alg:altmin}, greedily rounding the solution to the convex relaxation, and using this greedily rounded solution as a warm-start for a local method such as the method of \citet{burer2003nonlinear, burer2005local}. The success of this strategy can be explained because solving a relaxation and rounding is a global strategy which matches
    the low-order moments of the set of optimal solutions to obtain a solution near an optimal solution, while local methods polish a solution by seeking the best solution within some neighborhood of an initial point.
    \item As demonstrated in {\color{black}Section} \ref{ssec:impactofreg}, increasing the amount of regularization in a low-rank problem by decreasing $\gamma$ decreases the duality gap between a low-rank problem with Frobenius or spectral norm problem, and its convex relaxation. Therefore, increasing the amount of regularization makes the problem easier in a practical sense (although not necessarily in a complexity-theoretic sense).
    \item As demonstrated in Sections \ref{ssec:benchmarkalg1}, \ref{ssec:syntheticcoordrec}, Algorithm \ref{alg:cuttingPlaneMethod} scales to solve problems where $n$ is in the tens, i.e., hundreds or thousands of decision variables, in hours. Moreover, the main bottleneck inhibiting solving problems where $n$ is in the hundreds or thousands is that we solve our master problems using Gurobi, a {\color{black}QCQO} solver which translates the orthogonal projection matrix constraint into many piecewise linear constraints. This suggests that a custom branch-and-bound solver which explicitly models orthogonal projection matrices constitutes a promising area for future work.
\end{itemize}

\section{Conclusion}
In this paper, we introduced Mixed-Projection Conic Optimization, a new framework for modeling rank constrained optimization problems that, for the first time, solves low-rank problems to certifiable optimality at moderate problem sizes. We also provided a characterization of the complexity of rank constraints, and proposed new convex relaxations and rounding methods that lead to viable and more accurate solutions than those obtained via existing techniques such as the log-det or nuclear norm heuristics.
Inspired by the collective successes achieved in mixed-integer optimization, we hope that MPCO constitutes an exciting new research direction for the optimization community. {For instance, we believe that
custom branch-and-bound solvers that
explicitly model orthogonal projection matrices could further enhance the scalability of the MPCO framework.}

{\color{black}
\subsubsection*{Acknowledgments:} We thank the two anonymous referees and the associate editor for many valuable comments which improved the paper substantially.
}

{
\bibliographystyle{abbrvnat}
\setlength{\bibsep}{0pt plus 0.0ex}

}
\begin{appendices}
{
\section{Omitted Proofs}\label{sec:omittedproofs}
In this section, we supply omitted proofs of the results stated in the manuscript, in the order in which the results were stated.
\subsection{Proof of Theorem \ref{existentialtheoryofrealshard}}\label{ssec:existentialtheoryofrealshard}
\proof{Proof of Theorem \ref{existentialtheoryofrealshard}}
By \citep[Theorem 2.49]{blekherman2012semidefinite}, a set of distances $d_{i,j}$ can be embedded in a Euclidean space of dimension $k$ if and only if there exists some Gram matrix $\bm{G}$ such that $\bm{D}=\mathrm{Diag}(\bm{G})\bm{e}^\top +\bm{e}\mathrm{Diag}(\bm{G})^\top -2\bm{G}$ where $\bm{G} \succeq \bm{0}$ and $\mathrm{Rank}(\bm{G}) \leq k$. Therefore, Proposition \ref{prop:graphreduction}'s $\exists \mathbb{R}$-hard problem is a special case of:
\begin{align*}
    \min_{\bm{G} \in S_+^n} \ \mathrm{Rank}(\bm{G})\quad
    \text{s.t.} \quad \mathrm{Diag}(\bm{G})\bm{e}^\top +\bm{e}\mathrm{Diag}(\bm{G})^\top -2\bm{G}= \bm{D},
\end{align*}
where $k=2$, $D_{i,j}=1 \ \forall (i,j) \in E$, and we do not impose the $(i,j)th$ equality otherwise.
\hfill\Halmos \endproof
\subsection{Proof of Theorem \ref{existentialtheorycompleteness}}\label{ssec:proofofthmexistentialcomp}
\proof{Proof of Theorem \ref{existentialtheorycompleteness}}
To establish this result, it suffices to perform a reduction from the following feasibility system to a polynomially sized system of polynomial equalities and inequalities:
\begin{align*}
    \exists \bm{X}: \quad & \mathrm{Rank}(\bm{X}) \leq k,\quad \langle \bm{A}_i, \bm{X} \rangle = b_i \quad \forall i \in [m],  \quad \bm{X} \succeq \bm{0},
\end{align*}
because testing the feasibility of the above system a polynomial number of times for $k \in [n]$ certainly solves Problem \eqref{rankminproblem_proj}.

By Proposition \ref{prop:rankconstraintreformulation}, this system is feasible if and only if the following system is also feasible:
\begin{align*}
    \exists \bm{X}, \bm{Y}: \quad & \bm{Y}^2=\bm{Y},\quad \bm{X}=\bm{X}\bm{Y},\quad  \mathrm{tr}(\bm{Y}) \leq k,\  \langle \bm{A}_i, \bm{X} \rangle = b_i \quad \forall i \in [m],  \quad \bm{X} \succeq \bm{0},
\end{align*}

The result then follows by observing that semidefinite constraints are indeed semialgebraic constraints, by the Tarski-Seidenberg theorem \citep[see][Chapter A.4.4]{blekherman2012semidefinite}, and therefore this system is equivalent to a polynomially sized system of polynomial equalities and inequalities.
\hfill\Halmos \endproof

\subsection{Complexity of Low-Rank Integer Optimization and Hilbert's $\nth{10}$ Problem}\label{ssec:undecidableproof}
We now demonstrate that imposing the constraint $\bm{X} \in \mathbb{Z}^{n \times n}$ changes Problem \eqref{rankminproblem}'s complexity status, by making it undecidable. Formally, we have:
\begin{theorem}\label{rankminintegerundecidable}
The following problem is undecidable, even when $k=1$ and its objective is binary:
\begin{align*}
    \min_{\bm{X} \in \mathbb{Z}^{n \times n}} \quad & \langle \bm{C}, \bm{X} \rangle \quad \text{s.t.} \quad \langle \bm{A}_i, \bm{X} \rangle = b_i \quad \forall i \in [m],\ \mathrm{Rank}(\bm{X}) \leq k,\ \bm{X} \succeq \bm{0}.
\end{align*}
\end{theorem}
\begin{remark}
Theorem \ref{rankminintegerundecidable}'s reduction does not hold in the presence of a regularizer. Indeed, imposing either a spectral or Frobenius norm regularizer ensures the boundedness of the problem's level sets, which allows the problem to be solved in finite time via branch-and-bound. Nonetheless, Theorem \ref{rankminintegerundecidable} shows that imposing an integrality constraint makes Problem \eqref{rankminproblem} \textit{much} harder.
\end{remark}

\proof{Proof of Theorem \ref{rankminintegerundecidable}}
We perform a reduction from integer optimization with quadratic constraints, which is undecidable when the objective is binary \cite[][]{jeroslow1973there}, by reduction from Hilbert's $\nth{10}$ problem. Recall that integer optimization with quadratic constraints is definitionally:
\begin{equation}
\begin{aligned}
  \min_{\bm{x} \in \mathbb{Z}^n} \quad \langle \bm{c}, \bm{x}\rangle\quad \text{s.t.} \quad \langle \bm{x},\bm{Q}_i \bm{x}\rangle+\langle \bm{a}_i, \bm{x}\rangle \leq b_i \quad \forall i \in [m].
\end{aligned}
\end{equation}
We now show that this problem is equivalent to:
\begin{align*}
\min_{\bm{x} \in \mathbb{Z}^n, \bm{X} \in \mathbb{Z}^{n \times n}} \quad & \langle \bm{c}, \bm{x} \rangle \quad \text{s.t.} \quad \langle \bm{Q}_i, \bm{X}\rangle+\langle \bm{a}_i,  \bm{x}\rangle\leq b_i \quad \forall i \in [m],\quad \mathrm{Rank}     \begin{pmatrix}
    1 & \bm{x}^\top\\
    \bm{x} & \bm{X}
    \end{pmatrix} \leq 1,\  \begin{pmatrix}
    1 & \bm{x}^\top\\
    \bm{x} & \bm{X}
    \end{pmatrix} \succeq \bm{0}.
\end{align*}

To establish the result, it suffices to show that $\bm{X}=\bm{x}\bm{x}^\top$ in any feasible solution to the second problem. We now show this, by appealing to the Guttman rank identity (see Lemma \ref{guttmanrankidentity}).
We remind the reader that an equivalent form of the Guttman rank formula is the identity:
\begin{align*}
    \mathrm{Rank}& \begin{pmatrix}
    1 & \bm{x}^\top\\
    \bm{x} & \bm{X}
    \end{pmatrix}=\mathrm{Rank}(\bm{X})+\mathrm{Rank}(\bm{X}
    -\bm{x}\bm{x}^\top),
\end{align*}
Since $\bm{X} \succeq \bm{0}$ and the left hand side of the above expression is either $0$ or $1$, this formula implies that either $\bm{X}=\bm{x}\bm{x}^\top$ or $\bm{X}=\bm{0}$. However, the latter case can only hold if $\bm{x}=\bm{0}$, since $\bm{X} \succeq \bm{x}\bm{x}^\top$ by Schur complements. Therefore, $\bm{X}=\bm{x}\bm{x}^\top$ and the result holds.
\hfill\Halmos \endproof

}

\subsection{Proof of {\color{black}Lemma} \ref{thm:perspectivereformulation}}\label{subsec:perspectivecutproof}
\proof{Proof of {\color{black}Lemma} \ref{thm:perspectivereformulation}}
Let us fix $\bm{Y} \in \mathrm{Conv}\left(\mathcal{Y}_n\right)$. Then, we have that:
\begin{align*}
\max_{\bm{\alpha}} \: h(\bm{\alpha} ) - \dfrac{\gamma}{2} \left\langle \bm{\alpha}, \bm{Y}\bm{\alpha}\right\rangle
&= \max_{\bm{\alpha}, \bm{\beta}} \: h(\bm{\alpha}) - \dfrac{\gamma}{2} \left\langle \bm{\beta}, \bm{Y}\bm{\beta}\right\rangle
\mbox{  s.t.  } \bm{\beta} = \bm{\alpha}, \\
&= \max_{\bm{\alpha}, \bm{\beta}} \: \min_{\bm{X}} \: h(\bm{\alpha}) - \dfrac{\gamma}{2} \left\langle \bm{\beta}, \bm{Y}\bm{\beta}\right\rangle - \left \langle \bm{X}, \bm{\beta} - \bm{\alpha}\right\rangle, \\
&= \min_{\bm{X}} \: \underbrace{ \max_{\bm{\alpha}} \:  \left[ h(\bm{\alpha}) + \langle \bm{X}, \bm{\alpha}\rangle \right]}_{(-h)^\star(\bm{X}) = g(\bm{X})} + \max_{\bm{\beta}}\left[\frac{-\gamma}{2}\left\langle \bm{Y}\bm{\beta}, \bm{\beta}\right\rangle-\left\langle \bm{X}, \bm{\beta}\right\rangle \right].
\end{align*}
Finally, the optimality condition with respect to $\bm{\beta}$ is
$
    \bm{Y}\bm{\beta}=\frac{-1}{\gamma}\bm{X},
$
which implies the later term is
\begin{align*}
    \max_{\bm{W}} \left[\frac{1}{2\gamma}\langle \bm{X}, \bm{Y}^\dag \bm{X}\rangle -\frac{1}{2}\langle \bm{X}, (\mathbb{I}-\bm{Y}^\dag \bm{Y})\bm{W}\rangle \right]& =\max_{\bm{W}} \left[\frac{1}{2\gamma}\langle \bm{X}, \bm{Y}^\dag \bm{X}\rangle -\frac{1}{2}\langle \bm{W}, (\mathbb{I}-\bm{Y}^\dag \bm{Y})\bm{X}\rangle \right]\\ & = \begin{dcases}
\frac{1}{2\gamma}\left\langle \bm{X}, \bm{Y}^\dag \bm{X}\right\rangle &\mbox{  if  } \bm{Y} \in \mathrm{Span}(\bm{X}),\\
+\infty &\mbox{  otherwise}.
\end{dcases}
\end{align*}
We therefore conclude that the later term is equal to $\frac{1}{2\gamma}\left\langle \bm{X}, \bm{Y}^\dag \bm{X}\right\rangle$ whenever the constraint $\bm{Y}^\dag\bm{Y}\bm{X}=\bm{X}$ holds. By the generalized Schur complement lemma \ref{generalizedschur}, this expression is equivalent to introducing a new matrix $\bm{\theta}$, imposing the term $\frac{1}{2\gamma}\mathrm{tr}(\bm{\theta})$ and requiring that
$\begin{pmatrix} \bm{\theta} & \bm{X} \\ \bm{X}^\top & \bm{Y}\end{pmatrix} \succeq \bm{0}$. \qquad \hfill\Halmos
\endproof

\subsection{Proof of {\color{black}Lemma} \ref{thm:spectralcovrelax}}\label{proof:thm:spectralcovrelax}
\proof{Proof of {\color{black}Lemma} \ref{thm:spectralcovrelax}}
Let us fix $\bm{Y} \in \mathrm{Conv}\left(\mathcal{Y}_n\right)$. Then, we have that:
\begin{align*}
& \max_{\bm{\alpha} \in S^n, \bm{W}_+, \bm{W}_- \succeq \bm{0}} \: h(\bm{\alpha} ) - M \left\langle \bm{Y}, \bm{W}_{+}-\bm{W}_{-}\right\rangle \quad \text{s.t.} \quad \bm{\alpha}=\bm{W}_{+}-\bm{W}_{-}\\
&= \max_{\bm{\alpha}\in S^n, \bm{W}_+, \bm{W}_- \succeq \bm{0}} \min_{\bm{X} \in S^n} \: h(\bm{\alpha}) - M \left\langle \bm{Y}, \bm{W}_{+}-\bm{W}_{-}\right\rangle +\langle \bm{X}, \bm{\alpha}-\bm{W}_{+}+\bm{W}_{-}\rangle\\
&= \min_{\bm{X} \in S^n} \max_{\bm{\alpha}\in S^n, \bm{W}_+, \bm{W}_- \succeq \bm{0}}  \: h(\bm{\alpha}) - M \left\langle \bm{Y}, \bm{W}_{+}-\bm{W}_{-}\right\rangle +\langle \bm{X}, \bm{\alpha}-\bm{W}_{+}+\bm{W}_{-}\rangle\\
&= \min_{\bm{X} \in S^n} \max_{\bm{\alpha}\in S^n}  \: \underbrace{\left[ h(\bm{\alpha}) + \langle \bm{X}, \bm{\alpha}\rangle \right]}_{(-h)^\star(\bm{X}) = g(\bm{X})}+\max_{\bm{W}_+, \bm{W}_- \succeq \bm{0}}\left[- M \left\langle \bm{Y}, \bm{W}_{+}-\bm{W}_{-}\right\rangle +\langle \bm{X}, -\bm{W}_{+}+\bm{W}_{-}\rangle \right].
\end{align*}
Finally, the optimality conditions with respect to $\bm{W}_{+}, \bm{W}_{-}$ imply that $-M \bm{Y} \preceq \bm{X} \preceq M \bm{Y}.$ \hfill\Halmos
\endproof

\subsection{Proof of {\color{black}Lemmas} \ref{prop:rectcase1} and \ref{prop:rectcase2}}\label{ssec:proofofnucnormprops}

\proof{Proof of {\color{black}Lemma} \ref{prop:rectcase1}}
In Problem \eqref{prob:specprob1}, it is not too hard to see that for any $\bm{X}$ an optimal choice of $\bm{Y}$ is $\bm{Y}=\frac{1}{M}\bm{X}_+ + \frac{1}{M}\bm{X}_-$, where $\bm{X}_+, \bm{X}_-$ are orthogonal positive semidefinite matrices such that $\bm{X}=\bm{X}_+ - \bm{X}_-$. Therefore, the result follows by observing that $\mathrm{tr}(\bm{X}_+ + \bm{X}_-)=\Vert \bm{X}\Vert_*$.
\hfill\Halmos \endproof

\proof{Proof of {\color{black}Lemma} \ref{prop:rectcase2}}
In Problem \eqref{prob:spectprob3}, for any feasible $\bm{X}$ we have $\Vert \bm{X}\Vert_\sigma \leq M$. Under this constraint, it follows that for any $\bm{X}$ an optimal choice of $\bm{Y}, \bm{Y}'$ is $\bm{Y}=\bm{U}\bm{\Sigma}\bm{U}^\top$, ${\bm{Y}'=\bm{V}\bm{\Sigma}\bm{V}^\top}$, where $\bm{X}=\bm{U}\bm{\Sigma}\bm{V}^\top$ is an SVD of $\bm{X}$. The result follows as $\mathrm{tr}(\bm{Y})=\mathrm{tr}(\bm{\Sigma})=\Vert \bm{X}\Vert_*$.
\hfill\Halmos \endproof

\subsection{Proof of {\color{black}Lemma} \ref{prop:generalizedrevhuber}}\label{ssec:appendgenrevhuber}

\proof{Proof of {\color{black}Lemma} \ref{prop:generalizedrevhuber}}
Observe that, by the Generalized Schur Complement Lemma (see, e.g., lemma \ref{generalizedschur}), an optimal choice of $\bm{\theta}$ in Problem \eqref{prob:revhub1} is $\bm{\theta}=\bm{X}\bm{Y}^\dag \bm{X}^\top$. Therefore, we can eliminate $\bm{\theta}$ from Problem \eqref{prob:revhub1}, to obtain the equivalent objective:
\begin{align*}
     \min_{\bm{Y} \in \mathrm{Conv}(\mathcal{Y}_n)}  \ \min_{\bm{X} \in \mathbb{R}^{n \times m}} \quad & \lambda \cdot \mathrm{tr}(\bm{Y})+g(\bm{X})+\frac{1}{2\gamma}\langle \bm{X}\bm{X}^\top, \bm{Y}^\dag\rangle.
\end{align*}
Moreover, by the rank-nullity theorem \citep[see, e.g.,][Chapter 0.2.3]{johnson1985matrix}, we can split the columns of $\bm{Y}$ into columns in the span of the columns of $\bm{X}$ and columns orthogonal to the columns of $\bm{X}$. Since the columns orthogonal to the columns of $\bm{X}$ do not affect the objective value, it follows that we can write $\bm{Y}^\dag=\sum_{i=1}^n \frac{1}{\theta_i}\bm{u}_i \bm{u}_i^\top$ without loss of optimality, where $\bm{X}\bm{X}^\top=\bm{U}\bm{\Sigma}\bm{U}^\top$ is an SVD of $\bm{X}\bm{X}^\top$, and $0 \leq \theta_i \leq 1$ for each $\theta_i$, because $\bm{Y} \in \mathrm{Conv}(\mathcal{Y}_n)$. Problem \eqref{prob:revhub1} then becomes:
\begin{align*}
\min_{\bm{X} \in \mathbb{R}^{n \times m}, \bm{\theta} \in \mathbb{R}^n: \ \bm{0} \leq \bm{\theta} \leq \bm{e}} \quad & g(\bm{X})+\sum_{i=1}^n \left(\lambda \theta_i +\frac{\sigma_i(\bm{X})^2}{2\gamma \theta_i}\right).
\end{align*}
The result then follows because, for any $\lambda > 0$, \citep[c.f.][Equation (30)]{pilanci2015sparse}
\begin{align*}
    \min_{0 \leq \theta \leq 1} \left[ \lambda \theta +\frac{t^2}{\theta}\right]=\begin{cases} 2\sqrt{\lambda}\vert t \vert, & \text{if} \ \vert t \vert \leq \sqrt{\lambda},\\
    t^2+\lambda, & \text{otherwise.}\end{cases} \hfill\Halmos
\end{align*}
\endproof

\subsection{Proof of Theorem \ref{proofofconvergence}}
\label{ssec:proofofconv}
\proof{Proof of Theorem \ref{proofofconvergence}}
We only detail the proof of $\epsilon$-optimality; the proof of $\epsilon$-feasibility is almost identical \cite[see][]{mutapcic2009cutting}.

Suppose that at some iteration $k>1$, Algorithm \ref{alg:cuttingPlaneMethod} has not converged. Then, $$\theta_k-f(\bm{Y}_k) < -\epsilon, \qquad
\text{and} \quad \theta_k \geq f(\bm{Y}_i) + \langle \bm{H}_{i}, \bm{Y}_k-\bm{Y}_i\rangle\ \forall i < k.$$ But $\theta_k \leq f(\bm{Y}_i),$ since $\theta_k$ and $f(\bm{Y}_i)$ are respectively valid lower and upper bounds on the optimal objective. Therefore, $\langle \bm{H}_{i}, \bm{Y}_k-\bm{Y}_i\rangle \geq 0$. Putting the two inequalities together then implies that
$$f(\bm{Y}_k)-f(\bm{Y}_i) > \epsilon+\langle H_{i}, \bm{Y}_k-\bm{Y}_i\rangle \geq \epsilon, \ \text{or equivalently} \quad \epsilon < f(\bm{Y}_k)-f(\bm{Y}_i) \leq L \Vert \bm{Y}_i-\bm{Y}_k\Vert_F,$$
where the second inequality holds by Lipschitz continuity. Rearranging this inequality implies that ${\Vert \bm{Y}_i-\bm{Y}_k\Vert_F > \frac{\epsilon}{L}}$, i.e., Algorithm \ref{alg:cuttingPlaneMethod} never visits any point within a ball of radius $\frac{\epsilon}{L}$ (with respect to the Frobenius norm) twice. Moreover, by iteration $k$, Algorithm \ref{alg:cuttingPlaneMethod} visits $k$ points within non-overlapping balls with combined volume $$k \frac{\pi^\frac{n^2}{2}}{\Gamma(\frac{n^2}{2}+1)} \left(\frac{\epsilon}{L}\right)^{n^2}, $$
and these balls are centered at feasible points, i.e., contained within a ball of radius $K+\frac{\epsilon}{L}$ with volume $$ \frac{\pi^\frac{n^2}{2}}{\Gamma(\frac{n^2}{2}+1)} \left(K+\frac{\epsilon}{L}\right)^{n^2}. $$ That is, if Algorithm \ref{alg:cuttingPlaneMethod} has not converged at iteration $k$, we have:
$
    k < \left(\frac{L K}{\epsilon}+1\right)^{n^2},
$
which implies that we converge to an $\epsilon$-optimal solution within $k\leq \big(\frac{L K}{\epsilon}+1\big)^{n^2}$ iterations.
\hfill\Halmos \endproof

\subsection{Proof of {\color{black}Lemma} \ref{lemma:optygivenx}}\label{ssec:proofoflemmaaltmin1}
\proof{Proof of {\color{black}Lemma} \ref{lemma:optygivenx}}
The equality $\bm{\theta}^\star=\bm{X}_t^\top (\bm{Y}^\star)^\dag \bm{X}_t$ is immediate from the Schur complement lemma \ref{generalizedschur}. Therefore, we focus on deriving an optimal $\bm{Y}$ for a fixed $\bm{X}_t$.

From the second equality in the Schur complement lemma \ref{generalizedschur}, we must have $\bm{X}=\bm{Y}\bm{Y}^\dag \bm{X}$ for feasibility. Therefore, $\mathrm{span}(\bm{X}) \subseteq \mathrm{span}(\bm{Y})$. Moreover, columns of $\bm{Y}$ which are in $\mathrm{null}(\bm{X})$ do not contribute to the optimal objective and can therefore be omitted without loss of optimality. Therefore, $\bm{Y}=\sum_{i=1}^n \rho_i \bm{u}_i \bm{u}_i^\top$ for some $\bm{\rho}$, where $\bm{X}_t=\bm{U}\bm{\Sigma}\bm{V}^\top$ is an SVD of $\bm{X}_t$. The result follows from observing that $\bm{0} \leq \bm{\rho} \leq \bm{e}$ and $\bm{e}^\top \bm{\rho} \leq k$, since $\bm{Y} \in \mathrm{Conv}(\mathcal{Y}^k_n)$.
\hfill\Halmos \endproof

\subsection{Proof of {\color{black}Lemma} \ref{prop:conv_fullmax}}\label{sec:proofofpropositionconvminimax}
\proof{Proof of {\color{black}Lemma} \ref{prop:conv_fullmax}}
As Assumption \ref{strongduality} holds, we can exchange the minimization and maximization operators in Problem \eqref{minmaxconv_rankkcase}. Therefore, Problem \eqref{minmaxconv_rankkcase} has the same optimal objective as:
\begin{align}\label{conv_maxminrankk}
     \max_{\bm{\alpha}} h(\bm{\alpha})-\max_{\bm{Y} \in \mathrm{Conv}\left(\mathcal{Y}_n^k\right)}\frac{\gamma}{2}\sum_{i=1}^n \sum_{j=1}^n Y_{i,j}\langle \bm{\alpha}_i, \bm{\alpha}_j \rangle.
\end{align}
Therefore, to establish the result, it suffices to show that we obtain Problem \eqref{max_conv_rankk} after taking the dual of Problem \eqref{conv_maxminrankk}'s inner problem. This is indeed the case, because $\mathrm{Conv}(\mathcal{Y}_n^k)$ is a convex compact set with non-empty relative interior, and therefore strong duality holds between the following two problems:
\begin{align*}
    \max_{\bm{Y} \succeq \bm{0}} \quad & \frac{\gamma}{2}\sum_{i=1}^n \sum_{j=1}^n Y_{i,j}\langle \bm{\alpha}_i, \bm{\alpha}_j \rangle \quad & \text{s.t.} \quad & \bm{Y} \preceq \mathbb{I}, \ [\bm{U}]\  \langle \mathbb{I}, \bm{Y} \rangle \leq k, \ [t],\\
\min_{\bm{U} \succeq \bm{0}, t \geq 0} \quad & \mathrm{tr}(\bm{U})+kt\quad & \text{s.t.} \quad & \bm{U}+\mathbb{I}t \succeq \frac{\gamma}{2}\bm{\alpha}\bm{\alpha}^\top. \quad \hfill\Halmos
\end{align*}
\endproof

\subsection{Proof of Theorem \ref{thm:randrounding2}}\label{ssec:proofoftheoremrandround}
\proof{Proof of Theorem \ref{thm:randrounding2}}
{\color{black}To establish the result, we establish the first half of the following inequalites:}
\begin{align*}
    {\color{black}f(\bm{Y}_{rounded})}-f(\bm{Y}^\star) {\color{black}\leq \frac{\gamma}{2} L^2
    \max_{\bm{\beta} \geq 0 : \| \bm{\beta} \|_\infty \leq 1} \: \sum_{i \in \mathcal{R}} \left(\Lambda_{i,i}^\star-\Lambda_{i,i}^{rounded}\right) \beta_{i}} \leq \frac{\gamma}{2}L^2 \vert\mathcal{R}\vert \max_{\bm{\beta} \geq \bm{0}: \Vert \bm{\beta} \|_1 \leq 1} \sum_{i \in \mathcal{R}} (\Lambda_{i,i}^\star-{\color{black}\Lambda_{i,i}^{rounded}})\beta_i,
\end{align*}
under the Frobenius penalty and
\begin{align*}
   {\color{black}f(\bm{Y}_{rounded})}-f(\bm{Y}^\star)  {\color{black}\leq M L
    \max_{\bm{\beta} \geq 0 : \| \bm{\beta} \|_\infty \leq 1} \: \sum_{i \in \mathcal{R}} \left(\Lambda_{i,i}^\star-\Lambda_{i,i}^{rounded}\right) \beta_{i}} \leq M L \vert\mathcal{R}\vert \max_{\bm{\beta} \geq \bm{0}: \Vert \bm{\beta} \|_1 \leq 1} \sum_{i \in \mathcal{R}} (\Lambda_{i,i}^\star-{\color{black}\Lambda_{i,i}^{rounded}})\beta_i,
\end{align*}
for the spectral penalty—{\color{black} the second half of both inequalities follows readily from the fact that $\Vert \bm{\beta}\Vert_1 \leq \vert \mathcal{R}\vert \Vert \bm{\beta}\Vert_\infty \leq \vert \mathcal{R}\vert$ which allows us to replace $\Vert \bm{\beta}\Vert_\infty \leq 1$ with $\Vert \bm{\beta}\Vert_1 \leq \vert \mathcal{R}\vert$ and move $\vert \mathcal{R}\vert$ outside the bound}. {\color{black}Indeed, after establishing these inequalities, the result follows by observing that $\bm{Y}_{greedy}$ minimizes the right-hand-side of \eqref{eqn:concentrationbound1}-\eqref{eqn:eqn:concentrationbound2} over the class of projection matrices $\bm{Y}_{rounded}$.}

Under a Frobenius penalty, by Lipschitz continuity, we have
\begin{align*}
    {\color{black}f(\bm{Y}_{rounded})}-f(\bm{Y}^\star) \leq \frac{\gamma}{2}\langle \bm{\alpha}^\star(\bm{Y}) \bm{\alpha}^{\star}(\bm{Y})^\top, \bm{U}(\bm{\Lambda}^\star-{\color{black}\bm{\Lambda}_{rounded}})\bm{U}^\top\rangle=\frac{\gamma}{2}\langle \bm{U}^\top \bm{\alpha}^\star(\bm{Y}) \bm{\alpha}^{\star}(\bm{Y})^\top\bm{U}, \bm{\Lambda}^\star-{\color{black}\bm{\Lambda}_{rounded}}\rangle.
\end{align*}
Moreover, since $\bm{\Lambda}^\star-{\color{black}\bm{\Lambda}_{rounded}}$ is a diagonal matrix we need only include the diagonal terms in the inner product. Therefore, since $$(\bm{U}^\top \bm{\alpha}^\star(\bm{Y}) \bm{\alpha}^{\star}(\bm{Y})^\top\bm{U})_{i,i}=\langle \bm{\alpha}^{\star}(\bm{Y})^\top\bm{\alpha}^\star(\bm{Y}), \bm{U}_i\bm{U}_i^\top \rangle\leq \lambda_{\max}(\bm{\alpha}^{\star}(\bm{Y})^\top \bm{\alpha}^{\star}(\bm{Y})) \leq L^2,$$
where the second-to-last inequality holds because $\Vert \bm{U}_i\Vert_2=1$, the bound on ${\color{black}f(\bm{Y}_{rounded})}-f(\bm{Y}^\star)$ holds.

Alternatively, under spectral norm regularization, by Lipschitz continuity we have
\begin{align*}
    {\color{black}f(\bm{Y}_{rounded})}-f(\bm{Y}^\star) & \leq M \langle \bm{V}_{11}^\star(\bm{Y})+\bm{V}_{22}^\star(\bm{Y}), \bm{U}(\bm{\Lambda}^\star-{\color{black}\bm{\Lambda}_{rounded}})\bm{U}^\top\rangle\\
    & =M\langle \bm{U}^\top (\bm{V}_{11}^\star(\bm{Y})+\bm{V}_{22}^\star(\bm{Y}))\bm{U}, \bm{\Lambda}^\star-{\color{black}\bm{\Lambda}_{rounded}}\rangle.
\end{align*}
Moreover, $\bm{\Lambda}^\star-{\color{black}\bm{\Lambda}_{rounded}}$ is a diagonal matrix and therefore
$$
 (\bm{U}^\top (\bm{V}_{11}^\star(\bm{Y})+\bm{V}_{22}^\star(\bm{Y}))\bm{U})_{i,i}=\langle\bm{U}_i \bm{U}_i^\top, \bm{V}_{11}^\star(\bm{Y})+\bm{V}_{22}^\star(\bm{Y})\rangle\leq \lambda_{\max}(\bm{\alpha}^{\star}(\bm{Y})) \leq L,
$$
where the last inequality follows since $\bm{V}_{11},\bm{V}_{22}$ are orthogonal at optimality, meaning $\bm{V}_{11}+\bm{V}_{22}$'s leading eigenvalue equals $\bm{\alpha}^\star$'s leading singular value. Therefore, the bound on ${\color{black}f(\bm{Y}_{rounded})}-f(\bm{Y}^\star)$ holds. \hfill\Halmos
\endproof

{\color{black}
\section{Derivations for the conjugate of the regularizer}\label{sec:A.linear.conjugate}}
In this section, we derive the conjugates of the penalties stated in Table \ref{tab:reg_conj}, in order to complete our proof of Lemma \ref{lemma:conj.reg}. We first derive our results for rectangular matrices under the formulation $\bm{X}=\bm{Y}\bm{X}\bm{Y}'$ for appropriate projection matrices $\bm{Y}, \bm{Y}'$, before specializing our results by setting $\bm{Y}'=\mathbb{I}$. For completeness, we first prove that this reformulation is indeed a valid reformulation of a rank constraint.

\begin{proposition}\label{prop:rectangularprojreformulation}
For any $\bm{X} \in \mathbb{R}^{n \times m},\ \text{Rank}( \bm{X} ) \leq k \iff
  \exists \bm{Y} \in\mathcal{Y}_n, \bm{Y}' \in \mathcal{Y}_m \ \mbox{ s.t. } \ \text{tr}(\bm{Y}),\text{tr}(\bm{Y}') \leq k \mbox{ and } {\bm{X} = \bm{Y}\bm{X}\bm{Y}'},$ $\text{where} \ \mathcal{Y}_n:=\{\bm{P} \in S^n: \bm{P}^2 =\bm{P}\} \ \text{is the set of projection matrices.}$
\end{proposition}

\proof{Proof of Proposition \ref{prop:rectangularprojreformulation}} We prove the two implications successively.
\begin{itemize}
    \item[$(\implies)$] Let $\bm{X}=\bm{U}\bm{\Sigma}\bm{V}^\top$ be a singular value decomposition of $\bm{X}$. Since $\mathrm{Rank}(\bm{X}) \leq k$, we can let $\bm{U} \in \mathbb{R}^{n \times k}, \bm{\Sigma} \in \mathbb{R}^{k \times k}, \bm{V} \in \mathbb{R}^{k \times m}$ without loss of generality. Define $\bm{Y}=\bm{U}(\bm{U}^\top \bm{U})^{-1}\bm{U}$ and $\bm{Y}'=\bm{V}(\bm{V}^\top \bm{V})^{-1}\bm{V}$.  By construction, $\bm{Y} \bm{X}\bm{Y}' = \bm{X}$. In addition, $\text{tr}(\bm{Y}) = \text{rank}(\bm{Y}) =  \text{rank}((\bm{U})) \leq k$, and $\text{tr}(\bm{Y}') = \text{rank}(\bm{Y}') =  \text{rank}((\bm{V})) \leq k$.
    \item[$(\impliedby)$] Since $\bm{X} = \bm{Y} \bm{X}\bm{Y}'$,
 $\text{rank}(\bm{X}) \leq \min(\text{rank}(\bm{Y}), \text{rank}(\bm{Y}')) \leq \min(\text{tr}(\bm{Y}), \text{tr}(\bm{Y}'))\leq k$.\hfill\Halmos
 \end{itemize}
\endproof

\subsection{Rectangular Case}
In this section, we derive a dual reformulation for the conjugate of the regularization term in \eqref{rankminproblem_proj}. More precisely, for all regularizers of interest $\Omega(\cdot)$, we show that for any matrix $\bm{\alpha}$ of similar dimension as $\bm{X}$ and any projection matrices $\bm{Y} \in \mathcal{Y}_n^k$ and $\bm{Y}' \in \mathcal{Y}_m^k$,
\begin{align*}
  \min_{\bm{X}} \left[ \Omega(\bm{Y}\bm{X}\bm{Y}') + \langle \bm{\alpha}, \bm{Y}\bm{X}\bm{Y}' \rangle \right] = \max_{\bm{V}_{11}, \bm{V}_{22}} \: - \Omega^\star(\bm{\alpha}, \bm{Y}, \bm{Y}', \bm{V}_{11}, \bm{V}_22),
\end{align*}
where $\Omega^\star(\cdot)$ is notably linear in $\bm{Y}, \bm{Y}'$.

\begin{appendixlemma} \label{lemma:conj_frob_rect}
Let $\bm{A} \in \mathbb{R}^{n \times m}$ be a rectangular matrix, $\bm{Y} \in \mathcal{Y}_n$, $\bm{Y}' \in \mathcal{Y}_m$ be projection matrices and $\gamma > 0$ be a positive scalar. Then
\begin{align*}
\min_{\bm{X} \in \mathbb{R}^{n \times m} } \: \dfrac{1}{2 \gamma} \| \bm{Y} \bm{X} \bm{Y}' \|_F^2 + \langle \bm{A}, \bm{Y} \bm{X} \bm{Y}' \rangle =   -\dfrac{\gamma}{2}  \mathrm{tr}(\bm{Y} \bm{A} \bm{Y}' \bm{A}^\top).
\end{align*}
\end{appendixlemma}
\proof{Proof of Lemma \ref{lemma:conj_frob_rect}} Any solution to the minimization problem satisfies the first-order condition $\tfrac{1}{\gamma} \bm{Y} \bm{X} \bm{Y}' + \bm{Y} \bm{A} \bm{Y}' = 0$.
Hence, $\bm{X}^\star = - \gamma \bm{Y} \bm{A} \bm{Y}'$ is optimal and the objective is $-\tfrac{\gamma}{2}  \mathrm{tr}(\bm{Y} \bm{A} \bm{Y}' \bm{A}^\top)$.
\hfill\Halmos \endproof
\begin{appendixlemma} \label{lemma:conj_spectral_rect}
Let $\bm{A} \in \mathbb{R}^{n \times m}$ be a rectangular matrix, $\bm{Y} \in \mathcal{Y}_n$, $\bm{Y}' \in \mathcal{Y}_m$ be projection matrices and $M > 0$ be a positive scalar. Then
\begin{align*}
\min_{\bm{X} \in \mathbb{R}^{n \times m} : \| \bm{Y} \bm{X} \bm{Y}' \|_\sigma \leq M } \: \langle \bm{A}, \bm{Y} \bm{X} \bm{Y}' \rangle =     \max_{\bm{V}_{11}, \bm{V}_{22}} \quad & -\tfrac{M}{2} \langle \bm{V}_{11}, \bm{Y} \rangle - \tfrac{M}{2} \langle \bm{V}_{22}, \bm{Y}'\rangle\\
    \text{s.t.} \quad & \begin{pmatrix} \bm{V}_{11} & \bm{A} \\ \bm{A}^\top & \bm{V}_{22}\end{pmatrix} \succeq \bm{0},
\end{align*}
\end{appendixlemma}
\proof{Proof of Lemma \ref{lemma:conj_spectral_rect}} Since the trace is invariant by cyclic permutation and the matrices $\bm{Y}$ and $\bm{Y}'$ are symmetric,
\begin{align*}
\min_{\bm{X} \in \mathbb{R}^{n \times m} : \| \bm{Y} \bm{X} \bm{Y}' \|_\sigma \leq M } \: \langle \bm{A}, \bm{Y} \bm{X} \bm{Y}' \rangle & = \quad
\min_{\bm{X} \in \mathbb{R}^{n \times m} : \| \bm{Y} \bm{X} \bm{Y}' \|_\sigma \leq M } \: \langle \bm{Y} \bm{A} \bm{Y}', \bm{X} \rangle.
\end{align*}
The spectral norm penalty is equivalent to
\begin{align*}
\| \bm{Y} \bm{X} \bm{Y}' \|_\sigma \leq M & \iff \bm{Y} \bm{X} \bm{Y}' \bm{Y}' \bm{X}^\top \bm{Y} \preceq M^2 \bm{I}_n \\
&\iff (\bm{Y} \bm{X} \bm{Y}') \bm{Y}' (\bm{Y}' \bm{X}^\top \bm{Y}) \preceq M^2 \bm{I}_n \\
&\iff (\bm{Y} \bm{X} \bm{Y}') \bm{Y}' (\bm{Y}' \bm{X}^\top \bm{Y}) \preceq M^2 \bm{Y},
\end{align*}
where the last inequality follows from the fact that the matrix on the left-hand side is equal to zero over $\text{Im}(\bm{Y})^\top$. By Schur complements, the final semidefinite inequality is equivalent to
\begin{align*}
\begin{pmatrix}
  M \bm{Y} & \bm{Y} \bm{X} \bm{Y}' \\ \bm{Y}' \bm{X}^\top \bm{Y} & M \bm{Y}'
\end{pmatrix} \succeq \bm{0} \quad & \left[\begin{pmatrix}
  \bm{W}_{11} & \bm{W}_{12}\\ \bm{W}_{12}^\top & \bm{W}_{22}
\end{pmatrix} \right].
\end{align*}
We associate a matrix of dual variables in square brackets. Therefore, the dual problem is
\begin{align*}
  \max_{\bm{W}_{11}, \bm{W}_{12}, \bm{W}_{22}} \quad & - M \langle \bm{W}_{11}, \bm{Y} \rangle -M \langle  \bm{W}_{22}, \bm{Y}' \rangle\\
    \text{s.t.} \quad & \bm{Y} \bm{A}\bm{Y}' = 2  \bm{Y} \bm{W}_{12} \bm{Y}',\ \begin{pmatrix}
  \bm{W}_{11} & \bm{W}_{12}\\ \bm{W}_{12}^\top & \bm{W}_{22}
\end{pmatrix} \succeq \bm{0},\\
  \text{or equivalently,} \quad \max_{\bm{V}_{11}, \bm{V}_{22}} \quad & -\tfrac{M}{2} \langle \bm{V}_{11}, \bm{Y} \rangle - \tfrac{M}{2} \langle \bm{V}_{22}, \bm{Y}'\rangle\
    \text{s.t.} \quad \begin{pmatrix} \bm{V}_{11} & \bm{A} \\ \bm{A}^\top & \bm{V}_{22}\end{pmatrix} \succeq \bm{0}. \qquad \hfill\Halmos
\end{align*}
\endproof


\subsection{Square and Symmetric Case}
When the matrix $\bm{X}$ is square and symmetric we can take $\bm{Y} = \bm{Y'}$ and apply the previous results.
Alternatively, for the spectral norm penalty, we can further simplify
\begin{appendixlemma}
\label{lemma:conj_spectral_symmetric}
Let $\bm{A} \in \mathbb{R}^{n \times n}$ be a square symmetric matrix and $\bm{Y} \in \mathcal{Y}_n$ be a projection matrix. Then
\begin{align*}
\min_{\bm{X} \in \mathbb{S}_n : \| \bm{Y} \bm{X} \bm{Y} \|_\sigma \leq M } \: \langle \bm{A}, \bm{Y} \bm{X} \bm{Y} \rangle =     \max_{\bm{V}_{+}, \bm{V}_{-} \succeq 0} \quad & - M \langle \bm{V}_{+} + \bm{V}_{-}, \bm{Y} \rangle\quad \text{s.t.} \quad \bm{A} = \bm{V}_{+} - \bm{V}_{-}.
\end{align*}
\end{appendixlemma}
\proof{Proof of Lemma \ref{lemma:conj_spectral_symmetric}}
The constraint $\| \bm{Y} \bm{X} \bm{Y} \|_\sigma \leq M$ can be rewritten as $- M \bm{Y} \preceq \bm{Y} \bm{X} \bm{Y}  \preceq M \bm{Y}$. By strong semidefinite duality (which holds as the minimization problem has non-empty interior; see \citep[][Chapter 4.1]{wolkowicz2012handbook}):
\begin{align*}
\min_{\bm{X} \in \mathbb{S}_n : \| \bm{Y} \bm{X} \bm{Y}  \|_\sigma \leq M } \: \langle \bm{A}, \bm{Y} \bm{X} \bm{Y} \rangle = \max_{\bm{W}_{+}, \bm{W}_{-} \succeq 0} - M \langle \bm{Y}, \bm{W}_{+}+\bm{W}_{-}\rangle \quad \text{s.t.}\quad \bm{Y}\bm{A} \bm{Y} =  \bm{Y}\bm{W}_{+} \bm{Y}- \bm{Y}\bm{W}_{-} \bm{Y}.
\end{align*}
The decision variables in the maximization problem decompose the symmetric matrix $\bm{Y} \bm{A} \bm{Y}$ into a positive and negative definite parts, $\bm{W}_+$ and $\bm{W}_-$ respectively. Without loss of optimality we can restrict our attention to $\bm{W}_+ = \bm{Y} \bm{V}_+  \bm{Y}$ and $\bm{W}_- = \bm{Y} \bm{V}_-  \bm{Y}$ where $\bm{A} = \bm{V}_+ - \bm{V}_- $.
\hfill\Halmos \endproof

\section{A Collection of Useful Matrix Identities}\label{append:auxresults}
In this work, we have repeatedly used some technical matrix identities. In order to keep this paper self contained, we now state these identities formally.


The following result generalizes the well-known Schur complement lemma to the case where neither on-diagonal block matrix is positive definite \citep[see][Equation 2.41]{boyd1994linear}
\begin{appendixlemma}[Generalized Schur Complement Lemma]\label{generalizedschur}
Let $\bm{A}, \bm{B}, \bm{C}$ be components of $$\bm{X}:=\begin{pmatrix} \bm{A} & \bm{B}\\ \bm{B}^\top &\bm{C}\end{pmatrix}$$ of appropriate dimension. Then, $\bm{X}$ is positive semidefinite if and only if the following conditions hold:
\begin{itemize}
    \item[(i)] $\bm{A} \succeq \bm{0}$,
    \item[(ii)] $(\mathbb{I}-\bm{A}\bm{A}^\dag)\bm{B}=\bm{0}$,
    \item[(iii)] $\bm{C}\succeq \bm{B}^\top \bm{A}^\dag \bm{B}$.
\end{itemize}
\end{appendixlemma}

The following result characterizes the rank of a block submatrix in terms of the rank of the entire matrix \citep[see, e.g.,][Section 0.9]{zhang2006schur}:
\begin{appendixlemma}[Guttman Rank Identity]\label{guttmanrankidentity}
Let $\bm{A}, \bm{B}, \bm{C}$ be components of the matrix $$\bm{X}:=\begin{pmatrix} \bm{A} & \bm{B}\\ \bm{B}^\top &\bm{C}\end{pmatrix}$$ of appropriate dimension. Then, we have the identity:
$
    \mathrm{Rank}(\bm{X})=\mathrm{Rank}(\bm{A})+\mathrm{Rank}(\bm{C}-\bm{B}\bm{C}^\dag \bm{B}^\top).
$
\end{appendixlemma}

In general, a product of positive semidefinite matrices may not be positive semidefinite, indeed, it may not even be symmetric. However, the following result demonstrates that a \textit{symmetric} product of PSD matrices is indeed PSD \citep[see][for a proof]{meenakshi1999product}:
\begin{appendixlemma}[A Symmetric Product of PSD Matrices is PSD]\label{symmproductispsd}
Let $\bm{X}, \bm{Y} \in S^n_+$ be matrices of appropriate size, and let $\bm{Z}:=\bm{X}\bm{Y}$. Suppose that $\bm{Z}=\bm{Z}^\top$ is a symmetric matrix. Then, $\bm{Z} \succeq \bm{0}$.
\end{appendixlemma}

The following result demonstrates that evaluating the nuclear norm of a matrix via semidefinite optimization yields a singular value decomposition \citep[see][Proposition 2.1, for a proof]{recht2010guaranteed}
\begin{appendixlemma}[Nuclear norm minimization and the SVD of a matrix]\label{lemma:svdnuclearnorm}
An optimal solution to
\begin{align*}
    \min_{\bm{W}_1 \in \mathbb{R}^{n \times n}, \bm{W}_2 \in \mathbb{R}^{m \times m}} \quad & \mathrm{tr}(\bm{W}_1)+\mathrm{tr}(\bm{W}_2)\quad \text{\rm s.t.} \quad \begin{pmatrix} \bm{W}_1 & \bm{X} \\
    \bm{X}^\top & \bm{W}_2 \end{pmatrix} \succeq \bm{0},
\end{align*}
is attained by $\bm{W}_1=\bm{U}\bm{\Sigma}\bm{U}^\top$ and $\bm{W}_2=\bm{V}\bm{\Sigma}\bm{V}^\top$, where $\bm{X}=\bm{U}\bm{\Sigma}\bm{V}^\top$ is a singular value decomposition.
\end{appendixlemma}

\section{Additional Results}
\subsection{Pseudocode for the In-out Method}\label{ssec:inout}
Our main points of difference from \citet{fischetti2016redesigning}'s implementation are twofold. First, we optimize the outer problem over $\mathrm{Conv}(\mathcal{Y}_n^k)$, rather than the Boolean polytope. Second, we recognize that the purpose of the method is to warm-start Algorithm \ref{alg:cuttingPlaneMethod}'s lower bound, rather than to solve Problem \eqref{prob:sdprelax}. In this spirit, we accelerate the \verb|in-out| method by first solving Problem \eqref{prob:sdprelax} in one shot using an interior point method and second using the optimal solution $\bm{Y}^\star$ as a stabilization point. Note that a similar method was proposed for sparse portfolio selection problems and MICO by \citet{bertsimas2018scalable}.
\begin{algorithm*}[h!]
\caption{The in-out method of \citet{ben2007acceleration}.}
\label{alg:CuttingPlaneMethodinout}
\begin{algorithmic}
\REQUIRE Stabilization point $\bm{Y}^\star$, $\epsilon \gets 10^{-10}, \lambda \gets 0.1, \delta \gets 2 \epsilon$, $t \leftarrow 1 $
\REPEAT
\STATE Compute $\bm{Y}_0, \theta_0$ solution of \vspace{-2mm}
\begin{align*}
\min_{\bm{Y} \in \mathrm{Conv}\left(\mathcal{Y}_n^k\right), \theta} \: \theta+\lambda \cdot \mathrm{tr}(\bm{Y}) \quad \mbox{ s.t. } \quad z_i\theta \geq h_i + \left\langle H_i, \bm{Y}-\bm{Y}_i\right\rangle \ \forall i \in [t].
\end{align*}\vspace{-5mm}
\IF{$\bm{Y}_0$ has not improved for $5$ consecutive iterations}
\STATE Set $\lambda=1$
\IF{$\bm{Y}_0$ has not improved for $10$ consecutive iterations}
\STATE Set $\delta=0$
\ENDIF
\ENDIF
\STATE Set $\bm{Y}_{t+1} \gets \lambda \bm{Y}_0+(1-\lambda)\bm{Y}^\star+\delta \mathbb{I}$; project $\bm{Y}_{t+1}$ onto $\mathrm{Conv}(\mathcal{Y}_n^k)$.
\STATE Compute $f(\bm{Y}_{t+1}), H_{t+1}, z_{t+1}, d_{t+1}$.
\UNTIL{$ f(\bm{Y}_0)-\theta_0\leq \varepsilon$}
\RETURN $\bm{Y}_t$
\end{algorithmic}
\end{algorithm*}
\FloatBarrier
{\color{black}
\subsection{Number of cuts generated vs. $\gamma$}
\begin{figure}[h!]

    \begin{subfigure}[t]{.45\linewidth}
    \centering
            \includegraphics[scale=0.4]{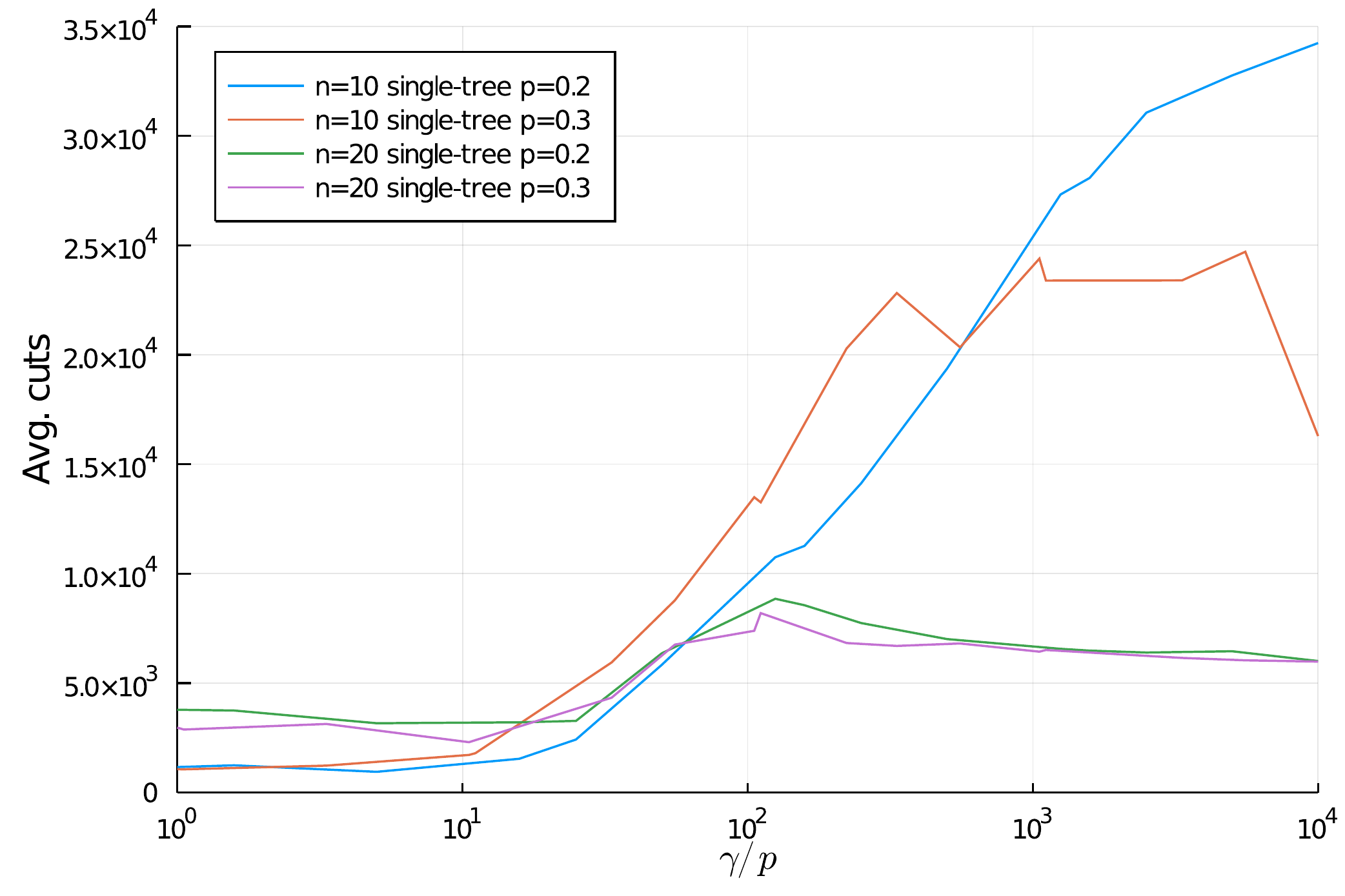}
    \end{subfigure}
        \begin{subfigure}[t]{.45\linewidth}
            \includegraphics[scale=0.4]{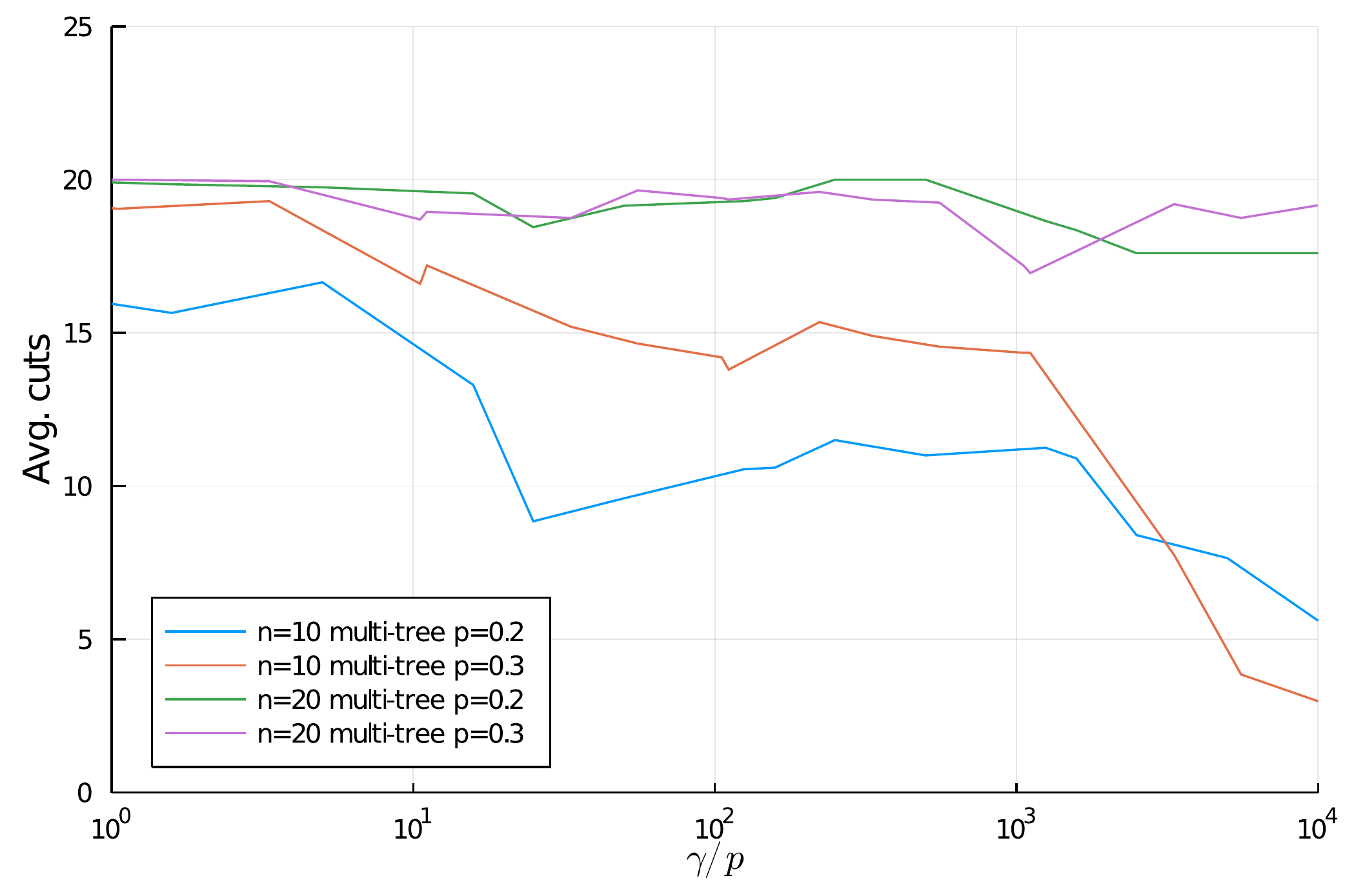}
    \end{subfigure}
   \caption{\color{black}Average number of cuts vs. $\gamma$ for Algorithm \ref{alg:cuttingPlaneMethod} single-tree (left) and multi-tree (right) implementations over $20$ synthetic matrix completion instances where $p \in \{0.2, 0.3\}$, $r=1$ and $n \in \{10, 20\}$. The same random seeds were used to generate random matrices completed by single-tree and multi-tree.}
   \label{fig:sensitivitytogamma2}
\end{figure}

}

\end{appendices}
\end{document}